\documentclass[11pt]{amsart}
\allowdisplaybreaks[4]
\usepackage{amssymb}
\usepackage{amsfonts}
\usepackage{graphicx}
\usepackage{epstopdf}
\usepackage{dcolumn}
\usepackage{amsmath}
\usepackage{enumerate}
\usepackage{latexsym,bm}
\usepackage{slashed}
\usepackage{float}
\usepackage{cite}
\usepackage{CJK}
\usepackage[all,cmtip]{xy}
\usepackage{appendix}
\usepackage{ulem}
\usepackage[colorlinks,
linkcolor=black,
citecolor=red
]{hyperref}

\DeclareMathOperator{\Span}{Span}

\DeclareMathOperator{\Supp}{Supp}

\DeclareMathOperator{\rr}{r}
\DeclareMathOperator{\sr}{sr}
\DeclareMathOperator{\id}{Id}

\DeclareMathOperator{\im}{im}

\DeclareMathOperator{\sign}{sign}
\newcommand{\A}{\mathcal A}
\newcommand{\Z}{\mathbb Z}

\newcommand{\K}{\mathbb K}

\theoremstyle{plain}
\newtheorem{theorem}{Theorem}[section]
\newtheorem{lemma}[theorem]{Lemma}
\newtheorem{lem/def}[theorem]{Lemma/Definition}
\newtheorem{cor/def}[theorem]{Corollary/Definition}
\newtheorem{proposition}[theorem]{Proposition}
\newtheorem{corollary}[theorem]{Corollary}
\newtheorem{claim}[theorem]{Claim}

\theoremstyle{definition}
\newtheorem{definition}[theorem]{Definition}
\newtheorem{example}[theorem]{Example}

\newtheorem{remark}[theorem]{Remark}
\numberwithin{equation}{section}

\newenvironment{acknowledgement}{\smallskip{\sc Acknowledgments.}\rm}{\smallskip}

\allowdisplaybreaks

\begin{document}

\title{Minimal Path and Acyclic Model in the Path Complex}
\author{Xinxing Tang and Shing-Tung Yau}

  \address{
X. Tang: Beijing Institute of Mathematical Sciences and Applications, Beijing, China;
}
\email{tangxinxing@bimsa.cn}
  \address{
S.-T. Yau: Yau Mathematical Sciences Center, Tsinghua University and Beijing Institute of Mathematical Sciences and Applications, Beijing, China;\\
Beijing Institute of Mathematical Sciences and Applications, Beijing, China;
}
\email{styau@tsinghua.edu.cn}

\maketitle

\begin{abstract} In this paper, firstly, we will study the structure of the path complex $(\Omega_*(G;\Z),\partial)$ of a digraph $G$ via the $\Z$-generators of $\Omega_*(G,\Z)$  under strongly regular condition, which is called the minimal path in \cite{HY}. In particular, we will study various examples of the minimal $3$-paths. Secondly, we will show that the supporting sub-digraph of minimal path has acyclic path homologies. Thirdly, we will consider the applications of such an acyclic model.
\end{abstract}

\tableofcontents

\section{Introduction}

In the past few years, there are several attempts to define the homotopy and (co)homology of (di)graphs, e.g. via cliques \cite{CYY}, or via Hochschild (co)homology \cite{Hochschild}. The homotopy and (co)homology theory of digraphs considered in this paper were introduced in \cite{GJMY,GLMY3,GLMY,GLMY2,GMY1,GMY}.
Such a homology-homotopy theory for digraphs is now known as GLMY theory. It has many advantages in comparison with the previously studied notions of graph homologies:
\begin{enumerate}
  \item The path complex could be regarded as a generalization of the simplicial complex. It contains not
only cliques but also binary hypercubes and many other subgraphs. (One can see more examples in Subsection \ref{exminmal3}.)
  \item It satisfies the properties that are analogous to Eilenberg-Steenrod axioms.
  \item It is well linked to graph-theoretical operations. For example, the K\"{u}nneth formula holds for both Cartesian product of two digraphs and join of two digraphs.
  \item There is a dual cohomology theory of digraphs which coincides with the one developed by Dimakis and M\"{u}ller-Hoissen in \cite{DM1,DM} as a reduction of the ``universal differential algebra".
\end{enumerate}

\subsubsection{Motivations}
There are rapidly increasing developments of GLMY theory, such as its relation with simplicial complex and cubical complex \cite{GMY1} under some conditions, homotopy theory of digraph \cite{GLMY}, the discrete Morse theory on digraph \cite{LWY}. Recently, Ivanov-Pavutnitskiy \cite{IP} develop a generalisation of the path homology theory in a general simplicial setting. However, there still remain many basic and important questions. For example, the concatenation of paths induces a well defined cup product on path cohomologies. But it is still an open problem whether the cup product is skew-symmetric, since one can not arbitrarily change the order of the path vertices (due to the restriction of the directions of edges). To give some positive answers for it becomes one of our motivations of this paper.

One can check that the skew-symmetry holds in some simple examples by following the definitions, that is,
\begin{itemize}
  \item Compute the path cohomologies of digraph $G$ first;
  \item Choose representatives $\varphi$, $\psi$, such that $[\varphi]\in H^p(G)$ and $[\psi]\in H^q(G)$. Check that $\varphi\cup \psi$ and $(-1)^{pq}\psi\cup\varphi$ are cohomologously equivalent.
\end{itemize}

We find that both steps are complicated (even in the case that $G$ is the Cartesian product of monotone cycles) and uninspiring. Then we switch to a traditional but also more powerful method, which is based on the acyclic model. Thus, a natural question arises: what is the acyclic model in the path complex of a digraph? Under the strongly regular condition, the answer is given by the minimal path (see, Definition \ref{minimaldef}) and its supporting digraph (see Definition \ref{supp}).

The new question forces us to have a better understanding of the structure of the path complex.

\subsubsection{Main theorems}
The definition of the minimal path here was introduced in Huang-Yau's paper\cite{HY}, where the minimal path and minimal relation played an important role to obtain an integral basis of the path complex. In this paper we reprove and improve their result by deeply studying the definition of path complex.

During the research, we first study many minimal $3$-paths. They have a rich combinatorial structure which we summarize and generalize to the minimal paths of arbitrary length. We prove the following theorem by double/lexicographic induction first on the length of minimal path and then on the number of points in $E_1$ (see the definition in Subsection \ref{defbasicpropsubsection}). Our first main theorem is as follows.

\begin{theorem}[Theorem \ref{structurethm}] Let $P\in\Omega_n(G;\Z)$ be a minimal path with the starting vertex $S$ and ending vertex $E$, $\Supp(P)$ be its supporting digraph and $d_S$, $d_E$ be two distance functions (see the definition in Subsection \ref{defbasicpropsubsection}).

(1) Let $S_1=d_S^{-1}(1)$ and $E_1=d_E^{-1}(1)$. Then
$P$ is a linear combination of s-regular allowed elementary paths from $S$ to $E$, with coefficients being either 1 or -1. And in $\Supp(P)$,
\begin{equation*}
\begin{aligned}
\partial P=\sum_{\alpha\in E_1}P_{S,n-1,\alpha}+\sum_{\beta\in S_1}P_{\beta,n-1,E}+\sum_{k\in I_P} P_{S,n-1,E}^k,
\end{aligned}
\end{equation*}
where
\begin{itemize}
  \item  $P_{S,n-1,\alpha},P_{\beta,n-1,E},P_{S,n-1,E}^k\in\Omega_{n-1}(\Supp(P);\Z)$, and moreover
  \begin{itemize}
          \item $P_{S,n-1,\alpha}$ is the minimal $(n-1)$-path starting from $S$ and ending with $\alpha$;
          \item $P_{\beta,n-1,E}$ is the minimal $(n-1)$-path starting from $\beta$ and ending with $E$;
          \item $P_{S,n-1,E}^k$ is the minimal $(n-1)$-path starting from $S$ and ending with $E$.
        \end{itemize}
  \item Such $P_{S,n-1,\alpha}$, $P_{\beta,n-1,E}$ are unique (up to a sign) in $\Supp(P)$ for each $\alpha\in E_1$, $\beta\in S_1$.
  \item The set $I_P$ in the last summand depends on $P$, and $|I_P|\leq 1$. That is, there exists at most one (up to a sign) minimal $(n-1)$-path in $\Supp(P)$ starting from $S$ and ending with $E$.
\end{itemize}

(2) For any $v\in d_E^{-1}(k)\cap\Supp(P)$, in $\Supp(P)$,
\begin{itemize}
  \item There is a unique minimal $(n-k)$-path (up to a sign) starting from $S$ and ending with $v$, denoted by $P_{S,n-k,v}$, as well as a unique minimal $k$-path (up to a sign) starting from $v$ and ending with $E$, denoted by $P_{v,k,E}$.
  \item There is at most one minimal $(n-k-1)$-path (up to a sign) starting from $S$ and ending with $v$, denoted by $P_{S,n-k-1,v}$, as well as at most one minimal $(k-1)$-path starting from $v$ and ending with $E$, denoted by $P_{v,k-1,E}$.
\end{itemize}

(3) Each minimal $2$-path with fixed starting and ending vertices in $\Supp(P)$ ($n\geq2$) is unique.
\end{theorem}

We also study the homotopic and homological properties of several examples of minimal $3$-paths. We find that some of their supporting digraphs are contractible in the sense of digraph homotopy (see definition in Subsubsection \ref{Homotopy}), while some are not. What is more important to us is that their supporting digraphs have acyclic path homologies. A general acyclic result holds for the supporting digraph of any minimal path, which is our second main theorem.

\begin{theorem}[Theorem \ref{acyclicresult}] Let $P\in\Omega_n(G;\Z)$ be a minimal path, and $\Supp(P)$ be its supporting digraph, then
$$H_i(\Supp(P);\Z)=0,~i>0;\quad H_0(\Supp(P);\Z)=\Z.$$
\end{theorem}

The proof is also based on the lexicographic induction on the length of the minimal path and the number of points in $E_1$. To deal with the induction, we study some operations on the supporting digraph $\Supp(P)$. The idea is as follows:
\begin{itemize}
  \item Firstly, we add new edges to obtain a larger digraph $\widehat{\Supp(P)}$, and split $\widehat{\Supp(P)}$ into two smaller digraphs with decreasing $|E_1|$.
  \item Secondly, we apply Mayer-Vietoris method on the above decomposition for $\widehat{\Supp(P)}$.
  \item Thirdly, we compare the path homologies of $\Supp(P)$ and $\widehat{\Supp(P)}$.
\end{itemize}

It follows immediately from the acyclic model theorem, we obtain

\begin{theorem}[Theorem \ref{skewsymm}] For $\varphi\in H^p(G)$, $\psi\in H^q(G)$ defined under the strongly regular condition, we have
$$\varphi\cup\psi=(-1)^{pq}\psi\cup\varphi.$$
\end{theorem}

\subsubsection{Organizations}
Our paper is organized as follows.
In Section \ref{pathcomplexsection}, we will briefly recall the definition of path complex of digraphs, where we restricted to a strongly regular condition.

In Subsection \ref{defbasicpropsubsection}, we introduce the definitions of minimal path and its supporting digraph. And then we make an intensive study of its structure step by step. In Subsections \ref{strthmsubsection} and \ref{strproofsubsection}, we arrive at our structure theorem and its proof. In Subsection \ref{exminmal3}, we will give various examples of minimal $3$-paths and study their homotopic and homological results, where the homotopy theory of digraphs will be briefly recalled in advance.

In Section \ref{acyclicsection}, we state and prove our second main theorem. In Subsection \ref{MVsubsection}, we will recall the Mayer-Vietoris exact sequence in the path complex which is one of our main tools in the proof. The whole proof of the acyclic result will be shown in Subsection \ref{Proof}.

In Section \ref{applicationsection}, we give the applications of the acyclic model. In Subsection \ref{cohomology}, we first recall the definitions of path cohomology of digraphs. In Subsection \ref{cup}, we reformulate the cup product through the diagonal approximation map, and then apply the acyclic result to prove the skew-symmetry property of cup product on the path cohomology.

In Appendix \ref{A}, we give an example that the path homologies defined under the regular condition and the strongly regular condition respectively are not the same.

In Appendix \ref{B}, we give an example of minimal $4$-path and explain the proof idea of our second main theorem via such an example.

\begin{acknowledgement} The authors thank Alexander Grigoryan, Sergei Ivanov, Yifan Li and Xiaoyu Su for useful discussions. The work of X.T. is partially supported by research funding at BIMSA.
\end{acknowledgement}
\vskip 0.2cm

\section{Path Complex of Digraphs}\label{pathcomplexsection}

Let $V$ be a finite set. For any $n\geq0$, an elementary $n$-path is any ordered sequence $i_0,\ldots, i_n$ of $n+1$ vertices of $V$, write it as $e_{i_0\ldots i_n}$. Let $\Lambda_n(V;\Z)$ be the $\Z$-module generated by all such elementary $n$-paths. We also call $n$ the length of the path in $\Lambda_n(V;\Z)$.

The $\Z$-homomorphism $\partial:\Lambda_n(V;\Z)\rightarrow \Lambda_{n-1}(V;\Z)$ is defined via the generators:
$$\partial e_{i_0\ldots i_n}=\sum_{j=0}^n(-1)^je_{i_0\ldots \widehat{i_j}\ldots i_n}.$$
Clearly, $\partial\left(\Lambda_n(V;\Z)\right)\subset\Lambda_{n-1}(V;\Z)$, and $\partial^2=0$.

\begin{definition}\label{strongregular}
Let $G=(V,E)$ be a finite simple\footnote{Here simple condition means that there are no multiple edges and self loops.} digraph, where $V$ is the set of vertices and $E$ is the set of directed edges.
\begin{enumerate}
  \item An elementary path $e_{i_0i_1\ldots i_n}$ is called \textbf{allowed}, if each directed edge $i_{k-1}\rightarrow i_k$ belongs to $E$, $k=1,2\ldots, n$.
  \item An elementary path $e_{i_0i_1\ldots i_n}$ is called \textbf{strongly regular} (abbreviated as \textbf{s-regular}), if all $i_0, i_1,\ldots,i_n$ are distinct.
\end{enumerate}
\end{definition}

Let $\mathcal{A}_n(G;\Z)$ be the free $\Z$-module generated by all s-regular allowed elementary $n$-paths. Note that, $\mathcal{A}_n(G;\Z)=0$, if $n\geq\big|V(G)\big|$.

If we take the boundary operator $\partial$ on the s-regular allowed elementary path $e_{i_0\ldots i_n}$, in general, $\partial e_{i_0\ldots i_n}$ may not be in $\mathcal{A}_{n-1}(G;\Z)$. Then furthermore, we consider the submodule $\Omega_n(G;\Z)$ of $\partial$-invariant s-regular allowed $n$-paths.
That is
$$\Omega_n(G;\Z)=\{p\in \mathcal{A}_{n}(G;\Z)| \partial p\in \mathcal{A}_{n-1}(G;\Z)\}.$$

\begin{definition}\label{pathcomplex} The complex $(\Omega_*(G;\Z), \partial)$ is called the path complex of a digraph $G$ with $\Z$-coefficients (under the s-regular condition). For each $n\geq0$,
$$H_n(G,\Z):=\frac{\ker(\partial:\Omega_n(G;\Z)\rightarrow \Omega_{n-1}(G;\Z))}{\im(\partial:\Omega_{n+1}(G;\Z)\rightarrow \Omega_{n}(G;\Z))},$$
it is called the $n$-th path homology of $G$.
\end{definition}

\begin{remark}\label{regular} The s-regular condition is more restrictive than the one in \cite{GLMY3}.
In \cite{GLMY3}, the set of regular $n$-paths is defined as a quotient space $\mathcal{R}_n(V)$:
$$\mathcal{R}_n(V)=\Lambda_n(V)/I_n(V),$$
where $I_n(V)$ is the sub-module generated by the irregular path $e_{i_0i_1\ldots i_n}$, where $i_{k-1}=i_k$ for some $k=1,\ldots,n$. One can easily check that $(\mathcal{R}_*(V),\partial)$ forms a quotient chain complex with the induced boundary operator, we still denote by $\partial$.

The authors of \cite{GLMY3} consider the submodule of allowed paths $\mathcal{A}_n(G)\subset\mathcal{R}_n(G)$ and $\partial$-invariant submodule $\Omega_n(G)$, then give the definition of path complex of a digraph $G$.
\begin{itemize}
  \item The path complex in Definition \ref{pathcomplex} is a sub-complex of the original one. If the digraph does not contain a monotone cycle, then the s-regular condition is the same as the regular condition. It follows that the corresponding homologies are the same.
  \item As \cite{HY} write, with this strongly regular condition, the homology groups (see, Definition \ref{pathcomplex}) are now obviously bounded above. But there exist simple finite digraphs that have unbounded path homologies under the original regular condition. We will give such an example in Appendix \ref{A}. Thus, the two kinds of path homologies are not always the same.
\end{itemize}
\end{remark}

\begin{remark} (1) Similarly, one can define the reduced path complex $(\widetilde{\Omega}_n(G;\Z),\partial)_{n\geq-1}$ by
\begin{itemize}
  \item $\widetilde{\Omega}_n(G;\Z)=\Omega_n(G;\Z)$, for $n\geq0$, and $\Omega_{-1}(G;\Z)=\mathcal{A}_{-1}(G;\Z)=\Z$,
  \item $\partial: \Omega_{0}(G;\Z)\rightarrow\Omega_{-1}(G;\Z)$, $\partial(\sum_{v}a_vv)=\sum_{v}a_v$, where $a_v\in\Z$ and $v\in V(G)$.
\end{itemize}
We denote $\tilde{H}_*(G;\Z)$ for the reduced path homology of $G$. The relation between $\tilde{H}_*(G;\Z)$ and $H_*(G;\Z)$ is the same as the usual one.

(2) One can also consider the path homology over a more general coefficient. For example, Grigoryan-Jimenez-Muranov\cite{GJM}, Grigoryan-Jimenez-Muranov-Yau\cite{GJMY} studied the homology theory of path complex over any abelian group $K$. Grigoryan-Muranov-Yau\cite{GMY} studied the K\"unneth formula over the field $\K$ with characteristic 0.
\end{remark}
\vskip 0.2cm

\section{Minimal Path}\label{minimalpathsection}

In this section, we continue focusing on the $\mathbb{Z}$-coefficient and abbreviate $\Omega_*(G;\Z)$ (respectively, $\mathcal{A}_*(G;\Z)$) as $\Omega_*(G)$ (respectively, $\mathcal{A}_*(G)$). We will recall the definition of the minimal path in \cite{HY} and prove some structure theorems. And then, we will study the minimal $3$-paths in details.

\subsection{Definitions and basic properties}\label{defbasicpropsubsection}

For any $P=\sum_{p=1}^mc_pe_p\in \Omega_n(G)$ with $e_p\in \mathcal{A}_n(G)$ being elementary, we define
$$w(P)=\sum_{p=1}^m |c_p|,$$
and call it the width of the path $P$.

\begin{definition}\label{minimaldef} Let $P=\sum_{p=1}^mc_pe_p\in\Omega_n(G;\Z)$, if there exists
$$P'=\sum_{p=1}^md_pe_p\in\Omega_n(G;\Z)\setminus\{0\},$$
satisfying
\begin{enumerate}
  \item $|c_p-d_p|\leq |c_p|$ and $|d_p|\leq |c_p|$ for each $p=1,\ldots,m$;
  \item $w(P')<w(P)$;
\end{enumerate}
then we call $P'$ is smaller than $P$, and denote by $P'<P$. If there does not exist such $P'$, we call $P$ a minimal $n$-path.
\end{definition}
The following observations are immediate from this definition.

\begin{itemize}
  \item If $P$ is minimal, so is $-P$.
  \item If $P'<P$, then $P-P'\in\Omega_n(G;\Z)$ is also smaller than $P$.
  \item Any element in $\Omega_n(G;\Z)$ is a $\Z$-linear combination of minimal $n$-paths.
  \item Minimal $0$-path is just represented by a single point (up to a sign); minimal $1$-path is just represented by a directed edge $e_{12}$ (up to a sign).
\end{itemize}

\begin{example} Let us look at the minimal paths in the following digraphs $G_1$ and $G_2$.
\begin{figure}[h]
	\centering
	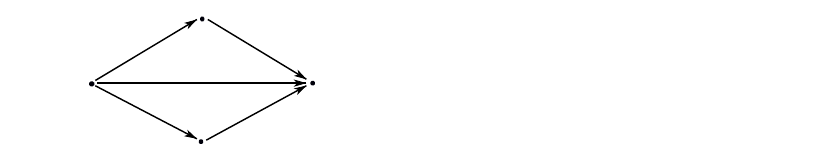
	%\caption[]{$\Supp(e_{0123})$}
	%\label{Fig:Q1}
\end{figure}
\begin{enumerate}
  \item For the digraph $G_1$, $P_1=e_{013}-e_{023}\in\Omega_2(G_1)$. But $P_1$ is not minimal, $e_{013}$ is smaller than $P_1$, and so is $-e_{023}=P_1-e_{013}$. All the minimal $2$-paths in $G_1$ are given by
      $$\pm e_{013}, ~\pm e_{023}.$$
  \item For the digraph $G_2$, $P_2=e_{014}-2e_{024}+e_{034}\in\Omega_2(G_2)$. But $P_2$ is not minimal, $e_{014}-e_{024}$ is smaller than $P_2$, so is $e_{024}-e_{034}=P_2-(e_{014}-e_{024})$. All the minimal $2$-paths in $G_2$ are given by
      $$\pm(e_{014}-e_{024}),~ \pm(e_{014}-e_{034}),~\pm(e_{024}-e_{034}).$$
\end{enumerate}

More examples will be studied in Subsection \ref{exminmal3} and Appendix \ref{B}.
\end{example}

The following result of minimal path is elementary but important for our later argument.

\begin{lemma}\label{minimalpath} Any minimal path is a $\Z$-linear combination of s-regular allowed elementary paths with the same starting vertex and ending vertex.
\end{lemma}

\begin{proof} Let $P$ be a minimal $n$-path of $G$, we write $P$ as
\begin{equation}\label{P}
P=P_1+P_2+\cdots+P_m\in\Omega_n(G),
\end{equation}
where for each $i=1,\ldots,m$, $P_i\in\mathcal{A}_n(G)$ is a $\Z$-linear combination of elementary $n$-paths with the same starting vertex $S_i$, and
$$S_i\neq S_j, \text{ for } i\neq j.$$
Let $\delta_0$ be $0$-th component of $\partial$, that is
$$\delta_0(e_{i_0i_1\ldots i_n})=e_{i_1i_2\ldots i_n}.$$
It is of course a map from $\mathcal{A}_*(G)$ to $\mathcal{A}_{*-1}(G)$. Thus by definition,
\begin{equation}\label{starting}
\partial P-\delta_0P=\sum_{i=1}^m(\partial-\delta_0)P_i\in\mathcal{A}_{n-1}(G).
\end{equation}
Note that $(\partial-\delta_0)P_i$ has the unique starting vertex $S_i$, and $S_i\neq S_j$ for $i\neq j$, then $(\partial-\delta_0)P_i$'s can not be cancelled by each other. Thus, \eqref{starting} implies
$$(\partial-\delta_0)P_i\in\mathcal{A}_{n-1}(G),\quad i=1,2,\ldots,m.$$
Furthermore,
$$P_i\in\Omega_n(G), \quad i=1,2,\ldots,m.$$
Again since different $P_i$'s have different starting vertices, we have $P_i\leq P$, for all $i$. The minimal condition for $P$ and the decomposition \eqref{P} imply that
$$P=P_i,\text{ for some } i,\quad \text{ and}\quad  P_j=0, \text{ for }j\neq i.$$
Thus the starting vertex of $P$ is unique. One can do the similar argument for the ending vertex.
\end{proof}

Since the minimal path depends on the digraph, to study its deep structure, we introduce the following definition of supporting digraphs.

\begin{definition}\label{supp} For each minimal path $P$ in the digraph $G$, we define $\Supp(P)$ to be the minimal sub-digraph of $G$ such that $P\in\Omega_*(\Supp(P))$.
\end{definition}

One can obtain the sub-digraph $\Supp(P)$ as follows:
\begin{enumerate}
  \item First, we express $P$ uniquely in terms of the elementary paths, that is
        $$P=\sum_{p=1}^lc_pe_p,\quad e_p\in\mathcal{A}_n(G),~ c_p\in\Z\setminus\{0\}.$$
        Then we get a digraph $\Supp(P)^{\mbox{pre}}$ given by $\{e_p\}_{p=1}^l$.
  \item Second, similarly, we express $\partial P$ uniquely in terms of the elementary paths, that is
        $$\partial P=\sum_{q=1}^md_qe_q,\quad e_q\in\mathcal{A}_{n-1}(G),~d_q\in\Z\setminus\{0\}.$$
        Then we add necessary edges from $\{e_q\}_{q=1}^m$ to $\Supp(P)^{\mbox{pre}}$ and obtain $\Supp(P)$.
\end{enumerate}

\begin{example}\label{exsuppdef} In the following 3-simplex digraph $G$, $P=e_{0123}$ is minimal, and its supporting sub-digraph is given by
\begin{figure}[H]
	\centering
	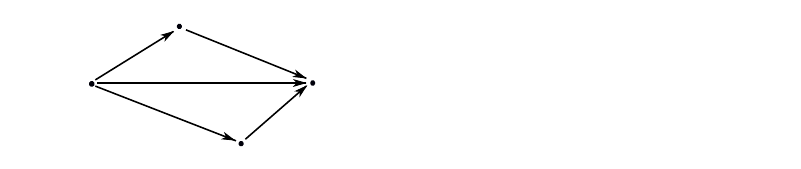
	%\caption[]{}
	%\label{Fig:Q1}
\end{figure}
\end{example}

For each s-regular allowed elementary path $e_I=e_{i_0i_1\ldots i_n}$, we have the position function
$$f_I:V(e_I)~\longrightarrow~\mathbb{N},\quad f_I(i_a)=a.$$
Let $P\in\Omega_n(G)$ be a minimal path with the starting vertex $S$ and the ending vertex $E$. For the convenience of the following discussion, we introduce the following two ``distance functions" on $V(P)$:
\begin{align*}
&d_S:V(P)\rightarrow \mathbb{N},\quad d_S(v)=\min\{f_I(v)|e_I\text{ is a component of }P\};\\
&d_E:V(P)\rightarrow \mathbb{N},\quad d_E(v)=\min\{n-f_I(v)|e_I\text{ is a component of }P\}.
\end{align*}
One can understand $d_S$ and $d_E$ as the distance from $v$ to $S$ and $E$ in $P$ (not in $\Supp(P)$) respectively. Clearly,
$$d_S(S)=0, \quad d_S(E)=n, \quad d_{E}(S)=n,\quad d_E(E)=0,\quad d_S(v)+d_E(v)\leq n.$$

\begin{example}\label{SabcdeE} For the following digraph $G$, the $3$-path $P=e_{SacE}-e_{SaeE}+e_{SdeE}+e_{SbdE}-e_{SbcE}$ is a minimal path.
\begin{figure}[H]
	\centering
	%% Creator: Inkscape 1.0.1 (3bc2e813f5, 2020-09-07), www.inkscape.org
%% PDF/EPS/PS + LaTeX output extension by Johan Engelen, 2010
%% Accompanies image file 'EXfSfE.pdf' (pdf, eps, ps)
%%
%% To include the image in your LaTeX document, write
%%   \input{<filename>.pdf_tex}
%%  instead of
%%   \includegraphics{<filename>.pdf}
%% To scale the image, write
%%   \def\svgwidth{<desired width>}
%%   \input{<filename>.pdf_tex}
%%  instead of
%%   \includegraphics[width=<desired width>]{<filename>.pdf}
%%
%% Images with a different path to the parent latex file can
%% be accessed with the `import' package (which may need to be
%% installed) using
%%   \usepackage{import}
%% in the preamble, and then including the image with
%%   \import{<path to file>}{<filename>.pdf_tex}
%% Alternatively, one can specify
%%   \graphicspath{{<path to file>/}}
%% 
%% For more information, please see info/svg-inkscape on CTAN:
%%   http://tug.ctan.org/tex-archive/info/svg-inkscape
%%
\begingroup%
  \makeatletter%
  \providecommand\color[2][]{%
    \errmessage{(Inkscape) Color is used for the text in Inkscape, but the package 'color.sty' is not loaded}%
    \renewcommand\color[2][]{}%
  }%
  \providecommand\transparent[1]{%
    \errmessage{(Inkscape) Transparency is used (non-zero) for the text in Inkscape, but the package 'transparent.sty' is not loaded}%
    \renewcommand\transparent[1]{}%
  }%
  \providecommand\rotatebox[2]{#2}%
  \newcommand*\fsize{\dimexpr\f@size pt\relax}%
  \newcommand*\lineheight[1]{\fontsize{\fsize}{#1\fsize}\selectfont}%
  \ifx\svgwidth\undefined%
    \setlength{\unitlength}{172.73074894bp}%
    \ifx\svgscale\undefined%
      \relax%
    \else%
      \setlength{\unitlength}{\unitlength * \real{\svgscale}}%
    \fi%
  \else%
    \setlength{\unitlength}{\svgwidth}%
  \fi%
  \global\let\svgwidth\undefined%
  \global\let\svgscale\undefined%
  \makeatother%
  \begin{picture}(1,0.45592596)%
    \lineheight{1}%
    \setlength\tabcolsep{0pt}%
    \put(0,0){\includegraphics[width=\unitlength,page=1]{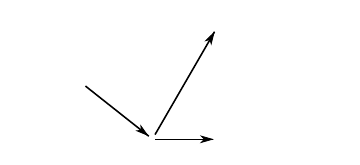}}%
    \put(0.16066841,0.19840504){\makebox(0,0)[lt]{\lineheight{1.25}\smash{\begin{tabular}[t]{l}$S$\end{tabular}}}}%
    \put(0.40156266,0.41961927){\makebox(0,0)[lt]{\lineheight{1.25}\smash{\begin{tabular}[t]{l}$a$\end{tabular}}}}%
    \put(0.40024298,0.00636223){\makebox(0,0)[lt]{\lineheight{1.25}\smash{\begin{tabular}[t]{l}$b$\end{tabular}}}}%
    \put(0.58477661,0.42404832){\makebox(0,0)[lt]{\lineheight{1.25}\smash{\begin{tabular}[t]{l}$c$\end{tabular}}}}%
    \put(0.59230915,0.00445021){\makebox(0,0)[lt]{\lineheight{1.25}\smash{\begin{tabular}[t]{l}$d$\end{tabular}}}}%
    \put(0.71992686,0.14505418){\makebox(0,0)[lt]{\lineheight{1.25}\smash{\begin{tabular}[t]{l}$e$\end{tabular}}}}%
    \put(0,0){\includegraphics[width=\unitlength,page=2]{EXfSfE.pdf}}%
    \put(-0.00393576,0.1948735){\makebox(0,0)[lt]{\lineheight{1.25}\smash{\begin{tabular}[t]{l}$G=$\end{tabular}}}}%
    \put(0,0){\includegraphics[width=\unitlength,page=3]{EXfSfE.pdf}}%
    \put(0.93143075,0.20313626){\makebox(0,0)[lt]{\lineheight{1.25}\smash{\begin{tabular}[t]{l}$E$\end{tabular}}}}%
  \end{picture}%
\endgroup%

	%\caption[]{$\Supp(e_{0123})$}
	%\label{Fig:Q1}
\end{figure}
We have
\begin{align*}
&d_S(a)=1, d_E(a)=2;\quad d_S(b)=1, d_E(b)=2;\quad d_S(c)=2, d_E(c)=1\\
&d_S(d)=1, d_E(d)=1;\quad d_S(e)=2, d_E(e)=1.
\end{align*}
\end{example}

Now let us decompose the differential $\partial=\partial_n: \Omega_n(G)\rightarrow\Omega_{n-1}(G)$ according to the positions of the vertices in the path
$$\partial_n=\delta_0+\delta_1+\cdots+\delta_{n}.$$

Immediately, we have

\begin{lemma}\label{deltaiP} Let $P\in\Omega_n(G)$ be a minimal path, then
$$\delta_iP\in \mathcal{A}_{n-1}(\Supp(P))\subset\mathcal{A}_{n-1}(G),\quad \text{for each }i=0,\ldots,n.$$
Moreover, for any s-regular allowed elementary $n$-path $A$, $A\in\Omega_n(G)$ if and only if
$$\delta_iA\in \mathcal{A}_{n-1}(G),\quad \text{for each }i=0,\ldots,n.$$
\end{lemma}

\begin{proof} For $n=1$, the result holds naturally. Let us assume that the result holds for any $k<n$, i.e. for any minimal $k$-path $P_k$ in $G$,
$$\delta_iP_k\in\mathcal{A}_{k-1}(\Supp(P_k))\subset\mathcal{A}_{k-1}(G),\quad \text{for each }i=0,\ldots,k.$$
Now let us consider the case $k=n$, first since $\partial P\in \Omega_{n-1}(\Supp(P))$, we can expand it in terms of minimal paths in $\Omega_{n-1}(\Supp(P))$ and organize it according to the starting and ending points
\begin{equation}\label{minimaldecomp}
\partial P=\sum_{\alpha\in E_1}\sum_{j=1,\ldots,i_{\alpha}}c_{\alpha}^jP_{S,n-1,\alpha}^j+\sum_{\beta\in S_1}A_{\beta,n-1,E}+A_{S,n-1,E},
\end{equation}
where
\begin{itemize}
  \item $E_1=d_E^{-1}(1)$, $S_1=d_S^{-1}(1)$;
  \item For each $\alpha\in E_1$, and $j=1,\ldots,i_{\alpha}$,
  $$P_{S,n-1,\alpha}^j\in\Omega_{n-1}(\Supp(P);\Z)$$
  is the minimal $(n-1)$-path starting from $S$ and ending with $\alpha$.
  \item $A_{\beta,n-1,E}\in\Omega_{n-1}(\Supp(P))$ is starting from $\beta$ and ending with $E$;
  \item $A_{S,n-1,E}\in\Omega_{n-1}(\Supp(P))$ is starting from $S$ and ending with $E$.
\end{itemize}
It is easy to obtain
$$\delta_nP=\sum_{\alpha\in E_1}\sum_{j=1,\ldots,i_{\alpha}}c_{\alpha}^jP_{S,n-1,\alpha}^j.$$
By Lemma \ref{minimalpath}, we have
$$P=(-1)^n\sum_{\alpha\in E_1}\sum_{j=1,\ldots,i_{\alpha}}c_{\alpha}^jP_{S,n-1,\alpha}^jE,$$
where $P_{S,n-1,\alpha}^jE$ means the concatenation of the path $P_{S,n-1,\alpha}^j$ and $e_{\alpha E}$. That is, if we write $P_{S,n-1,\alpha}^j$ as $P_{S,n-1,\alpha}^j=\sum_v d_ve_{Sv_1\ldots v_{n-2}\alpha}$, then
$$P_{S,n-1,\alpha}E=\sum_v d_v e_{Sv_1\ldots v_{n-2}\alpha E}.$$

Now for $i=0,\ldots,n-1$, we can compute
\begin{equation}\label{deltai}
\delta_iP=(-1)^{n}\sum_{\alpha\in E_1}\sum_{j=1,\ldots,i_{\alpha}}\bigg(c_{\alpha}^j\delta_iP_{S,n-1,\alpha}^j\bigg)E.
\end{equation}
By induction hypothesis, for each $\alpha\in E_1$, $j=1,\ldots, i_{\alpha}$ and for each $i=0,\ldots,n-1$ we have
$$\delta_iP_{S,n-1,\alpha}^j\in \mathcal{A}_{n-2}(\Supp(P_{S,n-1,\alpha}^j))\subset\mathcal{A}_{n-2}(\Supp(P))\subset\mathcal{A}_{n-2}(G).$$
In particular, all the non-zero path components of $\delta_0P_{S, n-1,\alpha}^j,\delta_1P_{S,n-1,\alpha}^j,\ldots,\delta_{n-2}P_{S,n-1,\alpha}^j$ end with $\alpha$, which means for each $k=0,\ldots,n-2$,
$$\delta_kP=(-1)^{n}\sum_{\alpha\in E_1}\sum_{j=1,\ldots,i_{\alpha}}\bigg(c_{\alpha}^j\delta_kP_{S,n-1,\alpha}^j\bigg)E\in \mathcal{A}_{n-1}(\Supp(P))\subset\mathcal{A}_{n-1}(G).$$
It is obvious that $\delta_nP\in \mathcal{A}_{n-1}(\Supp(P))\subset\mathcal{A}_{n-1}(G)$, thus
$$\delta_{n-1}P=\partial P-\sum_{k=0}^{n-2}\delta_kP-\delta_nP\in\mathcal{A}_{n-1}(\Supp(P))\subset\mathcal{A}_{n-1}(G).$$

For any $A\in\Omega_n(G)$, we can write $A=\sum_P a_PP$. Then we are done.
\end{proof}

Moreover, let us decompose the minimal path $P$ according to the $k$-th vertices in each elementary path component of $P$. That is,
\begin{equation}\label{PiA}
P=\sum_Ic_Ie_I=\sum_{v\in \{f_I^{-1}(k)\}_I}A_{SvE},
\end{equation}
where $A_{SvE}\in\mathcal{A}_n(\Supp(P))$ is the sum of elementary path components of $P$ with the $k$-th vertex being $v$. Then we have
\begin{lemma}\label{PiA2} Let $P\in\Omega_n(G)$ be a minimal path with the decomposition \eqref{PiA}, then for each $v\in \{f_I^{-1}(k)\}_I$,
$$\delta_j A_{SvE}\in\mathcal{A}_{n-1}(\Supp(P)), \quad j=0,1,\ldots,k-1,k+1,\ldots,n.$$
\end{lemma}

\begin{proof} The proof is almost the same as above. For $j\neq k$, by Lemma \ref{deltaiP}, we have
\begin{equation}\label{deltajP}
\delta_jP=\sum_{v\in \{f_I^{-1}(k)\}_I}\delta_j A_{SvE}\in\mathcal{A}_{n-1}(\Supp(P)).
\end{equation}
Since $\delta_jA_{SvE}$ does not kill the $k$-th vertex $v$ for $j\neq k$, then any different $v,v'\in\{f_I^{-1}(k)\}_I$, the paths $\delta_jA_{SvE}$ and $\delta_jA_{Sv'E}$ can not cancel each other, by \eqref{deltajP}, it means,
$$\delta_jA_{SvE}\in \mathcal{A}_{n-1}(\Supp(P)).$$
The remaining statement is obvious.
\end{proof}

According to our decomposition, we can write the above $A_{SvE}$ as, %by \eqref{ASvE}, we also write the above path $A_{SvE}$ as
$$A_{SvE}=A_{S,k,v}\bullet A_{v,n-k,E},$$
where $A_{S,k,v}$ and $A_{v,n-k,E}$ are the front $k$-path and back $(n-k)$-path of $A_{SvE}$ respectively, and $\bullet$ by connecting two paths $A_{S,k,v}$ and $A_{v,n-k,E}$ at the point $v$.

\begin{example} Let us look at Example \ref{SabcdeE} again, we can expand the minimal path $P$ according to the $1$-st and $2$-nd positions:
\begin{align*}
&P=e_{Sa}\bullet(e_{acE}-e_{aeE})-e_{Sb}\bullet(e_{bcE}-e_{bdE})+e_{Sd}\bullet e_{deE};\\
&P=(e_{Sac}-e_{Sbc})\bullet e_{cE}+e_{Sbd}\bullet e_{dE}-(e_{Sae}-e_{Sde})\bullet e_{eE}.
\end{align*}
\end{example}

Let us state another simple structure result for $\Supp(P)$ before our main structure theorem, which will be useful in our proof of acyclic result in Section \ref{acyclicsection}.

\begin{lemma}\label{suppedge} Let $P\in\Omega_n(G)$ be a minimal path and $f_I$ be the position function for each path component $e_I$ of $P$. Then there are not directed edges from points of $f_I^{-1}(k)$ to points of $f_I^{-1}(l)$, if $|l-k|>2$.
\end{lemma}

\begin{proof} We prove the result by induction. For minimal paths of length $0,1,2$, it is obvious. Assume the result holds for minimal path of length $k<n$. Now for minimal $n$-path $P$, we have seen that
$$P=(-1)^n\sum_{\alpha\in E_1}\sum_{j=1,\ldots,i_{\alpha}}c_{\alpha}^jP_{S,n-1,\alpha}^jE.$$
Similarly, one can also write $P$ as
$$P=\sum_{\beta\in S_1}\sum_{k=1,\ldots,j_{\beta}}d_{\beta}^kSP_{\beta,n-1,E}^k.$$
Applying the induction hypothesis on the two kinds of the expression for $P$, we are done.
\end{proof}

\subsection{The structure theorem}\label{strthmsubsection}
Now, let us improve the decomposition result for $\partial P$ and explore the deep structure of the minimal path. Let us state our first main theorem as follows.

\begin{theorem}[Structure Theorem]\label{structurethm} Let $P\in\Omega_n(G;\Z)$ be a minimal path with the starting vertex $S$ and ending vertex $E$, $\Supp(P)$ be its supporting digraph and $d_S$, $d_E$ be the functions defined above.

(1) Let $S_1=d_S^{-1}(1)$ and $E_1=d_E^{-1}(1)$. Then
$P$ is a linear combination of s-regular allowed elementary paths from $S$ to $E$, with coefficients being either 1 or -1. And in $\Supp(P)$,
\begin{equation*}
\begin{aligned}
\partial P=\sum_{\alpha\in E_1}P_{S,n-1,\alpha}+\sum_{\beta\in S_1}P_{\beta,n-1,E}+\sum_{k\in I_P} P_{S,n-1,E}^k,
\end{aligned}
\end{equation*}
where
\begin{itemize}
  \item  $P_{S,n-1,\alpha},P_{\beta,n-1,E},P_{S,n-1,E}^k\in\Omega_{n-1}(\Supp(P);\Z)$, and moreover
  \begin{itemize}
          \item $P_{S,n-1,\alpha}$ is the minimal $(n-1)$-path starting from $S$ and ending with $\alpha$;
          \item $P_{\beta,n-1,E}$ is the minimal $(n-1)$-path starting from $\beta$ and ending with $E$;
          \item $P_{S,n-1,E}^k$ is the minimal $(n-1)$-path starting from $S$ and ending with $E$.
        \end{itemize}
  \item Such $P_{S,n-1,\alpha}$, $P_{\beta,n-1,E}$ are unique (up to a sign)\footnote{such a sign is fixed by $P$.} in $\Supp(P)$ for each $\alpha\in E_1$, $\beta\in S_1$.
  \item The set $I_P$ in the last summand depends on $P$, and $|I_P|\leq 1$. That is, there exists at most one (up to a sign) minimal $(n-1)$-path in $\Supp(P)$ starting from $S$ and ending with $E$.
\end{itemize}

(2) For any $v\in d_E^{-1}(k)\cap\Supp(P)$, in $\Supp(P)$,
\begin{itemize}
  \item There is a unique minimal $(n-k)$-path (up to a sign) starting from $S$ and ending with $v$, denoted by $P_{S,n-k,v}$, as well as a unique minimal $k$-path (up to a sign) starting from $v$ and ending with $E$, denoted by $P_{v,k,E}$.
  \item There is at most one minimal $(n-k-1)$-path (up to a sign) starting from $S$ and ending with $v$, denoted by $P_{S,n-k-1,v}$, as well as at most one minimal $(k-1)$-path starting from $v$ and ending with $E$, denoted by $P_{v,k-1,E}$.
\end{itemize}

(3) Each minimal $2$-path with fixed starting and ending vertices in $\Supp(P)$ ($n\geq2$) is unique.
\end{theorem}

We will prove the results (1) (2) (3) at the same time by induction. Before the general proof, let us look at the minimal $2$-path. There are only two kinds of $2$-minimal paths: write the minimal path $P$ as
$$P=\sum_{a\in V(P)\setminus\{S,E\}}c_ae_{SaE},\quad c_a\in\mathbb{Z}.$$
By Lemma \ref{deltai}, the requirement $\delta_1P=-\sum_ac_ae_{SE}\in\mathcal{A}_1(G)$ tells us
\begin{itemize}
  \item If $S\rightarrow E$ in $G$, then for $a$ with $c_a\neq0$, $\widetilde{P}:=\sign(c_a)e_{SaE}\leq P$. By the minimal condition for $P$, the equality holds if and only if $c_a=\pm 1$, and $c_{b}=0$ for $b\neq a$. The corresponding sub-digraph $\Supp(P)$ is the triangle
      \begin{equation*}
      \Supp(P)=%
      \begin{array}{ccccc}
      &  & a &  &  \\
      & \nearrow &  & \searrow &  \\
      S &  & \longrightarrow &  & E%
      \end{array}%
      \end{equation*}
  \item If $S\nrightarrow E$ in $G$, then
  $$\delta_1P=-\sum_ac_ae_{SE}=0.$$
  which means there exist $a,b$ such that $c_a,c_b\neq0$, thus $\widetilde{P}:=\sign(c_a)(e_{SaE}-e_{SbE})\leq P$. By the minimal condition for $P$, the equality holds if and only if $c_a=-c_b=\pm 1$, and other $c_{i}=0$ for $i\neq a,b$. The corresponding sub-digraph $\Supp(P)$ is the square
  \begin{equation*}
  \Supp(P)=%
  \begin{array}{ccc}
  a & \longrightarrow & E \\
  \uparrow &  & \uparrow \\
  S & \longrightarrow & b%
  \end{array}
  \end{equation*}
\end{itemize}

Now let us illustrate the proof idea of the theorem. First, by the rough decomposition \eqref{minimaldecomp}, we can write
\begin{equation}\label{Prough}
P=(-1)^n\sum_{\alpha\in E_1}\sum_{j=1,\ldots,i_{\alpha}}c_{\alpha}^jP_{S,n-1,\alpha}^jE,\quad c_{\alpha}^j\in \Z-\{0\}.
\end{equation}
each $P_{S,n-1,\alpha}^{j}\in\Omega_{n-1}(\Supp(P))$ is the minimal path starting from $S$ and ending with $\alpha$.
If $i_a>1$, or if there exists some $c_{\alpha}^j$, $|c_{\alpha}^j|\geq2$, we will construct a path $\widetilde{P}\in\Omega_{n}(\Supp(P))$ which is strictly smaller than $P$, as we did in the length $2$ case. To simplify the combinatorial discussion in the proof, we generalize the partial order on $\Omega_n(G)$ in Definition \ref{minimaldef} to the partial order on $\mathcal{A}_n(G)$ and introduce the definition of $\partial$-invariant completion.

\begin{definition}\label{partialorder} Let $u=\sum_{p=1}^mu_pe_p$, $v=\sum_{p=1}^mv_pe_p\in\mathcal{A}_n(G)$, $u_p,v_p\in\Z$. We say $u< v$ if
\begin{enumerate}
  \item $|u_p|\leq|v_p|$ and $|u_p-v_p|\leq|v_p|$, for each $p=1,\ldots,m$;
  \item $\sum_{p=1}^m|u_p|<\sum_{p=1}^m|v_p|$.
\end{enumerate}
We write $u\leq v$ if $u<v$ or $u=v$.
\end{definition}

Immediately, we have
\begin{proposition} For each $n\in\mathbb{N}$, the above $\leq$ on $\mathcal{A}_n(G)$ defines a partial order on $\mathcal{A}_n(G)$. In particular, for $u=\sum_{n=1}^mu_ne_n\in\mathcal{A}_n(G)$, with $u_n\in\Z\setminus\{0\}$, we have
$$\sign(u_n)e_n\leq u,\quad n=1,\ldots, m.$$
\end{proposition}

\begin{definition}[$\partial$-invariant completion]\label{invcompletion} Let $u\in\mathcal{A}_n(G)$. We say $u$ admits a $\partial$-invariant completion, if there exists a path $\tilde{u}\in\Omega_n(G)$ such that $u\leq \tilde{u}$. Furthermore, if $u$ admits $\partial$-invariant completions, let
$$C_{\min}(G,u)=\{\tilde{u}\in\Omega_n(G)|u\leq \tilde{u},\text{there does not exist $v\in\Omega_n(G)$, such that $v<\tilde{u}$.}\}$$
be the subset of the minimal $\partial$-invariant completions of $u$.
\end{definition}

It is obvious that if $u_{S\cdots E}\in\mathcal{A}_n(G)$ is a s-regular allowed path with starting vertex $S$ and ending vertex $E$ and admits $\partial$-invariant completions, then any $\tilde{u}\in C_{min}(G,u_{S\cdots E})$ has the same starting vertex $S$ and the same ending vertex $E$.

\begin{example} For the path $P=ae_{024}\in\mathcal{A}_2(G)$, $a\in\mathbb{N}_+$ in the following digraph $G$,
\begin{figure}[H]
	\centering
	%% Creator: Inkscape 1.0.1 (3bc2e813f5, 2020-09-07), www.inkscape.org
%% PDF/EPS/PS + LaTeX output extension by Johan Engelen, 2010
%% Accompanies image file 'invcompletion.pdf' (pdf, eps, ps)
%%
%% To include the image in your LaTeX document, write
%%   \input{<filename>.pdf_tex}
%%  instead of
%%   \includegraphics{<filename>.pdf}
%% To scale the image, write
%%   \def\svgwidth{<desired width>}
%%   \input{<filename>.pdf_tex}
%%  instead of
%%   \includegraphics[width=<desired width>]{<filename>.pdf}
%%
%% Images with a different path to the parent latex file can
%% be accessed with the `import' package (which may need to be
%% installed) using
%%   \usepackage{import}
%% in the preamble, and then including the image with
%%   \import{<path to file>}{<filename>.pdf_tex}
%% Alternatively, one can specify
%%   \graphicspath{{<path to file>/}}
%% 
%% For more information, please see info/svg-inkscape on CTAN:
%%   http://tug.ctan.org/tex-archive/info/svg-inkscape
%%
\begingroup%
  \makeatletter%
  \providecommand\color[2][]{%
    \errmessage{(Inkscape) Color is used for the text in Inkscape, but the package 'color.sty' is not loaded}%
    \renewcommand\color[2][]{}%
  }%
  \providecommand\transparent[1]{%
    \errmessage{(Inkscape) Transparency is used (non-zero) for the text in Inkscape, but the package 'transparent.sty' is not loaded}%
    \renewcommand\transparent[1]{}%
  }%
  \providecommand\rotatebox[2]{#2}%
  \newcommand*\fsize{\dimexpr\f@size pt\relax}%
  \newcommand*\lineheight[1]{\fontsize{\fsize}{#1\fsize}\selectfont}%
  \ifx\svgwidth\undefined%
    \setlength{\unitlength}{175.33805547bp}%
    \ifx\svgscale\undefined%
      \relax%
    \else%
      \setlength{\unitlength}{\unitlength * \real{\svgscale}}%
    \fi%
  \else%
    \setlength{\unitlength}{\svgwidth}%
  \fi%
  \global\let\svgwidth\undefined%
  \global\let\svgscale\undefined%
  \makeatother%
  \begin{picture}(1,0.4171721)%
    \lineheight{1}%
    \setlength\tabcolsep{0pt}%
    \put(0,0){\includegraphics[width=\unitlength,page=1]{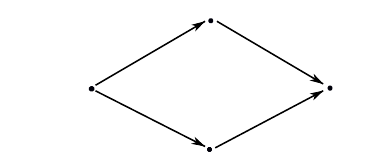}}%
    \put(0.16988147,0.15704767){\makebox(0,0)[lt]{\lineheight{1.25}\smash{\begin{tabular}[t]{l}$0$\end{tabular}}}}%
    \put(0.53764957,0.3835993){\makebox(0,0)[lt]{\lineheight{1.25}\smash{\begin{tabular}[t]{l}$1$\end{tabular}}}}%
    \put(0.54194331,0.05005086){\makebox(0,0)[lt]{\lineheight{1.25}\smash{\begin{tabular}[t]{l}$3$\end{tabular}}}}%
    \put(0.92738823,0.15399381){\makebox(0,0)[lt]{\lineheight{1.25}\smash{\begin{tabular}[t]{l}$4$\end{tabular}}}}%
    \put(-0.00431638,0.15748622){\makebox(0,0)[lt]{\lineheight{1.25}\smash{\begin{tabular}[t]{l}$G=$\end{tabular}}}}%
    \put(0,0){\includegraphics[width=\unitlength,page=2]{invcompletion.pdf}}%
    \put(0.54085256,0.20806403){\makebox(0,0)[lt]{\lineheight{1.25}\smash{\begin{tabular}[t]{l}$2$\end{tabular}}}}%
  \end{picture}%
\endgroup%

	%\caption[]{$\Supp(e_{0123})$}
	%\label{Fig:Q1}
\end{figure}
By definition, we see that
$$\tilde{P}:=be_{024}-(b+c)e_{014}+ce_{034},\quad \text{for any } c\in\Z, b\in\mathbb{N}_+,b\geq a$$
is a $\partial$-invariant completion for $P$. In particular, we have
$$C_{\min}(G,ae_{024})=\{ae_{024}-ae_{014}, ae_{024}-ae_{034}\}.$$
\end{example}

\subsection{Proof of Theorem \ref{structurethm}}\label{strproofsubsection}

To prove the theorem, we consider one more claim as follows.

\begin{claim}\label{uniqueclaim} (1) For $2\leq j\leq n$, each front elementary $j$-path component $e_{Si_1\ldots i_{j-1}u}$ of $P$ admits a unique minimal $\partial$-invariant completion in $\Supp(P)$, denote it by $P_{S,j,u}$.

(2) For $0\leq k\leq n-2$, each back elementary $(n-k-1)$-path component $e_{vi_{k+1}\ldots i_{n-1}E}$ of $P$ admits a unique minimal $\partial$-invariant completion in $\Supp(P)$, denote it by $P_{v,n-k-1,E}$.
\end{claim}

We prove our result as well as the claim by inductions, while the induction for the claim is in the decreasing order.
\begin{description}
  \item[\textbf{Initial result 1}] It is clear that the theorem holds for minimal paths of length $0,1,2$.
  \item[\textbf{Initial result 2 (for the claim)}] For $j=n$, the minimal $\partial$-invariant completion of $e_{Si_1\ldots i_{n-1}E}$ is exactly $P$, since $P$ is minimal. (The proof for the second part is similar.)
\end{description}

\begin{description}
  \item[\textbf{Hypothesis 1}] Assume that the result holds for minimal paths of length $k<n$. In particular, for each $P_{S,n-1,\alpha}^j$ in \eqref{Prough}, we have
\begin{equation}\label{strassumpn-1}
P_{S,n-1,\alpha}^j=(-1)^{n-1}\sum_{\alpha'\in d_{\alpha,j}^{-1}(1)}P_{S,n-2,\alpha'}^j\alpha,
\end{equation}
where $d_{\alpha,j}$ is the distance function on $V(P_{S,n-1,\alpha}^j)$, $P_{S,n-2,\alpha'}^j\in \Omega_{n-2}(\Supp(P_{S,n-1,\alpha}^j))$ is the unique minimal $(n-2)$-path starting from $S$ and ending with $\alpha'\in d_{\alpha,j}^{-1}(1)$.
  \item[\textbf{Hypothesis 2}] Assume that the claim holds for $k+1\leq j\leq n$.
\end{description}

Now let us look at the second induction first. For the case $j=k$, let $e_{Si_1\ldots i_{k-1}u}$ be such a front $k$-path component of $P$, by Lemma \ref{PiA2}, it admits a $\partial$-invariant completion in $\Supp(P)$. Now assume that $P_{S,k,u}^1$ and $P_{S,k,u}^2$ are two minimal $\partial$-invariant completions of $e_{Si_1\ldots i_{k-1}u}$ in $\Supp(P)$. By Hypothesis 1, both $P_{S,k,u}^1$ and $P_{S,k,u}^2$ are also minimal paths.

Now let us $e_{Si_1\ldots i_{k-1}uv}$ be the front $(k+1)$-path component of $P$. Then
\begin{itemize}
  \item by Hypothesis 2, it admits a unique minimal $\partial$-invariant completion $P_{S,k+1,v}$ in $\Supp(P)$.
  \item by Hypothesis 1, there is only one minimal $k$-path starting from $S$ and ending with $u$ in $\Supp(P_{S,k+1,v})$.
\end{itemize}
The two uniqueness results imply that $P_{S,k,u}^1=P_{S,k,u}^2$. Then we finish our proof of the claim.

Now let us turn to our first induction. For the minimal $n$-path $P$, recall we can write it as
\begin{equation}\label{Prough2}
P=(-1)^n\sum_{\alpha\in E_1}\sum_{j=1,\ldots,i_{\alpha}}c_{\alpha}^jP_{S,n-1,\alpha}^jE,\quad c_{\alpha}^j\in \Z-\{0\}.
\end{equation}

\subsubsection{The case $|E_1|=1$}
Let us start with the simplest case that $|E_1|=1$ with $E_1=\{\alpha\}$. By \eqref{strassumpn-1} and \eqref{Prough2}, we have
\begin{equation}\label{PE11}
P=(-1)^n\sum_{j=1,\ldots,i_{\alpha}}c_{\alpha}^jP_{S,n-1,\alpha}^jE=-\sum_{j=1,\ldots,i_{\alpha}}\sum_{\gamma\in d_{\alpha,j}^{-1}(1)}c_{\alpha}^jP_{S,n-2,\gamma}^j\alpha E.
\end{equation}
\begin{itemize}
  \item If there exists $j\in\{1,\ldots,i_{\alpha}\}$ such that all vertices $\gamma\in d_{\alpha,j}^{-1}(1)$, $\gamma\rightarrow E$, then
      $$\widetilde{P}:=P_{S,n-1,\alpha}^jE\in\Omega_n(\Supp(P)).$$
      Thus, $\widetilde{P}\leq P$, or $\widetilde{P}\leq -P$, and the equality holds if and only if $i_{\alpha}=1$, $c_{\alpha}^j=\pm 1$.
  \item If for all $j\in\{1,\ldots,i_{\alpha}\}$, there exists a vertex $\gamma\in d_{\alpha,j}^{-1}(1)$ such that $\gamma\nrightarrow E$. Set
      \begin{align*}
      &I_j=\left\{\gamma\in d_{\alpha,j}^{-1}(1)~\big|~ \gamma\nrightarrow E\right\},\\
      &I_j^c=d_{\alpha,j}^{-1}(1)-I_j.
      \end{align*}
      Then we consider $\delta_{n-1}P$, by \eqref{PE11},
      $$\delta_{n-1}P=(-1)^n\sum_{j=1,\ldots,i_{\alpha}}\sum_{\gamma\in I_j^c}c_{\alpha}^jP_{S,n-2,\gamma}^jE+(-1)^n\sum_{j=1,\ldots,i_{\alpha}}\sum_{\gamma\in I_j}c_{\alpha}^jP_{S,n-2,\gamma}^jE,$$
      where the first summands are allowed paths; while for each $\gamma\in I_j$, $P_{S,n-2,\gamma}^jE\notin\mathcal{A}_{n-1}(\Supp(P))$. \\
      Since $\delta_{n-1}P\in \mathcal{A}_{n-1}(\Supp(P))$, then we get
      $$\sum_{j=1,\ldots,i_{\alpha}}\sum_{\gamma\in I_j}c_{\alpha}^jP_{S,n-2,\gamma}^jE=0,$$
      which also means
      $$\sum_{j=1,\ldots,i_{\alpha}}\sum_{\gamma\in  I_j}c_{\alpha}^jP_{S,n-2,\gamma}^j\alpha E=0.$$
      Then we have
      $$P=\sum_{j=1,\ldots,i_{\alpha}}\sum_{\gamma\in I_j^c}c_{\alpha}^jP_{S,n-2,\gamma}^j\alpha E.$$
      Then we reduce to the first situation.
\end{itemize}

Let us finish the argument of the case $|E_1|=1$. For $P=P_{S,n-1,\alpha}E$, first by Hypothesis 1, we have
\begin{align*}
\partial P_{S,n-1,\alpha}=~&\delta_{n-1}P_{S,n-1,\alpha}+\delta_0 P_{S,n-1,\alpha}+(\delta_1+\cdots+\delta_{n-2})P_{S,n-1,\alpha}\\
                         =~&\sum_{\gamma\in d_{\alpha}^{-1}(1)} P_{S,n-2,\gamma}\alpha+\sum_{\beta\in d_S^{-1}(1)}P_{\beta,n-2,\alpha} +P_{S,n-2,\alpha},
\end{align*}
where in particular, either $P_{S,n-2,\alpha}$ is a minimal $(n-1)$-path or $P_{S,n-2,\alpha}=0$.
Similarly, we have
\begin{align*}
\partial P=~&\partial (P_{S,n-1,\alpha}E)=\delta_n (P_{S,n-1,\alpha}E)+\delta_0 (P_{S,n-1,\alpha}E)+(\delta_1+\cdots+\delta_{n-1})(P_{S,n-1,\alpha}E)\\
          =~&(-1)^n P_{S,n-1,\alpha}+\sum_{\beta\in d_S^{-1}(1)}P_{\beta,n-2,\alpha}E + \left(P_{S,n-2,\alpha}E-\sum_{\gamma\in d_{\alpha}^{-1}(1)} P_{S,n-2,\gamma}E\right).
\end{align*}
The remaining statements are trivial except the following one:
$$P_{S,n-1,E}:=P_{S,n-2,\alpha}E-\sum_{\gamma\in d_{\alpha}^{-1}(1)} P_{S,n-2,\gamma} E$$
is a minimal $(n-1)$-path. (It is obvious non-zero).

Clearly, $P_{S,n-1,E}\in\Omega_{n-1}(\Supp(P))$, we only need to prove it is minimal. If it is not minimal, we can write it as a $\Z$-linear combination of minimal paths in $\Supp(P)$:
$$P_{S,n-1,E}=P_{S,n-1,E}^1+P_{S,n-1,E}^2+\cdots+P_{S,n-1,E}^k,\quad k\geq2.$$
According to the decomposition at the $(n-2)$-th vertices, then there exists some $i=1,\ldots,k$, and $I_i\subset d_{\alpha}^{-1}(1)$ such that
$$P_{S,n-1,E}^i=\sum_{\gamma\in I_i}P_{S,n-2,\gamma} E.$$
By our construction, we can find that
$$\widetilde{P}:=\sum_{\gamma\in I_i}P_{S,n-2,\gamma} \alpha E\in\Omega_n(G),\quad \widetilde{P}<P.$$
Contradiction. Thus, $P_{S,n-1,E}$ must be minimal. The statements (2) and (3) in this theorem follow from (1) naturally.

\subsubsection{The case $|E_1|=2$}

Now let us consider the case $|E_1|=2$, i.e. $E_1=\{\alpha_1,\alpha_2\}$. We write $P$ as
\begin{equation}\label{PE12}
\begin{aligned}
P=~&(-1)^n\sum_{j=1,\ldots,i_{\alpha_1}}c_{\alpha_1}^jP_{S,n-1,\alpha_1}^jE+(-1)^n\sum_{j=1,\ldots,i_{\alpha_2}}c_{\alpha_2}^jP_{S,n-1,\alpha_2}^jE\\
 =~&-\sum_{j=1,\ldots,i_{\alpha_1}}\sum_{\gamma_1\in d_{\alpha_1,j}^{-1}(1)}c_{\alpha_1}^jP_{S,n-2,\gamma_1}^j\alpha_1 E\\
 ~&-\sum_{j=1,\ldots,i_{\alpha_2}}\sum_{\gamma_1\in d_{\alpha_2,j}^{-1}(1)}c_{\alpha_2}^jP_{S,n-2,\gamma_2}^j\alpha_2 E.
\end{aligned}
\end{equation}

Now, without loss of generality, assume $i_{\alpha_1}\geq1$ or $c_{\alpha_1}^1\geq 1$. Let us first consider the path
$$P_0=(-1)^nP_{S,n-1,\alpha_1}^1E\in\mathcal{A}_n(\Supp(P)).$$
Then we have
\begin{itemize}
  \item $P_0$ is a $\Z$-linear combination of s-regular allowed elementary path starting from $S$ and ending with $E$ with coefficients $\pm 1$ by Hypothesis 1;
  \item $\delta_i P_0\in\mathcal{A}_{n-1}(\Supp(P))$, $i=0,1\ldots,n-2,n$;
  \item $P_0< P$.
\end{itemize}
We want to construct a $\partial$-invariant completion of $P_0$ that is smaller than $P$.

By Hypothesis 1, we can write $P_{S,n-1,\alpha_0}^1$ as
\begin{align*}
P_{S,n-1,\alpha_1}^1=~&(-1)^{n-1}\sum_{\gamma_1\in d_{\alpha_1,1}^{-1}(1)}P_{S,n-2,\gamma_1}^1\alpha_1\\
                    =~&(-1)^{n-1}\sum_{\gamma_1\in I_1^c}P_{S,n-2,\gamma_1}^1\alpha_1+(-1)^{n-1}\sum_{\gamma_1\in I_1}P_{S,n-2,\gamma_1}^1\alpha_1
\end{align*}
where as before,
\begin{align*}
&I_1=\left\{\gamma\in d_{\alpha_1,1}^{-1}(1)~\big|~ \gamma_1\nrightarrow E\right\},\\
&I_1^c=d_{\alpha_1,1}^{-1}(1)-I_1.
\end{align*}

Then, we have
\begin{align*}
\delta_{n-1}P_0=(-1)^n\sum_{\gamma_1\in I_1^c}P_{S,n-2,\gamma_1}^1E+(-1)^n\sum_{\gamma_1\in I_1}P_{S,n-2,\gamma_1}^1E,
\end{align*}
where the first summands are allowed paths; while for each $\gamma_1\in I_1$, $P_{S,n-2,\gamma_1}^1E$ is not allowed.

The following two observations help us modify $P_0$ to become $\delta_{n-1}$-invariant.
\begin{enumerate}
  \item Since $\delta_{n-1}P\in \mathcal{A}_{n-1}(\Supp(P))$, there must be terms in $\delta_{n-1}P$ which cancel each $P_{S,n-2,\gamma_1}^1E$.
  \item By Claim \ref{uniqueclaim}, $P_{S,n-2,\gamma_1}^1E$ must be canceled by terms coming from $\delta_{n-1}P_{S,n-1,\alpha_2}^jE$.
\end{enumerate}

Now let us consider the allowed path in $\Supp(P)$,
$$P_1=(-1)^{n-1}\sum_{\gamma_1\in I_1}P_{S,n-2,\gamma_1}^1\alpha_2E,$$
and modify the path $P_0$ to be $P_0+P_1$. Then by our construction, $P_0+P_1$ satisfies
\begin{itemize}
  \item $P_0+P_1$ is a $\Z$-linear combination of s-regular allowed elementary path starting from $S$ and ending with $E$ with coefficients $\pm 1$ by Hypothesis 1;
  \item $\delta_i(P_0+P_1)\in\mathcal{A}_{n-1}(\Supp(P))$, $i=0,\ldots,n-3,n-1,n$;
  \item $P_0+P_1\leq P$.
\end{itemize}

If $\delta_{n-2}P_1\in\mathcal{A}_{n-1}(\Supp(P))$, then we are done. Otherwise, we continue considering the minimal $\partial$-invariant completion of $P_0+P_1$.

Due to the fact that $\delta_{n-2}P\in\mathcal{A}_{n-2}(\Supp(P))$, Claim \ref{uniqueclaim} and Hypothesis 1, there exists a modified term $P_1'$ such that
\begin{itemize}
  \item $P_1+P_1'$ is a sum of some terms in $\{\pm P_{S,n-1,\alpha_2}^jE\}_{j=1}^{i_{\alpha_2}}$.
  \item $P_1+P_1'<P$.
\end{itemize}

Then we look at the path $P_0+P_1+P_1'$, it satisfies
\begin{itemize}
  \item $P_0+P_1+P_1'$ is a $\Z$-linear combination of s-regular allowed elementary path starting from $S$ and ending with $E$ with coefficients $\pm 1$ by Hypothesis 1;
  \item $\delta_i(P_0+P_1+P_1')\in\mathcal{A}_{n-1}(\Supp(P)),\quad i=0,\ldots,n-3,n-2,n$;
  \item $P_0+P_1+P_1'\leq P$.
\end{itemize}

If $\delta_{n-1}(P_0+P_1+P_1')\in\mathcal{A}_{n-1}(\Supp(P))$, we are done. Otherwise, we continue modifying our path, to get $P_2$, $P_2'$, $P_3$, $P_3'$,\ldots, there are two key points in the construction:
\begin{description}
  \item[Point 1] By Claim \ref{uniqueclaim}, the terms which cancel
$$\delta_{n-1}P_{S,n-2,\gamma_1}\alpha_1E,~\gamma_1\nrightarrow E\quad (\text{respectively, }\delta_{n-1}P_{S,n-2,\gamma_2}\alpha_2E, ~\gamma_2\nrightarrow E)$$
must come from
$$\delta_{n-1}P_{S,n-2,\gamma_1}\alpha_2E,~\gamma_1\nrightarrow E\quad (\text{respectively, }\delta_{n-1}P_{S,n-2,\gamma_2}\alpha_1E, ~\gamma_2\nrightarrow E).$$
At the $k$th time to preserve the sum $\delta_{n-1}$-invariant, we get the path $P_k$.
  \item[Point 2] As we construct the $P_1'$, we get $P_k'$. If $P_k'\neq0$, then new vertices in $E_2$ will come in.
\end{description}
Since we work with a finite digraph, it will stop at finite steps, which means, for some positive integer $l$, we obtain
$$\widetilde{P}=P_0+P_1+P_1'+\cdots +P_l +P_l'\in\Omega_n(\Supp(P)),\quad (P_l'\text{ may be }0),$$
is a $\Z$-linear combination of s-regular allowed elementary path starting from $S$ and ending with $E$ with coefficients $\pm1$; in particular, $\widetilde{P}\leq P$.

Moreover, let us look at our construction more carefully and see that the modification will stop before $P_2$. That is,
$$P_0+P_1+P_1'\in\Omega_n(\Supp(P)),\quad (P_1'\text{ may be }0).$$
Let $P_{S,n-3,\epsilon}\gamma_1\alpha_1E$, $\gamma_1\in I_1$ be a component of $P_0$. By our construction, first, we need to add a term $-P_{S,n-3,\epsilon}\gamma_1\alpha_2E$ (component of $P_1$) to do the modification.
\begin{itemize}
    \item If $\epsilon\rightarrow\alpha_2$, we do not need to modify this term again, done;
    \item If $\epsilon\nrightarrow\alpha_2$, we need to add the term of the form $P_{S,n-3,\epsilon}\gamma_2\alpha_2E$ (component of $P_1'$).
    \begin{itemize}
      \item If $\gamma_2\rightarrow E$, we do not need to modify this term again, done;
      \item If $\gamma_2\nrightarrow E$, and such a term $P_{S,n-3,\epsilon}\gamma_2\alpha_2E$ is not canceled by term from $P_0$. Then we need to add the term of the form $P_{S,n-3,\epsilon}\gamma_2\alpha_1$ (component of $P_2$). Then we will obtain the following sub-digraph.
          \begin{figure}[H]
	      \centering
	      %% Creator: Inkscape 1.0.1 (3bc2e813f5, 2020-09-07), www.inkscape.org
%% PDF/EPS/PS + LaTeX output extension by Johan Engelen, 2010
%% Accompanies image file 'contracdition.pdf' (pdf, eps, ps)
%%
%% To include the image in your LaTeX document, write
%%   \input{<filename>.pdf_tex}
%%  instead of
%%   \includegraphics{<filename>.pdf}
%% To scale the image, write
%%   \def\svgwidth{<desired width>}
%%   \input{<filename>.pdf_tex}
%%  instead of
%%   \includegraphics[width=<desired width>]{<filename>.pdf}
%%
%% Images with a different path to the parent latex file can
%% be accessed with the `import' package (which may need to be
%% installed) using
%%   \usepackage{import}
%% in the preamble, and then including the image with
%%   \import{<path to file>}{<filename>.pdf_tex}
%% Alternatively, one can specify
%%   \graphicspath{{<path to file>/}}
%% 
%% For more information, please see info/svg-inkscape on CTAN:
%%   http://tug.ctan.org/tex-archive/info/svg-inkscape
%%
\begingroup%
  \makeatletter%
  \providecommand\color[2][]{%
    \errmessage{(Inkscape) Color is used for the text in Inkscape, but the package 'color.sty' is not loaded}%
    \renewcommand\color[2][]{}%
  }%
  \providecommand\transparent[1]{%
    \errmessage{(Inkscape) Transparency is used (non-zero) for the text in Inkscape, but the package 'transparent.sty' is not loaded}%
    \renewcommand\transparent[1]{}%
  }%
  \providecommand\rotatebox[2]{#2}%
  \newcommand*\fsize{\dimexpr\f@size pt\relax}%
  \newcommand*\lineheight[1]{\fontsize{\fsize}{#1\fsize}\selectfont}%
  \ifx\svgwidth\undefined%
    \setlength{\unitlength}{164.05206876bp}%
    \ifx\svgscale\undefined%
      \relax%
    \else%
      \setlength{\unitlength}{\unitlength * \real{\svgscale}}%
    \fi%
  \else%
    \setlength{\unitlength}{\svgwidth}%
  \fi%
  \global\let\svgwidth\undefined%
  \global\let\svgscale\undefined%
  \makeatother%
  \begin{picture}(1,0.43324409)%
    \lineheight{1}%
    \setlength\tabcolsep{0pt}%
    \put(0,0){\includegraphics[width=\unitlength,page=1]{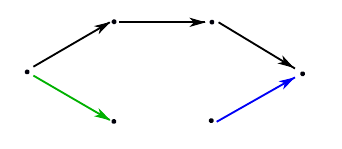}}%
    \put(-0.00360696,0.2060829){\makebox(0,0)[lt]{\lineheight{1.25}\smash{\begin{tabular}[t]{l}$\epsilon$\end{tabular}}}}%
    \put(0.29602835,0.39944377){\makebox(0,0)[lt]{\lineheight{1.25}\smash{\begin{tabular}[t]{l}$\gamma_1$\end{tabular}}}}%
    \put(0.58771189,0.400392){\makebox(0,0)[lt]{\lineheight{1.25}\smash{\begin{tabular}[t]{l}$\alpha_1$\end{tabular}}}}%
    \put(0.92318466,0.19960864){\makebox(0,0)[lt]{\lineheight{1.25}\smash{\begin{tabular}[t]{l}$E$\end{tabular}}}}%
    \put(0,0){\includegraphics[width=\unitlength,page=2]{contracdition.pdf}}%
    \put(0.29262194,0.00881757){\makebox(0,0)[lt]{\lineheight{1.25}\smash{\begin{tabular}[t]{l}$\gamma_2$\end{tabular}}}}%
    \put(0.58541266,0.01105614){\makebox(0,0)[lt]{\lineheight{1.25}\smash{\begin{tabular}[t]{l}$\alpha_2$\end{tabular}}}}%
    \put(0,0){\includegraphics[width=\unitlength,page=3]{contracdition.pdf}}%
  \end{picture}%
\endgroup%

          \end{figure}
    \end{itemize}
\end{itemize}

In this sub-digraph, we have a minimal $\partial$-invariant completion of $e_{\epsilon\gamma_1\alpha_1E}$, that is
$$P_{\epsilon,3,E}=e_{\epsilon\gamma_1\alpha_1E}-e_{\epsilon\gamma_2\alpha_1E}+e_{\epsilon\gamma_2\alpha_2E}-e_{\epsilon\gamma_1\alpha_2E}.$$

Since $\gamma_1\nrightarrow E$, $P_{\epsilon,3,E}$ is a minimal path. Apply Hypothesis 1 (3) for $P_{\epsilon,3,E}$, we have $\epsilon \nrightarrow \alpha_1$, otherwise $e_{\epsilon,\gamma_1,\alpha_1}$, $e_{\epsilon,\gamma_2,\alpha_1}$ become two different minimal $2$-faces with the same starting and ending vertices in $\Supp(P_{\epsilon,3,E})$. Since $P_{S,n-3,\epsilon,\gamma_1}\alpha_1<P_{S,n-1,\alpha_1}^1$, there must exist terms $-P_{S,n-3,\epsilon,\gamma_1'}\alpha_1<P_{S,n-1,\alpha_1}^1$. Then $e_{\epsilon\gamma_1\alpha_1}E-e_{\epsilon\gamma_1'\alpha_1}E$ has a minimal $\partial$-invariant completion different from $P_{\epsilon,3,E}$, contradiction.

Thus, we get
\begin{equation*}
\widetilde{P}=P_0+P_1+P_1'=\pm P_{S,n-1,\alpha_1}^1E\pm \sum_jP_{S,n-1,\alpha_2}^jE\in\Omega_n(\Supp(P)),\quad \widetilde{P}\leq P.
\end{equation*}
Furthermore, we can consider the $\partial$-invariant completion of $P_{S,n-1,\alpha_2}^jE$ in $\Supp(\widetilde{P})\subset\Supp(P)$, then we a smaller path
$$\widehat{P}=\pm P_{S,n-1,\alpha_2}^jE\pm P_{S,n-1,\alpha_1}^1E\in\Omega_n(\Supp(\widetilde{P})),\quad \widehat{P}\leq\widetilde{P}\leq P.$$

Then after the similar analysis of $\partial P$ as we do for the case $|E_1|=1$, using the structure result for $P_{S,n-1,\alpha_1}$ and $P_{S,n-1,\alpha_2}$, we finish the case of $|E_1|=2$.

\subsubsection{Induction on $|E_1|$}
Now we consider one more induction on $|E_1|$.
\begin{description}
  \item[\textbf{Initial result 3}] We have proven the result for the case that ``$P$ is of length $n$ and $|E_1|=1,2$" under Hypothesis 1.
  \item[\textbf{Hypothesis 3}] Assume that the result holds for $P$ of length $n$ with $|E_1|\leq m$.
\end{description}

For $|E_1|=m+1$, $E_1=\{\alpha_1,\ldots,\alpha_{m+1}\}$, first, as we have done for the case $|E_1|=2$, we start from $P_0=\pm P_{S,n-1,\alpha_1}E$, and obtain
$$\widetilde{P}=P_0+P_1+P_1'+\cdots+P_l+P_l'\in\Omega_{n}(\Supp(P)).$$

Similarly let us consider the component $P_{S,n-3,\epsilon}\gamma_1\alpha_1E$ of $P_0$, where $\gamma_1\nrightarrow E$. Then we repeat the argument as the case $|E_1|=2$: if the modification does not stop before $P_2$. Then $e_{\epsilon\gamma_1\alpha_1E}$ admits two different minimal $\partial$-invariant completions, contradiction. Then we are done.

\bigskip

For the convenience of the discussion in the next section, we study the intersections of some sub-digraphs in $\Supp(P)$. The result follows from Theorem \ref{structurethm} directly.
\begin{corollary}\label{intersection}
For any $v_k\in d_E^{-1}(n-k)$, $v_l\in d_E^{-1}(n-l)$, $v_k\neq v_l$, let $P_{S,k,v_k}\in\Omega_{k}(\Supp(P))$ and $P_{S,l,v_l}\in\Omega_{l}(\Supp(P))$ be the corresponding minimal paths, then the intersection
$$\Supp(P_{S,k,v_k})\cap \Supp(P_{S,l,v_l})$$
is a union of minimal paths of length less than $\min\{k,l\}$. More precisely, if $v_i\in\Supp(P_{S,k,v_k})\cap \Supp(P_{S,l,v_l})\cap d_S^{-1}(i)$, then
$$\Supp(P_{S,i,v_i})\subset\Supp(P_{S,k,v_k})\cap \Supp(P_{S,l,v_l}).$$
\end{corollary}

One can understand the above structure results in the following subsection of examples of minimal $3$-paths.

\subsection{Minimal $3$-paths and some properties}\label{exminmal3}

In this subsection, we will study the minimal $3$-paths as well as their basic homotopic and homological properties. We will see that unlike the length $2$ case, there are infinitely many kinds of minimal $3$-paths. In particular, we will see that some of the minimal paths are contractible, while some of them are not.

\subsubsection{Homotopy theory of digraphs}\label{Homotopy}

Before the discussion of the examples, let us recall necessary homotopy theory of digraphs developed by Grigoryan-Lin-Muranov-Yau \cite{GLMY}.

\begin{definition}[Digraph map] A morphism from a digraph $G=(V(G),E(G))$ to a digraph $H=(V(H),E(H))$ is a map $f:V(G)\rightarrow V(H)$ such that for any edge $v\rightarrow w$ on $G$, we have
$$f(v)\rightarrow f(w)\quad\text{or}\quad f(v)=f(w) \text{ in }H.$$
The requirement is also denoted by $f(v)\overset{\rightarrow}{=}f(w)$. We refer to such morphisms also as digraph maps, and also denote by $f:G\rightarrow H$.
\end{definition}

\begin{definition}[Cartesian product] For two digraphs $G=(V(G), E(G))$ and $H=(V(H), E(H))$, we define the Cartesian product $G\boxdot H$ as a digraph with the set of vertices $V(G)\times V(H)$ and with the set of edges as follows: for $v,v'\in V(G)$ and $w,w'\in V(H)$, we have
$(v,w)\rightarrow (v',w')$ in $G\boxdot H$ if and only if
$$\text{either $\{v=v',~w\rightarrow w'\}$, or $\{v\rightarrow v',~w=w'\}$}.$$
\end{definition}

Fix $n\geq0$. Denote by $I_n$ any digraph, with $V(I_n)=\{0,1,\ldots,n\}$ and $E(I_n)$ containing exactly one of the edges $i\rightarrow i+1$, $i+1\rightarrow i$ for any
$i=0,1,\ldots, n-1$. We call it the line digraph.

\begin{definition} Let $G$, $H$ be two digraphs. Two digraph maps $f,g:G\rightarrow H$
are called \textbf{homotopic} if there exists a line digraph $I_n$ for some $n\geq1$ and a digraph map $F:G\boxdot I_n\rightarrow H$, such that
$$F\big|_{G\boxdot\{0\}}=f\quad\text{and}\quad F\big|_{G\boxdot\{n\}}=g.$$
We shall write $f\simeq g$ and call $F$ a $n$-step homotopy between $f$ and $g$.
\end{definition}

\begin{definition} Two digraphs $G$ and $H$ are called \textbf{homotopy equivalent} if there exist digraph maps
$$f:G\rightarrow H,\quad g:H\rightarrow G$$
such that
$$f\circ g\simeq\id_H,\quad g\circ f\simeq\id_G.$$
A digraph $G$ is called contractible if $G\simeq\{*\}$ where $\{*\}$ is a single vertex digraph.
\end{definition}

Similarly, two homotopy equivalent digraphs have the isomorphic path homology groups, see in \cite{GLMY}. In particular, a contractible digraph has acyclic path homologies.

\begin{definition} Let $G$ be a digraph and $H$ be its sub-digraph.
\begin{enumerate}
  \item A \textbf{retraction} of $G$ onto $H$ is a digraph map $r:G\rightarrow H$ such that $r\big|_H=\id_H$.
  \item A retraction $r:G\rightarrow H$ is called a \textbf{deformation retraction} if $i\circ r\simeq\id_G$, where $i:H\rightarrow G$ is the natural inclusion digraph map.
\end{enumerate}
\end{definition}

\begin{proposition}[\cite{GLMY}]\label{nstepr} Let $r: G\rightarrow H$ be a retraction of a digraph $G$ onto a sub-digraph $H$. Assume that there exists a finite sequence $\{f_k\}_{k=0}^n$ of digraph maps $f_k: G\rightarrow G$ with the following properties
\begin{enumerate}
  \item $f_0=\id_G$, $f_n=i\circ r$;
  \item for any $k=1,\ldots,n$ either $f_{k-1}(x)\overset{\rightarrow}{=}f_k(x)$, any $x\in V(G)$, or $f_{k-1}(x)\overset{\leftarrow}{=}f_k(x)$, any $x\in V(G)$.
\end{enumerate}
Then $r$ is a deformation retraction, the digraphs $G$ and $H$ are homotopy equivalent.
\end{proposition}

\begin{corollary}[\cite{GLMY}]\label{onestepr} Let $r: G\rightarrow H$ be a retraction of a digraph $G$ onto a sub-digraph $H$ and
$$r(x)\overset{\rightarrow}{=}x\text{ for all }x\in V(G),\quad \text{or}\quad x\overset{\rightarrow}{=}r(x)\text{ for all }x\in V(G).$$
Then $r$ is a deformation retraction, the digraphs $G$ and $H$ are homotopy equivalent.
\end{corollary}

\subsubsection{Examples: minimal $3$-paths}

\begin{example} If an elementary path $P=e_{0123}$ is $\partial$-invariant, then the corresponding supporting digraph is as follows.
\begin{figure}[H]
	\centering
	%% Creator: Inkscape 1.0.1 (3bc2e813f5, 2020-09-07), www.inkscape.org
%% PDF/EPS/PS + LaTeX output extension by Johan Engelen, 2010
%% Accompanies image file 'EX313.pdf' (pdf, eps, ps)
%%
%% To include the image in your LaTeX document, write
%%   \input{<filename>.pdf_tex}
%%  instead of
%%   \includegraphics{<filename>.pdf}
%% To scale the image, write
%%   \def\svgwidth{<desired width>}
%%   \input{<filename>.pdf_tex}
%%  instead of
%%   \includegraphics[width=<desired width>]{<filename>.pdf}
%%
%% Images with a different path to the parent latex file can
%% be accessed with the `import' package (which may need to be
%% installed) using
%%   \usepackage{import}
%% in the preamble, and then including the image with
%%   \import{<path to file>}{<filename>.pdf_tex}
%% Alternatively, one can specify
%%   \graphicspath{{<path to file>/}}
%% 
%% For more information, please see info/svg-inkscape on CTAN:
%%   http://tug.ctan.org/tex-archive/info/svg-inkscape
%%
\begingroup%
  \makeatletter%
  \providecommand\color[2][]{%
    \errmessage{(Inkscape) Color is used for the text in Inkscape, but the package 'color.sty' is not loaded}%
    \renewcommand\color[2][]{}%
  }%
  \providecommand\transparent[1]{%
    \errmessage{(Inkscape) Transparency is used (non-zero) for the text in Inkscape, but the package 'transparent.sty' is not loaded}%
    \renewcommand\transparent[1]{}%
  }%
  \providecommand\rotatebox[2]{#2}%
  \newcommand*\fsize{\dimexpr\f@size pt\relax}%
  \newcommand*\lineheight[1]{\fontsize{\fsize}{#1\fsize}\selectfont}%
  \ifx\svgwidth\undefined%
    \setlength{\unitlength}{195.03154604bp}%
    \ifx\svgscale\undefined%
      \relax%
    \else%
      \setlength{\unitlength}{\unitlength * \real{\svgscale}}%
    \fi%
  \else%
    \setlength{\unitlength}{\svgwidth}%
  \fi%
  \global\let\svgwidth\undefined%
  \global\let\svgscale\undefined%
  \makeatother%
  \begin{picture}(1,0.41334657)%
    \lineheight{1}%
    \setlength\tabcolsep{0pt}%
    \put(-0.00290714,0.192108){\makebox(0,0)[lt]{\lineheight{1.25}\smash{\begin{tabular}[t]{l}$\Supp(P)=$\end{tabular}}}}%
    \put(0,0){\includegraphics[width=\unitlength,page=1]{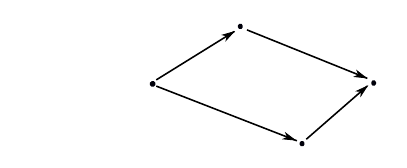}}%
    \put(0.29319116,0.18980587){\makebox(0,0)[lt]{\lineheight{1.25}\smash{\begin{tabular}[t]{l}$0$\end{tabular}}}}%
    \put(0.56168551,0.3794491){\makebox(0,0)[lt]{\lineheight{1.25}\smash{\begin{tabular}[t]{l}$1$\end{tabular}}}}%
    \put(0.72499958,0.00279948){\makebox(0,0)[lt]{\lineheight{1.25}\smash{\begin{tabular}[t]{l}$2$\end{tabular}}}}%
    \put(0.93950045,0.18672258){\makebox(0,0)[lt]{\lineheight{1.25}\smash{\begin{tabular}[t]{l}$3$\end{tabular}}}}%
    \put(0,0){\includegraphics[width=\unitlength,page=2]{EX313.pdf}}%
  \end{picture}%
\endgroup%

	%\caption[]{$\Supp(e_{0123})$}
	%\label{Fig:Q1}
\end{figure}
\vskip -0.2cm

We can compute $\partial e_{0123}$ and present it in terms of the minimal $2$-paths in $\Supp(P)$,
\begin{align*}
\partial e_{0123}=-e_{012}+(e_{013}-e_{023})+e_{123}.
\end{align*}
Let $r$ be the digraph map
\begin{figure}[H]
	\centering
	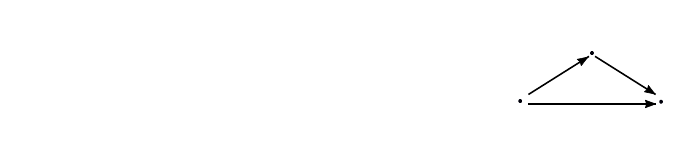
	%\caption[]{$\Supp(P)$}
	%\label{Fig:Q1}
\end{figure}
\noindent given by $r(3)=2$ and $r\big|_H=\id_H$. By Corollary \ref{onestepr}, we get
$$i\circ r\simeq \id_{\Supp(P)}$$
by a one-step homotopy. Furthermore, $H$ is obviously contractible, so is $\Supp(P)$. In particular, $\widetilde{H}_*(\Supp(P))=0$.
\end{example}

\begin{example} The minimal digraph which supports $P=e_{0134}-e_{0234}$ is as follows.
\begin{figure}[H]
	\centering
	%% Creator: Inkscape 1.0.1 (3bc2e813f5, 2020-09-07), www.inkscape.org
%% PDF/EPS/PS + LaTeX output extension by Johan Engelen, 2010
%% Accompanies image file 'EX314.pdf' (pdf, eps, ps)
%%
%% To include the image in your LaTeX document, write
%%   \input{<filename>.pdf_tex}
%%  instead of
%%   \includegraphics{<filename>.pdf}
%% To scale the image, write
%%   \def\svgwidth{<desired width>}
%%   \input{<filename>.pdf_tex}
%%  instead of
%%   \includegraphics[width=<desired width>]{<filename>.pdf}
%%
%% Images with a different path to the parent latex file can
%% be accessed with the `import' package (which may need to be
%% installed) using
%%   \usepackage{import}
%% in the preamble, and then including the image with
%%   \import{<path to file>}{<filename>.pdf_tex}
%% Alternatively, one can specify
%%   \graphicspath{{<path to file>/}}
%% 
%% For more information, please see info/svg-inkscape on CTAN:
%%   http://tug.ctan.org/tex-archive/info/svg-inkscape
%%
\begingroup%
  \makeatletter%
  \providecommand\color[2][]{%
    \errmessage{(Inkscape) Color is used for the text in Inkscape, but the package 'color.sty' is not loaded}%
    \renewcommand\color[2][]{}%
  }%
  \providecommand\transparent[1]{%
    \errmessage{(Inkscape) Transparency is used (non-zero) for the text in Inkscape, but the package 'transparent.sty' is not loaded}%
    \renewcommand\transparent[1]{}%
  }%
  \providecommand\rotatebox[2]{#2}%
  \newcommand*\fsize{\dimexpr\f@size pt\relax}%
  \newcommand*\lineheight[1]{\fontsize{\fsize}{#1\fsize}\selectfont}%
  \ifx\svgwidth\undefined%
    \setlength{\unitlength}{228.10694657bp}%
    \ifx\svgscale\undefined%
      \relax%
    \else%
      \setlength{\unitlength}{\unitlength * \real{\svgscale}}%
    \fi%
  \else%
    \setlength{\unitlength}{\svgwidth}%
  \fi%
  \global\let\svgwidth\undefined%
  \global\let\svgscale\undefined%
  \makeatother%
  \begin{picture}(1,0.26824515)%
    \lineheight{1}%
    \setlength\tabcolsep{0pt}%
    \put(0,0){\includegraphics[width=\unitlength,page=1]{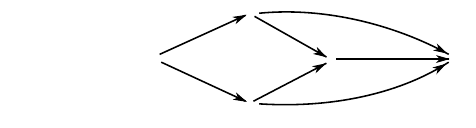}}%
    \put(0.50761068,0.25483439){\makebox(0,0)[lt]{\lineheight{1.25}\smash{\begin{tabular}[t]{l}$1$\end{tabular}}}}%
    \put(0.51475613,0.00199092){\makebox(0,0)[lt]{\lineheight{1.25}\smash{\begin{tabular}[t]{l}$2$\end{tabular}}}}%
    \put(0.61135128,0.13557957){\makebox(0,0)[lt]{\lineheight{1.25}\smash{\begin{tabular}[t]{l}$3$\end{tabular}}}}%
    \put(0.97073908,0.1388289){\makebox(0,0)[lt]{\lineheight{1.25}\smash{\begin{tabular}[t]{l}$4$\end{tabular}}}}%
    \put(0,0){\includegraphics[width=\unitlength,page=2]{EX314.pdf}}%
    \put(-0.00177002,0.13197324){\makebox(0,0)[lt]{\lineheight{1.25}\smash{\begin{tabular}[t]{l}$\Supp(P)=~ 0$\end{tabular}}}}%
  \end{picture}%
\endgroup%

	%\caption[]{$\Supp(P)$}
	%\label{Fig:Q2}
\end{figure}

Similarly, we have the decomposition
\begin{align*}
\partial (e_{0134}-e_{0234})=-(e_{013}-e_{023})+e_{134}-e_{234}+(e_{014}-e_{024}),
\end{align*}
and the deformation retraction
\begin{figure}[H]
	\centering
	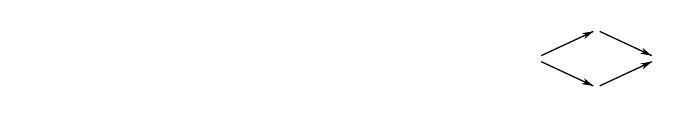
	%\caption[]{$\Supp(e_{0123})$}
	%\label{Fig:Q1}
\end{figure}
\noindent which is given by $r(4)=3$ and $r\big|_H=\id_H$. Thus, $\Supp(P)$ is contractible.
\end{example}

\begin{example} The minimal digraph which supports $P=e_{0135}-e_{0235}+e_{0245}$ is as follows.
\begin{figure}[H]
	\centering
	%% Creator: Inkscape 1.0.1 (3bc2e813f5, 2020-09-07), www.inkscape.org
%% PDF/EPS/PS + LaTeX output extension by Johan Engelen, 2010
%% Accompanies image file 'EX315.pdf' (pdf, eps, ps)
%%
%% To include the image in your LaTeX document, write
%%   \input{<filename>.pdf_tex}
%%  instead of
%%   \includegraphics{<filename>.pdf}
%% To scale the image, write
%%   \def\svgwidth{<desired width>}
%%   \input{<filename>.pdf_tex}
%%  instead of
%%   \includegraphics[width=<desired width>]{<filename>.pdf}
%%
%% Images with a different path to the parent latex file can
%% be accessed with the `import' package (which may need to be
%% installed) using
%%   \usepackage{import}
%% in the preamble, and then including the image with
%%   \import{<path to file>}{<filename>.pdf_tex}
%% Alternatively, one can specify
%%   \graphicspath{{<path to file>/}}
%% 
%% For more information, please see info/svg-inkscape on CTAN:
%%   http://tug.ctan.org/tex-archive/info/svg-inkscape
%%
\begingroup%
  \makeatletter%
  \providecommand\color[2][]{%
    \errmessage{(Inkscape) Color is used for the text in Inkscape, but the package 'color.sty' is not loaded}%
    \renewcommand\color[2][]{}%
  }%
  \providecommand\transparent[1]{%
    \errmessage{(Inkscape) Transparency is used (non-zero) for the text in Inkscape, but the package 'transparent.sty' is not loaded}%
    \renewcommand\transparent[1]{}%
  }%
  \providecommand\rotatebox[2]{#2}%
  \newcommand*\fsize{\dimexpr\f@size pt\relax}%
  \newcommand*\lineheight[1]{\fontsize{\fsize}{#1\fsize}\selectfont}%
  \ifx\svgwidth\undefined%
    \setlength{\unitlength}{234.99249796bp}%
    \ifx\svgscale\undefined%
      \relax%
    \else%
      \setlength{\unitlength}{\unitlength * \real{\svgscale}}%
    \fi%
  \else%
    \setlength{\unitlength}{\svgwidth}%
  \fi%
  \global\let\svgwidth\undefined%
  \global\let\svgscale\undefined%
  \makeatother%
  \begin{picture}(1,0.25034408)%
    \lineheight{1}%
    \setlength\tabcolsep{0pt}%
    \put(0,0){\includegraphics[width=\unitlength,page=1]{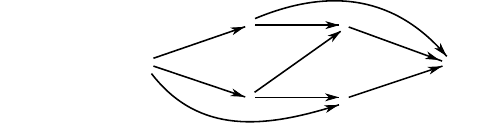}}%
    \put(-0.00161826,0.10272789){\makebox(0,0)[lt]{\lineheight{1.25}\smash{\begin{tabular}[t]{l}$\Supp(P)=~0$\end{tabular}}}}%
    \put(0.45816729,0.21803569){\makebox(0,0)[lt]{\lineheight{1.25}\smash{\begin{tabular}[t]{l}$1$\end{tabular}}}}%
    \put(0.48187275,0.07982951){\makebox(0,0)[lt]{\lineheight{1.25}\smash{\begin{tabular}[t]{l}$2$\end{tabular}}}}%
    \put(0.70293417,0.14779869){\makebox(0,0)[lt]{\lineheight{1.25}\smash{\begin{tabular}[t]{l}$3$\end{tabular}}}}%
    \put(0.70176179,0.0041847){\makebox(0,0)[lt]{\lineheight{1.25}\smash{\begin{tabular}[t]{l}$4$\end{tabular}}}}%
    \put(0.94374172,0.1046775){\makebox(0,0)[lt]{\lineheight{1.25}\smash{\begin{tabular}[t]{l}$5$\end{tabular}}}}%
    \put(0,0){\includegraphics[width=\unitlength,page=2]{EX315.pdf}}%
  \end{picture}%
\endgroup%

	%\caption[]{$\Supp(e_{0123})$}
	%\label{Fig:Q1}
\end{figure}

Similarly, we have the decomposition
\begin{align*}
\partial(e_{0135}-e_{0235}+e_{0245})=-(e_{013}-e_{023})-e_{024}+e_{135}+(e_{245}-e_{235})+(e_{015}-e_{045}).
\end{align*}
One can construct the following deformation retractions.

(1) $r_1(5)=3$, $r_1(4)=2$.
\begin{figure}[H]
	\centering
	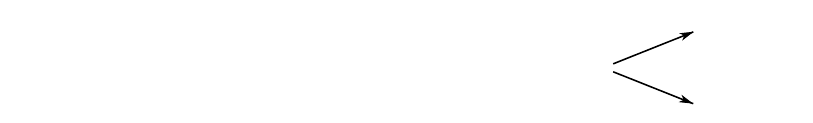
	%\caption[]{$\Supp(e_{0123})$}
	%\label{Fig:Q1}
\end{figure}

(2) $r_2(5)=4$, $r_2(3)=2$, $r_2(1)=0$.
\begin{figure}[H]
	\centering
	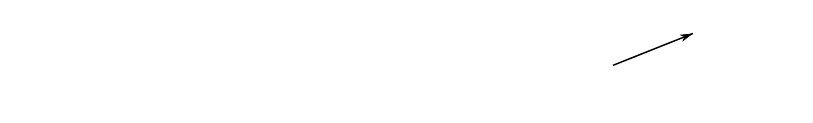
	%\caption[]{$\Supp(e_{0123})$}
	%\label{Fig:Q1}
\end{figure}
\noindent Thus, $\Supp(P)$ is contractible.
\end{example}

\begin{example} The minimal digraph which supports $P=e_{0146}-e_{0136}+e_{0236}-e_{0256}$ is as follows.
\begin{figure}[H]
	\centering
	%% Creator: Inkscape 1.0.1 (3bc2e813f5, 2020-09-07), www.inkscape.org
%% PDF/EPS/PS + LaTeX output extension by Johan Engelen, 2010
%% Accompanies image file 'EX316.pdf' (pdf, eps, ps)
%%
%% To include the image in your LaTeX document, write
%%   \input{<filename>.pdf_tex}
%%  instead of
%%   \includegraphics{<filename>.pdf}
%% To scale the image, write
%%   \def\svgwidth{<desired width>}
%%   \input{<filename>.pdf_tex}
%%  instead of
%%   \includegraphics[width=<desired width>]{<filename>.pdf}
%%
%% Images with a different path to the parent latex file can
%% be accessed with the `import' package (which may need to be
%% installed) using
%%   \usepackage{import}
%% in the preamble, and then including the image with
%%   \import{<path to file>}{<filename>.pdf_tex}
%% Alternatively, one can specify
%%   \graphicspath{{<path to file>/}}
%% 
%% For more information, please see info/svg-inkscape on CTAN:
%%   http://tug.ctan.org/tex-archive/info/svg-inkscape
%%
\begingroup%
  \makeatletter%
  \providecommand\color[2][]{%
    \errmessage{(Inkscape) Color is used for the text in Inkscape, but the package 'color.sty' is not loaded}%
    \renewcommand\color[2][]{}%
  }%
  \providecommand\transparent[1]{%
    \errmessage{(Inkscape) Transparency is used (non-zero) for the text in Inkscape, but the package 'transparent.sty' is not loaded}%
    \renewcommand\transparent[1]{}%
  }%
  \providecommand\rotatebox[2]{#2}%
  \newcommand*\fsize{\dimexpr\f@size pt\relax}%
  \newcommand*\lineheight[1]{\fontsize{\fsize}{#1\fsize}\selectfont}%
  \ifx\svgwidth\undefined%
    \setlength{\unitlength}{238.69811291bp}%
    \ifx\svgscale\undefined%
      \relax%
    \else%
      \setlength{\unitlength}{\unitlength * \real{\svgscale}}%
    \fi%
  \else%
    \setlength{\unitlength}{\svgwidth}%
  \fi%
  \global\let\svgwidth\undefined%
  \global\let\svgscale\undefined%
  \makeatother%
  \begin{picture}(1,0.30445378)%
    \lineheight{1}%
    \setlength\tabcolsep{0pt}%
    \put(0,0){\includegraphics[width=\unitlength,page=1]{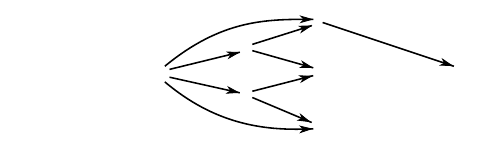}}%
    \put(0.65950576,0.14901177){\makebox(0,0)[lt]{\lineheight{1.25}\smash{\begin{tabular}[t]{l}$3$\end{tabular}}}}%
    \put(0,0){\includegraphics[width=\unitlength,page=2]{EX316.pdf}}%
    \put(0.47364739,0.15894658){\makebox(0,0)[lt]{\lineheight{1.25}\smash{\begin{tabular}[t]{l}$1$\end{tabular}}}}%
    \put(0.48211296,0.07201179){\makebox(0,0)[lt]{\lineheight{1.25}\smash{\begin{tabular}[t]{l}$2$\end{tabular}}}}%
    \put(0.62258352,0.28200176){\makebox(0,0)[lt]{\lineheight{1.25}\smash{\begin{tabular}[t]{l}$4$\end{tabular}}}}%
    \put(0.62246847,0.0029767){\makebox(0,0)[lt]{\lineheight{1.25}\smash{\begin{tabular}[t]{l}$5$\end{tabular}}}}%
    \put(0.95401583,0.15142087){\makebox(0,0)[lt]{\lineheight{1.25}\smash{\begin{tabular}[t]{l}$6$\end{tabular}}}}%
    \put(0,0){\includegraphics[width=\unitlength,page=3]{EX316.pdf}}%
    \put(-0.00293776,0.14148126){\makebox(0,0)[lt]{\lineheight{1.25}\smash{\begin{tabular}[t]{l}$\Supp(P)=~0$\end{tabular}}}}%
  \end{picture}%
\endgroup%

	%\caption[]{$\Supp(e_{0123})$}
	%\label{Fig:Q1}
\end{figure}

Similarly, we have
\begin{align*}
&~\partial(e_{0146}-e_{0136}+e_{0236}-e_{0256})\\
=~&-e_{014}+(e_{013}-e_{023})+e_{025}+(e_{146}-e_{136})+(e_{236}-e_{256})-(e_{046}-e_{056}).
\end{align*}
One can construct the following deformation retractions.

(1) $r_1(6)=3$, $r_1(5)=2$, $r_1(4)=1$.
\begin{figure}[H]
	\centering
	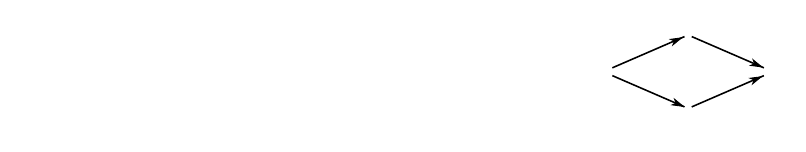
	%\caption[]{$\Supp(e_{0123})$}
	%\label{Fig:Q1}
\end{figure}

(2) $r_2(6)=4$, $r_2(3)=1$, $r_2(5)=r_2(2)=0$.
\begin{figure}[H]
	\centering
	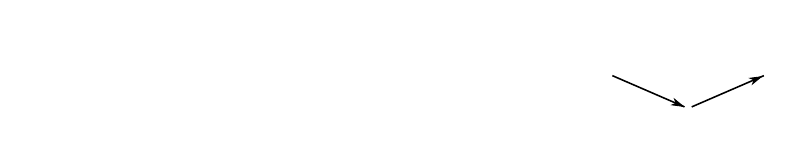
	%\caption[]{$\Supp(e_{0123})$}
	%\label{Fig:Q1}
\end{figure}

(3) $r_3(1)=4$, $r_3(2)=3$, $r_3(3)=6$.
\begin{figure}[H]
	\centering
	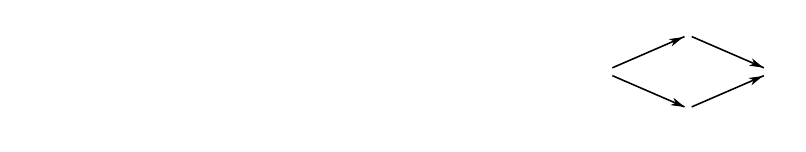
	%\caption[]{$\Supp(e_{0123})$}
	%\label{Fig:Q1}
\end{figure}

(4) $r_4(0)=1$, $r_4(2)=3$, $r_4(5)=3$.
\begin{figure}[H]
	\centering
	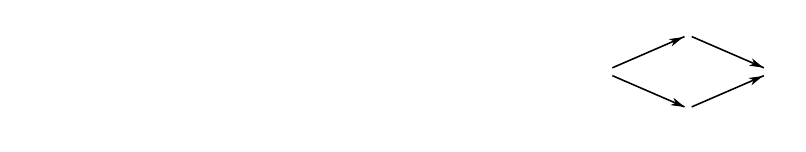
	%\caption[]{$\Supp(e_{0123})$}
	%\label{Fig:Q1}
\end{figure}

The relation $\id_{\Supp(P)}\simeq i_4\circ r_4$ can be realized by a $3$-step homotopy
$$F: \Supp(P)\boxdot I_3\rightarrow \Supp(P),\quad\text{with}\quad I_3= 0\leftarrow 1\rightarrow 2\leftarrow 3$$
with is defined to be
\begin{itemize}
  \item $F(\cdot,0)=\id_{\Supp(P)}$;
  \item $F(0,1)=0, F(1,1)=1, F(2,1)=2, F(3,1)=3$,\\
        $F(4,1)=1, F(5,1)=2, F(6,1)=3$.\qquad\qquad   (We have $F(x,0)\overset{\leftarrow}{=} F(x,1)$.)
  \item $F(0,2)=1, F(1,2)=4, F(2,2)=3, F(3,2)=6$,\\
        $F(4,2)=4, F(5,2)=3, F(6,2)=6$.\qquad\qquad  (We have $F(x,1)\overset{\rightarrow}{=} F(x,2)$.)
  \item $F(v,3)=i_4\circ r_4$, that is\\
        $F(0,3)=1, F(1,3)=1, F(2,3)=3, F(3,3)=3$,\\
        $F(4,3)=4, F(5,3)=3, F(6,3)=6$.\qquad\qquad (We have $F(x,2)\overset{\leftarrow}{=} F(x,3)$.)
\end{itemize}
Thus $\Supp(P)$ is contractible via alternative homotopies.
\end{example}

Note that in all the above examples of minimal $3$-paths, there is exactly one minimal $2$-path of type $P_{S,2,E}$. The following 3 examples tell us there may be no such a term.

\begin{example}\label{012345} The minimal digraph which supports $P=e_{0135}-e_{0145}+e_{0245}-e_{0235}$ is as follows.
\begin{figure}[H]
	\centering
	%% Creator: Inkscape 1.0.1 (3bc2e813f5, 2020-09-07), www.inkscape.org
%% PDF/EPS/PS + LaTeX output extension by Johan Engelen, 2010
%% Accompanies image file 'EX012345.pdf' (pdf, eps, ps)
%%
%% To include the image in your LaTeX document, write
%%   \input{<filename>.pdf_tex}
%%  instead of
%%   \includegraphics{<filename>.pdf}
%% To scale the image, write
%%   \def\svgwidth{<desired width>}
%%   \input{<filename>.pdf_tex}
%%  instead of
%%   \includegraphics[width=<desired width>]{<filename>.pdf}
%%
%% Images with a different path to the parent latex file can
%% be accessed with the `import' package (which may need to be
%% installed) using
%%   \usepackage{import}
%% in the preamble, and then including the image with
%%   \import{<path to file>}{<filename>.pdf_tex}
%% Alternatively, one can specify
%%   \graphicspath{{<path to file>/}}
%% 
%% For more information, please see info/svg-inkscape on CTAN:
%%   http://tug.ctan.org/tex-archive/info/svg-inkscape
%%
\begingroup%
  \makeatletter%
  \providecommand\color[2][]{%
    \errmessage{(Inkscape) Color is used for the text in Inkscape, but the package 'color.sty' is not loaded}%
    \renewcommand\color[2][]{}%
  }%
  \providecommand\transparent[1]{%
    \errmessage{(Inkscape) Transparency is used (non-zero) for the text in Inkscape, but the package 'transparent.sty' is not loaded}%
    \renewcommand\transparent[1]{}%
  }%
  \providecommand\rotatebox[2]{#2}%
  \newcommand*\fsize{\dimexpr\f@size pt\relax}%
  \newcommand*\lineheight[1]{\fontsize{\fsize}{#1\fsize}\selectfont}%
  \ifx\svgwidth\undefined%
    \setlength{\unitlength}{223.10847005bp}%
    \ifx\svgscale\undefined%
      \relax%
    \else%
      \setlength{\unitlength}{\unitlength * \real{\svgscale}}%
    \fi%
  \else%
    \setlength{\unitlength}{\svgwidth}%
  \fi%
  \global\let\svgwidth\undefined%
  \global\let\svgscale\undefined%
  \makeatother%
  \begin{picture}(1,0.31815132)%
    \lineheight{1}%
    \setlength\tabcolsep{0pt}%
    \put(0,0){\includegraphics[width=\unitlength,page=1]{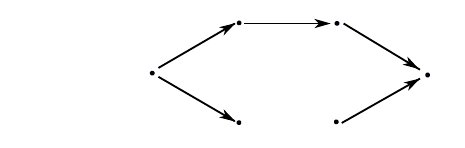}}%
    \put(-0.0021868,0.14789372){\makebox(0,0)[lt]{\lineheight{1.25}\smash{\begin{tabular}[t]{l}$\Supp(P)=$\end{tabular}}}}%
    \put(0.26636067,0.1485136){\makebox(0,0)[lt]{\lineheight{1.25}\smash{\begin{tabular}[t]{l}$0$\end{tabular}}}}%
    \put(0.48668306,0.2906922){\makebox(0,0)[lt]{\lineheight{1.25}\smash{\begin{tabular}[t]{l}$1$\end{tabular}}}}%
    \put(0.7011585,0.29138943){\makebox(0,0)[lt]{\lineheight{1.25}\smash{\begin{tabular}[t]{l}$3$\end{tabular}}}}%
    \put(0.94783225,0.14375306){\makebox(0,0)[lt]{\lineheight{1.25}\smash{\begin{tabular}[t]{l}$5$\end{tabular}}}}%
    \put(0,0){\includegraphics[width=\unitlength,page=2]{EX012345.pdf}}%
    \put(0.48417832,0.00346404){\makebox(0,0)[lt]{\lineheight{1.25}\smash{\begin{tabular}[t]{l}$2$\end{tabular}}}}%
    \put(0.69946788,0.00511006){\makebox(0,0)[lt]{\lineheight{1.25}\smash{\begin{tabular}[t]{l}$4$\end{tabular}}}}%
    \put(0,0){\includegraphics[width=\unitlength,page=3]{EX012345.pdf}}%
  \end{picture}%
\endgroup%

	%\caption[]{$\Supp(e_{0123})$}
	%\label{Fig:Q1}
\end{figure}
Similarly, we have
\begin{align*}
&\partial(e_{0135}-e_{0145}+e_{0245}-e_{0235})\\
=~&(e_{023}-e_{013})+(e_{014}-e_{024})+(e_{135}-e_{145})+(e_{245}-e_{235}).
\end{align*}
There is a deformation retraction as follows:
\begin{figure}[H]
	\centering
	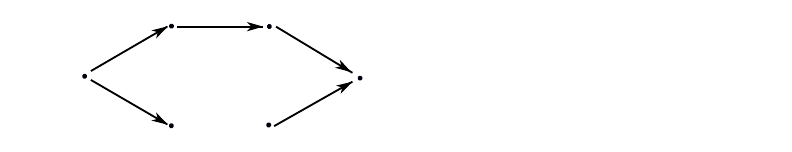
	%\caption[]{$\Supp(e_{0123})$}
	%\label{Fig:Q1}
\end{figure}

\noindent $r(5)=r(4)=3$, $r\big|_H=\id_H$,  where the homotopy
$$F: \Supp(P)\boxdot I_3\rightarrow \Supp(P),\quad\text{with}\quad I_3= 0\leftarrow 1\leftarrow 2\rightarrow 3$$
connecting $\id_{\Supp(P)}$ and $i\circ r$ is given by
\begin{itemize}
  \item $F(\cdot,0)=\id_{\Supp(P)}$;
  \item $F(0,1)=0, F(1,1)=1, F(2,1)=0$, \\
        $F(3,1)=1, F(4,1)=1, F(5,1)=3$.\qquad  $F(x,0)\overset{\leftarrow}{=} F(x,1)$.
        %(It's easy to check that $F(x,0)\xleftarrow{=} F(x,1)$);
  \item $F(0,2)=0, F(1,2)=1, F(2,2)=0$, \\
        $F(3,2)=1, F(4,2)=1, F(5,2)=1$.\qquad $F(x,1)\overset{\leftarrow}{=} F(x,2)$.
        %(It's easy to check that $F(x,1)\xleftarrow{=} F(x,2)$);
  \item $F(v,3)=i\circ r$, that is\\
        $F(0,3)=1, F(1,3)=1, F(2,3)=2$,\\
        $F(3,3)=3, F(4,3)=3, F(5,3)=3$.\qquad $F(x,2)\overset{\rightarrow}{=} F(x,3)$.
        %(It's easy to check that $F(x,2)\xrightarrow{=} F(x,3)$).
\end{itemize}

Note that one can also define $F(v,3)=0$ for all $v$ and then we see that the map $\Supp{P}\rightarrow\{0\}$ is a deformation retraction to $\{0\}$.
\end{example}

\begin{example}\label{EXfSfE1} The minimal digraph which supports
$$P=e_{0136}-e_{0156}+e_{0456}+e_{0246}-e_{0236}$$
is as follows.
\begin{figure}[H]
	\centering
	%% Creator: Inkscape 1.0.1 (3bc2e813f5, 2020-09-07), www.inkscape.org
%% PDF/EPS/PS + LaTeX output extension by Johan Engelen, 2010
%% Accompanies image file 'EXfSfE1.pdf' (pdf, eps, ps)
%%
%% To include the image in your LaTeX document, write
%%   \input{<filename>.pdf_tex}
%%  instead of
%%   \includegraphics{<filename>.pdf}
%% To scale the image, write
%%   \def\svgwidth{<desired width>}
%%   \input{<filename>.pdf_tex}
%%  instead of
%%   \includegraphics[width=<desired width>]{<filename>.pdf}
%%
%% Images with a different path to the parent latex file can
%% be accessed with the `import' package (which may need to be
%% installed) using
%%   \usepackage{import}
%% in the preamble, and then including the image with
%%   \import{<path to file>}{<filename>.pdf_tex}
%% Alternatively, one can specify
%%   \graphicspath{{<path to file>/}}
%% 
%% For more information, please see info/svg-inkscape on CTAN:
%%   http://tug.ctan.org/tex-archive/info/svg-inkscape
%%
\begingroup%
  \makeatletter%
  \providecommand\color[2][]{%
    \errmessage{(Inkscape) Color is used for the text in Inkscape, but the package 'color.sty' is not loaded}%
    \renewcommand\color[2][]{}%
  }%
  \providecommand\transparent[1]{%
    \errmessage{(Inkscape) Transparency is used (non-zero) for the text in Inkscape, but the package 'transparent.sty' is not loaded}%
    \renewcommand\transparent[1]{}%
  }%
  \providecommand\rotatebox[2]{#2}%
  \newcommand*\fsize{\dimexpr\f@size pt\relax}%
  \newcommand*\lineheight[1]{\fontsize{\fsize}{#1\fsize}\selectfont}%
  \ifx\svgwidth\undefined%
    \setlength{\unitlength}{204.63798703bp}%
    \ifx\svgscale\undefined%
      \relax%
    \else%
      \setlength{\unitlength}{\unitlength * \real{\svgscale}}%
    \fi%
  \else%
    \setlength{\unitlength}{\svgwidth}%
  \fi%
  \global\let\svgwidth\undefined%
  \global\let\svgscale\undefined%
  \makeatother%
  \begin{picture}(1,0.39427421)%
    \lineheight{1}%
    \setlength\tabcolsep{0pt}%
    \put(0,0){\includegraphics[width=\unitlength,page=1]{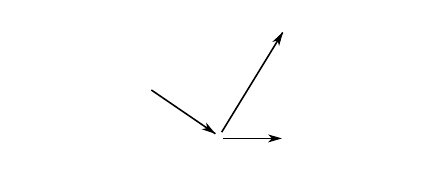}}%
    \put(0.29904814,0.17541918){\makebox(0,0)[lt]{\lineheight{1.25}\smash{\begin{tabular}[t]{l}$0$\end{tabular}}}}%
    \put(0.49887996,0.36225333){\makebox(0,0)[lt]{\lineheight{1.25}\smash{\begin{tabular}[t]{l}$1$\end{tabular}}}}%
    \put(0.49484763,0.00386232){\makebox(0,0)[lt]{\lineheight{1.25}\smash{\begin{tabular}[t]{l}$2$\end{tabular}}}}%
    \put(0.65352713,0.36571285){\makebox(0,0)[lt]{\lineheight{1.25}\smash{\begin{tabular}[t]{l}$3$\end{tabular}}}}%
    \put(0.65976945,0.00285283){\makebox(0,0)[lt]{\lineheight{1.25}\smash{\begin{tabular}[t]{l}$4$\end{tabular}}}}%
    \put(0.7676047,0.13021959){\makebox(0,0)[lt]{\lineheight{1.25}\smash{\begin{tabular}[t]{l}$5$\end{tabular}}}}%
    \put(0,0){\includegraphics[width=\unitlength,page=2]{EXfSfE1.pdf}}%
    \put(-0.00209816,0.17523738){\makebox(0,0)[lt]{\lineheight{1.25}\smash{\begin{tabular}[t]{l}$\Supp(P)=$\end{tabular}}}}%
    \put(0,0){\includegraphics[width=\unitlength,page=3]{EXfSfE1.pdf}}%
    \put(0.95114656,0.18653845){\makebox(0,0)[lt]{\lineheight{1.25}\smash{\begin{tabular}[t]{l}$6$\end{tabular}}}}%
  \end{picture}%
\endgroup%

	%\caption[]{$\Supp(e_{0123})$}
	%\label{Fig:Q1}
\end{figure}
Similarly, we have
\begin{align*}
&\partial(e_{0136}-e_{0156}+e_{0456}+e_{0246}-e_{0236})\\
=~&(e_{023}-e_{013})-e_{024}+(e_{015}-e_{045})+(e_{136}-e_{156})+(e_{246}-e_{236})+e_{456}.
\end{align*}
One can construct the following deformation retractions.

(1) $r_1(6)=5$, $r_1(3)=1$, $r_1(2)=0$.
\begin{figure}[H]
	\centering
	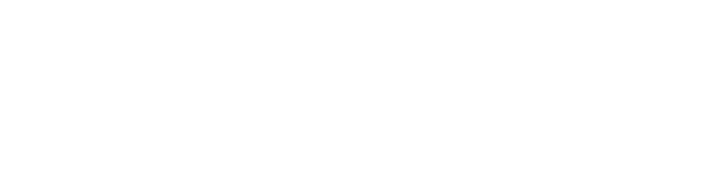
	%\caption[]{$\Supp(e_{0123})$}
	%\label{Fig:Q1}
\end{figure}

(2) $r_2(6)=r_2(5)=4$, $r_2(3)=2$, $r_2(1)=0$.
\begin{figure}[H]
	\centering
	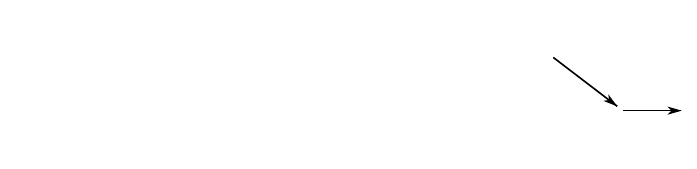
	%\caption[]{$\Supp(e_{0123})$}
	%\label{Fig:Q1}
\end{figure}

\end{example}

\begin{example}[$3$-Cube]\label{Cube}  The minimal digraph which supports
$$P=e_{0137}-e_{0237}+e_{0267}-e_{0467}+e_{0457}-e_{0157}$$
is as follows.
\begin{figure}[H]
	\centering
	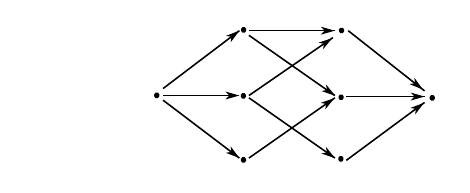
	%\caption[]{$\Supp(e_{0123})$}
	%\label{Fig:Q1}
\end{figure}
Similarly, we have
\begin{align*}
&\partial(e_{0137}-e_{0237}+e_{0267}-e_{0467}+e_{0457}-e_{0157})\\
=~&(e_{023}-e_{013})+(e_{015}-e_{045})+(e_{046}-e_{026})\\
&+(e_{137}-e_{157})-(e_{237}-e_{267})+(e_{457}-e_{467}).
\end{align*}
It is known as the $3$-cube digraph, which is isomorphic (there exists a bijective digraph map) to $I_1^{\boxdot3}$. Thus it is contractible. Also note that there are not minimal $2$-paths of type $P_{S,2,E}$ in this example as well as in Examples \ref{012345}, \ref{EXfSfE1}.
\end{example}

\begin{remark} We learn from the above examples that all of them can deform retract to some (actually all, one can construct one by one) minimal faces. By homotopy invariance, they have acyclic path homologies. In fact, one can also construct the digraph map $h: I_1^{\boxdot3}\rightarrow \Supp(P)$ such that the above $\Supp(P)$ is exactly the image of $h$. One can regard such a map as the analogue of the characteristic map in the CW complex of a topological space. We will do the further discussion in a separation paper \cite{TY2}.
\end{remark}

Now let us continue looking at two different minimal $3$-paths.

\begin{example}\label{firstnoncontractible} The minimal digraph which supports
$$P=e_{0137}-e_{0237}+e_{0267}-e_{0157}+e_{0457}$$
is as follows.
\begin{figure}[H]
	\centering
	%% Creator: Inkscape 1.0.1 (3bc2e813f5, 2020-09-07), www.inkscape.org
%% PDF/EPS/PS + LaTeX output extension by Johan Engelen, 2010
%% Accompanies image file 'EX317.pdf' (pdf, eps, ps)
%%
%% To include the image in your LaTeX document, write
%%   \input{<filename>.pdf_tex}
%%  instead of
%%   \includegraphics{<filename>.pdf}
%% To scale the image, write
%%   \def\svgwidth{<desired width>}
%%   \input{<filename>.pdf_tex}
%%  instead of
%%   \includegraphics[width=<desired width>]{<filename>.pdf}
%%
%% Images with a different path to the parent latex file can
%% be accessed with the `import' package (which may need to be
%% installed) using
%%   \usepackage{import}
%% in the preamble, and then including the image with
%%   \import{<path to file>}{<filename>.pdf_tex}
%% Alternatively, one can specify
%%   \graphicspath{{<path to file>/}}
%% 
%% For more information, please see info/svg-inkscape on CTAN:
%%   http://tug.ctan.org/tex-archive/info/svg-inkscape
%%
\begingroup%
  \makeatletter%
  \providecommand\color[2][]{%
    \errmessage{(Inkscape) Color is used for the text in Inkscape, but the package 'color.sty' is not loaded}%
    \renewcommand\color[2][]{}%
  }%
  \providecommand\transparent[1]{%
    \errmessage{(Inkscape) Transparency is used (non-zero) for the text in Inkscape, but the package 'transparent.sty' is not loaded}%
    \renewcommand\transparent[1]{}%
  }%
  \providecommand\rotatebox[2]{#2}%
  \newcommand*\fsize{\dimexpr\f@size pt\relax}%
  \newcommand*\lineheight[1]{\fontsize{\fsize}{#1\fsize}\selectfont}%
  \ifx\svgwidth\undefined%
    \setlength{\unitlength}{228.60046627bp}%
    \ifx\svgscale\undefined%
      \relax%
    \else%
      \setlength{\unitlength}{\unitlength * \real{\svgscale}}%
    \fi%
  \else%
    \setlength{\unitlength}{\svgwidth}%
  \fi%
  \global\let\svgwidth\undefined%
  \global\let\svgscale\undefined%
  \makeatother%
  \begin{picture}(1,0.33095961)%
    \lineheight{1}%
    \setlength\tabcolsep{0pt}%
    \put(0,0){\includegraphics[width=\unitlength,page=1]{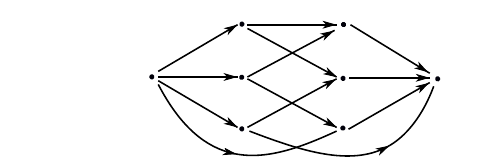}}%
    \put(-0.00221511,0.16009585){\makebox(0,0)[lt]{\lineheight{1.25}\smash{\begin{tabular}[t]{l}$\Supp(P)=$\end{tabular}}}}%
    \put(0.25612647,0.15918408){\makebox(0,0)[lt]{\lineheight{1.25}\smash{\begin{tabular}[t]{l}$0$\end{tabular}}}}%
    \put(0.48003564,0.30279932){\makebox(0,0)[lt]{\lineheight{1.25}\smash{\begin{tabular}[t]{l}$1$\end{tabular}}}}%
    \put(0.47778553,0.19701229){\makebox(0,0)[lt]{\lineheight{1.25}\smash{\begin{tabular}[t]{l}$2$\end{tabular}}}}%
    \put(0.47871027,0.09381995){\makebox(0,0)[lt]{\lineheight{1.25}\smash{\begin{tabular}[t]{l}$4$\end{tabular}}}}%
    \put(0.69728832,0.30351436){\makebox(0,0)[lt]{\lineheight{1.25}\smash{\begin{tabular}[t]{l}$3$\end{tabular}}}}%
    \put(0.7072981,0.19655258){\makebox(0,0)[lt]{\lineheight{1.25}\smash{\begin{tabular}[t]{l}$5$\end{tabular}}}}%
    \put(0.69338249,0.09444391){\makebox(0,0)[lt]{\lineheight{1.25}\smash{\begin{tabular}[t]{l}$6$\end{tabular}}}}%
    \put(0.9471566,0.15210846){\makebox(0,0)[lt]{\lineheight{1.25}\smash{\begin{tabular}[t]{l}$7$\end{tabular}}}}%
  \end{picture}%
\endgroup%

	%\caption[]{$\Supp(e_{0123})$}
	%\label{Fig:Q1}
\end{figure}
We can compute $\partial P$ and obtain
\begin{align*}
&\partial(e_{0137}-e_{0237}+e_{0267}-e_{0157}+e_{0457})\\
=~&(e_{023}-e_{013})+(e_{015}-e_{045})-e_{026}\\
&+(e_{137}-e_{157})-(e_{237}-e_{267})+e_{457}+(e_{047}-e_{067}).
\end{align*}
\end{example}

One can compute directly that $\widetilde{H}_*(\Supp(P))=0$.

\begin{example}[Exotic cube]\label{exotic} The minimal digraph which supports
$$P=e_{0158}-e_{0258}+e_{0268}-e_{0368}+e_{0378}-e_{0478}$$
is as follows.
\begin{figure}[H]
	\centering
	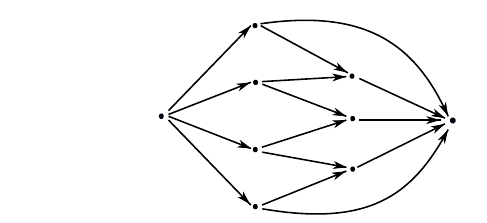
	%\caption[]{$\Supp(e_{0123})$}
	%\label{Fig:Q1}
\end{figure}

We can compute $\partial P$ and obtain
\begin{align*}
\partial P~=~&-(e_{015}-e_{025})-(e_{026}-e_{036})-(e_{037}-e_{047})\\
           ~&+e_{158}-(e_{258}-e_{268})-(e_{368}-e_{378})-e_{478} \\
           &+(e_{018}-e_{048}).
\end{align*}

One can compute directly that $\widetilde{H}_*(\Supp(P))=0$.

\end{example}

\begin{remark} Unlike the above other examples, there seem to be no deformation retractions from $\Supp(P)$ in Examples \ref{firstnoncontractible}, \ref{exotic} to its minimal $2$-faces. One can check one by one that such $\Supp(P)$ can not deform retract to any of its sub-digraph, which means that it is not contractible.

From another point of view, one can see that there are no digraph maps from the 3-cube $I_1^{\boxdot3}$ onto the above 2 examples by simple analysis:
\begin{itemize}
  \item For Example \ref{firstnoncontractible}, we have
  $$|E(\Supp(P))|=13>12=|E(I_1^{\boxdot 3})|,$$
  then this digraph can not be the image of $I_1^{\boxdot3}$;
  \item For Example \ref{exotic}, i.e., the exotic cube, we have
  $$|V(\text{exotic cube})|=9>8=|V(I_1^{\boxdot 3})|,$$
  then the exotic cube can also not be the image of $I_1^{\boxdot3}$.
\end{itemize}
Such an interpretation could be understood via the singular cubical homology of digraph, which is introduced in \cite{GJM} by Grigoryan-Jimenez-Muranov.
They proved that such a singular cubical homology is a homotopic invariant (Corollary 4.6). However, there exists a non-trivial 2-cycle in the singular cubical homology of the exotic cube, which also means that the exotic cube is not contractible.
\end{remark}

%Now learning from the above examples, we can see that there are infinitely many minimal paths of length $3$. For examples,

If we increase the number of the vertices and edges in an appropriate way, there are more examples of minimal $3$-paths.

\begin{example}\label{more1} Let us consider the path
\begin{align*}
P=~&e_{S05E}-e_{S15E}+e_{S17E}-e_{S37E}+e_{S39E}-e_{S49E}\\
   &+e_{S48E}-e_{S28E}+e_{S26E}-e_{S06E}.
\end{align*}
Its supporting digraph is given by
\begin{figure}[H]
	\centering
	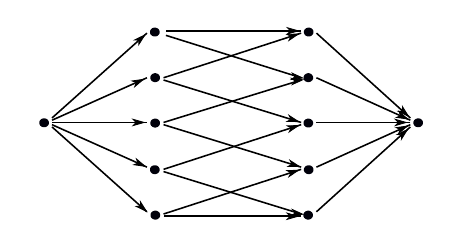
	%\caption[]{$\Supp(e_{0123})$}
	%\label{Fig:Q1}
\end{figure}

Note that one can think of this digraph as two $3$-cubes pasting together along
$$S\rightarrow2\rightarrow 7\rightarrow E$$
and then removing the edge $2\rightarrow 7$, as the following picture shows.
\begin{figure}[H]
	\centering
	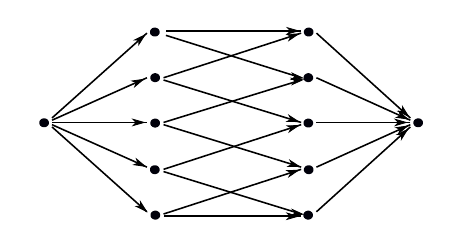
	%\caption[]{$\Supp(e_{0123})$}
	%\label{Fig:Q1}
\end{figure}
Note that the edge $2\rightarrow 7$ is necessary for the two paths
\begin{align*}
&P_1=e_{S05E}-e_{S15E}+e_{S17E}-e_{S27E}+e_{S26E}-e_{S06E}\\
&P_2=e_{S27E}-e_{S37E}+e_{S39E}-e_{S49E}+e_{S48E}-e_{S28E},
\end{align*}
to be $\partial$-invariant individually.
\end{example}

We can continue doing such a pasting-removing operation, and then we obtain infinitely many digraphs which play the role of supporting digraphs of the corresponding minimal $3$-paths.

\begin{example}\label{more2} Let us consider the path
\begin{align*}
P=~&e_{S06E}-e_{S16E}+e_{S17E}-e_{S27E}+e_{S28E}-e_{S38E}\\
   &+e_{S39E}-e_{S49E}+e_{S4(10)E}-e_{S5(10)E}.
\end{align*}
Its supporting digraph is given by
\begin{figure}[H]
	\centering
	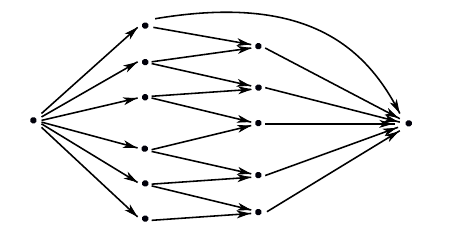
	%\caption[]{$\Supp(e_{0123})$}
	%\label{Fig:Q1}
\end{figure}

One can also think of this digraph as two smaller digraphs pasting together with some necessary edges removed, see, for examples
\begin{figure}[H]
	\centering
	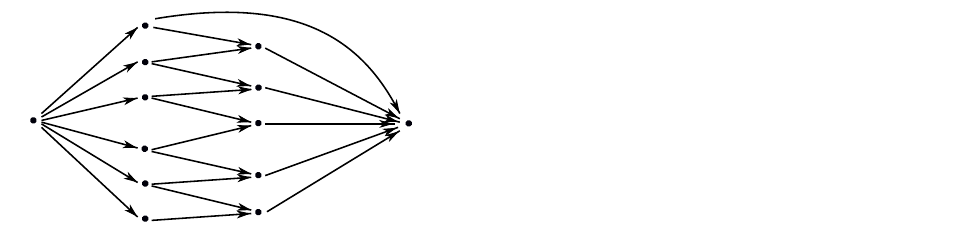
	%\caption[]{$\Supp(e_{0123})$}
	%\label{Fig:Q1}
\end{figure}

\end{example}

One can also compute directly that the supporting digraphs in Examples \ref{more1}, \ref{more2} have trivial reduced path homologies.
\vskip 0.2cm

\section{Acyclic Model}\label{acyclicsection}

We have learnt that all the examples of supporting digraphs in Subsection \ref{exminmal3} have acyclic path homologies. In this section, we will prove such an observation in the general case, which is our second main theorem.

\begin{theorem}[Acyclic Models]\label{acyclicresult} Let $P$ be a minimal path of $G$, and $\Supp(P)$ be its supporting subgraph. Then $\Supp(P)$ has acyclic path homologies, that is,
$$H_i(\Supp(P);\mathbb{Z})=0,~i>0;\qquad H_0(\Supp(P);\mathbb{Z})=\mathbb{Z}.$$
Or equivalently, $\widetilde{H}_*(\Supp(P);\mathbb{Z})=0$.
\end{theorem}

\subsection{Mayer-Vietoris exact sequences in the path complex}\label{MVsubsection}

We will use the technique of Mayer-Vietoris exact sequence in the path complex in \cite{GJMY}. Let us recall the necessary results first. It is worth to mention that such an exact sequence result holds for any abelian group as the coefficient. Here we still work with the $\Z$-coefficient.

Let $Y_1$ and $Y_2$ be two digraphs, and
$$X=Y_1\cup Y_2,\quad Z=Y_1\cap Y_2.$$
Let $i_1:Z\rightarrow Y_1$, $i_2:Z\rightarrow Y_2$ and $j_1:Y_1\rightarrow X$, $j_2:Y_2\rightarrow X$ be the corresponding inclusion digraph maps.
Suppose that any s-regular allowed elementary path on $X$ lies in $Y_1$ or $Y_2$. Then for any $p\geq-1$, there exists a short exact sequence of abelian
groups:
\begin{equation}\label{shortA}
0\longrightarrow \A_p(Z)\stackrel{\delta}{\longrightarrow}\A_p(Y_1)\oplus \A_p(Y_2)\stackrel{d}{\longrightarrow}\A_p(X)\longrightarrow0,
\end{equation}
where $\delta=(i_*^1,i_*^2)$, $d(a,b)=j^1_*(a)-j^2_*(b)$.

\begin{lemma}[\cite{GJMY} Mayer-Vietoris exact sequence]\label{MVexact} Let $Y_1$, $Y_2$, $X$, $Z$ be as above. If the homomorphism
$$\widetilde{\Omega}_p(Y_1)\oplus\widetilde{\Omega}_p(Y_2)\rightarrow\widetilde{\Omega}_p(X)$$
induced by $d$ is an epimorphism for any $p\geq-1$.  Then there is a short exact sequence of chain complexes
$$0\longrightarrow \widetilde{\Omega}_p(Z)\stackrel{\delta}{\longrightarrow}\widetilde{\Omega}_p(Y_1)\oplus \widetilde{\Omega}_p(Y_2)\stackrel{d}{\longrightarrow}\widetilde{\Omega}_p(X)\longrightarrow0.$$
And moreover we have the long exact sequence
$$\cdots\rightarrow \widetilde{H}_p(Z)\stackrel{\delta}{\longrightarrow} \widetilde{H}_p(Y_1)\oplus\widetilde{H}_p(Y_2)\stackrel{d}{\longrightarrow}\widetilde{H}_p(X) \stackrel{\partial'}{\longrightarrow}\widetilde{H}_{p-1}(Z) \stackrel{\delta}{\longrightarrow}\widetilde{H}_{p-1}(Y_1)\oplus \widetilde{H}_{p-1}(Y_2)\rightarrow\cdots$$
\end{lemma}

\subsection{Proof of Theorem \ref{acyclicresult}}\label{Proof}

This subsection is devoted to proving Theorem \ref{acyclicresult}. Before the proof, let us focus on another important observation.

\subsubsection{Non-degenerate minimal path}

Let $P$ be a minimal $n$-path of $G$, it is possible that there exists a minimal path $P^+$ of length bigger than $n$ in $\Supp(P)$. For example,

\begin{example} The supporting digraph of the minimal $3$-path $P=e_{0124}-e_{0134}+e_{0234}$ is as follows.
\begin{figure}[H]
	\centering
	%% Creator: Inkscape 1.0.1 (3bc2e813f5, 2020-09-07), www.inkscape.org
%% PDF/EPS/PS + LaTeX output extension by Johan Engelen, 2010
%% Accompanies image file '01234.pdf' (pdf, eps, ps)
%%
%% To include the image in your LaTeX document, write
%%   \input{<filename>.pdf_tex}
%%  instead of
%%   \includegraphics{<filename>.pdf}
%% To scale the image, write
%%   \def\svgwidth{<desired width>}
%%   \input{<filename>.pdf_tex}
%%  instead of
%%   \includegraphics[width=<desired width>]{<filename>.pdf}
%%
%% Images with a different path to the parent latex file can
%% be accessed with the `import' package (which may need to be
%% installed) using
%%   \usepackage{import}
%% in the preamble, and then including the image with
%%   \import{<path to file>}{<filename>.pdf_tex}
%% Alternatively, one can specify
%%   \graphicspath{{<path to file>/}}
%% 
%% For more information, please see info/svg-inkscape on CTAN:
%%   http://tug.ctan.org/tex-archive/info/svg-inkscape
%%
\begingroup%
  \makeatletter%
  \providecommand\color[2][]{%
    \errmessage{(Inkscape) Color is used for the text in Inkscape, but the package 'color.sty' is not loaded}%
    \renewcommand\color[2][]{}%
  }%
  \providecommand\transparent[1]{%
    \errmessage{(Inkscape) Transparency is used (non-zero) for the text in Inkscape, but the package 'transparent.sty' is not loaded}%
    \renewcommand\transparent[1]{}%
  }%
  \providecommand\rotatebox[2]{#2}%
  \newcommand*\fsize{\dimexpr\f@size pt\relax}%
  \newcommand*\lineheight[1]{\fontsize{\fsize}{#1\fsize}\selectfont}%
  \ifx\svgwidth\undefined%
    \setlength{\unitlength}{170.88346707bp}%
    \ifx\svgscale\undefined%
      \relax%
    \else%
      \setlength{\unitlength}{\unitlength * \real{\svgscale}}%
    \fi%
  \else%
    \setlength{\unitlength}{\svgwidth}%
  \fi%
  \global\let\svgwidth\undefined%
  \global\let\svgscale\undefined%
  \makeatother%
  \begin{picture}(1,0.30959558)%
    \lineheight{1}%
    \setlength\tabcolsep{0pt}%
    \put(0,0){\includegraphics[width=\unitlength,page=1]{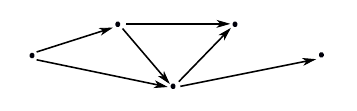}}%
    \put(-0.00202826,0.13375514){\makebox(0,0)[lt]{\lineheight{1.25}\smash{\begin{tabular}[t]{l}$0$\end{tabular}}}}%
    \put(0.30816414,0.28539158){\makebox(0,0)[lt]{\lineheight{1.25}\smash{\begin{tabular}[t]{l}$1$\end{tabular}}}}%
    \put(0.46474308,0.00287712){\makebox(0,0)[lt]{\lineheight{1.25}\smash{\begin{tabular}[t]{l}$2$\end{tabular}}}}%
    \put(0.65157022,0.28736846){\makebox(0,0)[lt]{\lineheight{1.25}\smash{\begin{tabular}[t]{l}$3$\end{tabular}}}}%
    \put(0.95859952,0.13790683){\makebox(0,0)[lt]{\lineheight{1.25}\smash{\begin{tabular}[t]{l}$4$\end{tabular}}}}%
    \put(0,0){\includegraphics[width=\unitlength,page=2]{01234.pdf}}%
  \end{picture}%
\endgroup%

	%\caption[]{$\Supp(e_{0123})$}
	%\label{Fig:Q1}
\end{figure}
It is also the supporting digraph of the length-$4$ path $e_{01234}$.
\end{example}

\begin{definition}\label{nondeg} If $P$ is the minimal path of maximal length in $\Supp(P)$, we call $P$ the non-degenerate minimal path.
\end{definition}

\subsubsection{Inductions}

In the following, we will use the non-degenerate minimal path $P$ to do the analysis for $\Supp(P)$.

The proof is done by several inductions on the length first and then on the number of the vertices in the same position of the elementary path component of $P$. The first induction is on the length of $P$. Below are the initial result (1) and the hypothesis (1).

\textbf{Initial result (1).} If $P$ is a minimal path of length $0$ or $1$ or $2$, then $\widetilde{H}_*(\Supp(P))=0$.

\textbf{Hypothesis (1).} We assume the acyclic result holds for any (non-degenerate) minimal path of length $k<n$.

To prove the result for the case $k=n$, we further do the second induction on the number $|E_1|$.

\textbf{Initial result (2)} Let us look at the simplest case, $|E_1|=1$. By Theorem \ref{structurethm}, if $|E_1|=1$, the minimal path $P$ is of the form
$$P=P_{S,n-1,\alpha}E, \quad E_1=d_E^{-1}(1)=\{\alpha\},$$
and the corresponding supporting digraph $\Supp(P)$ owns the property that each point $\beta\in d_{\alpha}^{-1}(1)=d_E^{-1}(2)$, we have the directed edge $\beta\rightarrow E$. Then we have the deformation retraction $r$:
\begin{align*}
r:&\Supp(P)\rightarrow\Supp(P_{S,n-1,\alpha}),\\
 &r(E)=\alpha,\quad  r\big|_{\Supp(P_{S,n-1,\alpha})}=\id_{\Supp(P_{S,n-1,\alpha})}.
\end{align*}
Thus
$$H_*(\Supp(P))\cong H_*(\Supp(P_{S,n-1,\alpha})).$$
By Hypothesis (1), we get our second initial result.

Next, let us continue discussing the case $|E_1|=2$. It is the key point to apply the induction method. We state the idea into the following three steps.
\begin{description}
  \item[Step 1] We embed $\Supp(P)$ into a larger digraph $\widehat{\Supp(P)}$ by adding some new edges from points in $E_2=d_E^{-1}(2)$ to $E$.
  \item[Step 2] Show that $\widehat{\Supp(P)}$ is acyclic by using the Mayer-Vietoris method, where we need the third induction for the intersection of the MV pair.
  \item[Step 3] Prove that the embedding $\Supp(P)\hookrightarrow\widehat{\Supp(P)}$ induces the isomorphism on the path homology groups.
\end{description}

Now assume that $E_1=\{\alpha_1,\alpha_2\}$. We embed $\Supp(P)$ to $\widehat{\Supp(P)}$ which is obtained from $\Supp(P)$ by adding the new edges
$$\beta\rightarrow E,\text{ whenever }\beta\in E_2=d_E^{-1}(2), \beta\nrightarrow E.$$
For the following example, we add $2\rightarrow 5$ (the blue edge) to get $\widehat{\Supp(P)}$.
\begin{figure}[H]
	\centering
	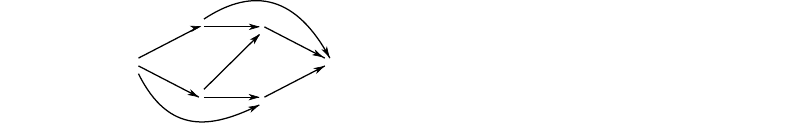
	%\caption[]{$\Supp(e_{0123})$}
	%\label{Fig:Q1}
\end{figure}

The following two observations for the new digraph $\widehat{\Supp(P)}$ are important for us.
\begin{itemize}
  \item The two paths
  $$P_{S,n-1,\alpha_1}E, P_{S,n-1,\alpha_2}E\in\Omega_n(\widehat{\Supp(P)})$$
  are two minimal paths. For simplicity, we denote them by $P_{\alpha_1}$, $P_{\alpha_2}$ respectively.
  \item We have the following relations among $\widehat{\Supp(P)}$, $\Supp(P_{\alpha_1})$ and $\Supp(P_{\alpha_2})$:
  \begin{align*}
  &\widehat{\Supp(P)}=\Supp(P_{\alpha_1})\cup\Supp(P_{\alpha_2})\\
  &\Supp(P_{\alpha_1})\cap\Supp(P_{\alpha_2})=\bigcup_{\beta\in d_E^{-1}(2),\beta\nrightarrow E \text{ in }\Supp(P)}\Supp(P_{S,n-2,\beta})E.
  \end{align*}
\end{itemize}

\begin{lemma}\label{Y1Y2MV} Let $X=\widehat{\Supp(P)}$, $Y_1=\Supp(P_{\alpha_1})$ and $Y_2=\Supp(P_{\alpha_2})$. Then the digraph pair $(Y_1,Y_2)$ forms a Mayer-Vietoris pair of the new digraph $X$.
\end{lemma}

\begin{proof} According to the above observations, it suffices to prove that the homomorphism in Lemma \ref{MVexact}
$$\widetilde{\Omega}_p(Y_1)\oplus\widetilde{\Omega}_p(Y_2)\rightarrow\widetilde{\Omega}_p(X)$$
is an epimorphism for any $p\geq-1$.

It is obvious for $p=-1,0,1,2$. By the construction of $X$, $Y_1$, $Y_2$, we only need to consider the minimal path $u\in\Omega_p(X)$ starting from some vertex $U$ and ending with $E$. We can write $u$ as
$$u=\sum_{\gamma_1\in V_1'}P_{U,p-1,\gamma_1}E+\sum_{\gamma_2\in V_2'}P_{U,p-1,\gamma_2}E+\sum_{\gamma_{12}\in V_{12}'}P_{U,p-1,\gamma_{12}}E,$$
where
\begin{itemize}
  \item $V_1'\subset V_1=\{\gamma_1\in V(Y_1)\big|\gamma_1\rightarrow E \text{ in }\Supp(P),\gamma_1\notin V(Y_2)\}$;
  \item $V_2'\subset V_2=\{\gamma_2\in V(Y_1)\big|\gamma_2\rightarrow E \text{ in }\Supp(P),\gamma_2\notin V(Y_1)\}$;
  \item $V_{12}'\subset V_{12}=\{\gamma_{12}\in V(Y_1)\cap V(Y_2)|\gamma_{12}\rightarrow E \text{ in }X\}$.
\end{itemize}
By our construction and the structure theorem, note that
\begin{itemize}
  \item Each $P_{U,p-1,\gamma}\in\Omega_{p-1}(\Supp(P))$ is the unique minimal $(p-1)$-path starting from $U$ and ending with $\gamma$, which also lives in $Y_1$ or $Y_2$ depending on $\gamma\in V(Y_1)$ or $V(Y_2)$, respectively;
  \item $V_{12}=\{\beta\in d_E^{-1}(2)\big|\beta\nrightarrow E\text{ in } \Supp(P)\}$, i.e.
  $$P_{U,p-1,\gamma_{12}}E\in\mathcal{A}_p(X)\setminus\mathcal{A}_{p}(\Supp(P)).$$
\end{itemize}
%$$V_{12}=\{\beta\in d_E^{-1}(2)\big|\beta\nrightarrow E\text{ in } \Supp(P)\}.$$

Without loss of generality, assume that $u\notin\Omega_p(Y_1)$ and $u\notin\Omega_p(Y_2)$, which corresponds to $V_1',V_2'\neq\emptyset$. Now for $i=1,2$, we consider the minimal $\partial$-invariant completion of $\sum_{\gamma_i\in V_i'}P_{U,p-1,\gamma_i}$ in $Y_i$, we denote it as $P_{U,p,E}^i$.

\begin{claim}\label{MVclaim} $P_{U,p,E}^1-\sum_{\gamma_1\in V_1'}P_{U,p-1,\gamma_1}E$ or $P_{U,p,E}^2-\sum_{\gamma_2\in V_2'}P_{U,p-1,\gamma_2}E$ must be of the form
$$\sum_{\gamma_{12}'\in V_{12}''\subset V_{12}}P_{U,p-1,\gamma_{12}'}E\in\mathcal{A}_{p-1}(Y_1)\cap\mathcal{A}_{p-1}(Y_2).$$
\end{claim}

\begin{proof}[Proof of Claim \ref{MVclaim}] As the minimal $\partial$-invariant completion, we can decompose $P_{U,p,E}^1\in\Omega_p(Y_1)$ and $P_{U,p,E}^2\in\Omega_p(Y_2)$ as
\begin{align*}
P_{U,p,E}^1=~&\sum_{\gamma_1\in V_1'}P_{U,p-1,\gamma_1}E+\sum_{\gamma_1'\in V_1''}P_{U,p-1,\gamma_1'}E+\sum_{\gamma_{12}'\in V_{12}''\subset V_{12}}P_{U,p-1,\gamma_{12}'}E,\\
P_{U,p,E}^2=~&\sum_{\gamma_2\in V_2'}P_{U,p-1,\gamma_2}E+\sum_{\gamma_2'\in V_2''}P_{U,p-1,\gamma_2'}E+\sum_{\gamma_{12}'\in V_{12}'''\subset V_{12}}P_{U,p-1,\gamma_{12}'}E,
\end{align*}
where $V_i''\subset V_i$ and $V_i'\cap V_i''=\emptyset$, $i=1,2$.

We want to prove that $V_1''=\emptyset$ or $V_2''=\emptyset$. If $V_1''$, $V_2''\neq \emptyset$, it means that there exists $\gamma_i'\in V_i''\setminus V_i'$ (since $V_i'\cap V_i''=\emptyset$), such that
$P_{U,p-1,\gamma_i'}E\in\mathcal{A}_{p}(\Supp(P))$, and it will cancel some terms in $P_{U,p-1,\gamma_i}E$ for some $\gamma_i\in V_i'$.

If we can not decompose $u$ into two $\partial$-invariant paths in $Y_1$ and $Y_2$ respectively, it means that there must exist a term $P_{U,p-1,\gamma_1}E$ for some $\gamma_1\in V_1'$, which cancels some terms in $P_{U,p-1,\gamma_2}$ for some $\gamma_2\in V_2'$. Then in $\Supp(P)$, there will exist a point $\beta$, such that
$$e_{\beta\gamma_1E}-e_{\beta\gamma_2E}\text{ and } e_{\beta\gamma_1E}-e_{\beta\gamma_1'E}\in\Omega_2(\Supp(P)).$$
It is in contradiction with Theorem \ref{structurethm} (3), since $P$ is a minimal path.\qedhere
\end{proof}

It follows from the claim that
\begin{align*}
u=~&P_{U,p,E}^1-\sum_{\gamma_{12}'\in V_{12}''\subset V_{12}}P_{U,p-1,\gamma_{12}'}E+\sum_{\gamma_2\in V_2'}P_{U,p-1,\gamma_2}E+\sum_{\gamma_{12}\in V_{12}'}P_{U,p-1,\gamma_{12}}E\\
(\text{ or}~=~&P_{U,p,E}^2-\sum_{\gamma_{12}'\in V_{12}''\subset V_{12}}P_{U,p-1,\gamma_{12}'}E+\sum_{\gamma_1\in V_1'}P_{U,p-1,\gamma_1}E+\sum_{\gamma_{12}\in V_{12}'}P_{U,p-1,\gamma_{12}}E~)\\
 \in~&\Omega_p(Y_1)\oplus \Omega_p(Y_2).
\end{align*}
Then we are done.
\end{proof}

Let $Z=Y_1\cap Y_2$, by Lemmas \ref{Y1Y2MV} \ref{MVexact}, we obtain the long exact sequence
$$\cdots\rightarrow \widetilde{H}_p(Z)\stackrel{\delta}{\longrightarrow} \widetilde{H}_p(Y_1)\oplus\widetilde{H}_p(Y_2)\stackrel{d}{\longrightarrow}\widetilde{H}_p(X) \stackrel{\partial'}{\longrightarrow}\widetilde{H}_{p-1}(Z) \stackrel{\delta}{\longrightarrow}\widetilde{H}_{p-1}(Y_1)\oplus \widetilde{H}_{p-1}(Y_2)\rightarrow\cdots$$
By the argument of initial result (2), we know that $\widetilde{H}(Y_1)=\widetilde{H}(Y_2)=0$, thus
\begin{equation}\label{newsupp}
\widetilde{H}_*(\widehat{\Supp(P)})\cong \widetilde{H}_*(Z).
\end{equation}
We want to show that $\widehat{\Supp(P)}$ is acyclic through $Z$. We will understand the structure of $Z$ as follows.
\begin{description}
  \item[(a)] $Z$ could be regard as a proper subgraph of $\Supp(P_{S,n-1,\alpha_1})$ or $\Supp(P_{S,n-1,\alpha_2})$ (It can embed to the two digraphs by mapping $E$ to $\alpha_1$ or $\alpha_2$ with other points fixed). For convenience, we call $n-1$ the length of the digraph $Z$.
  \item[(b)] If we remove the point $E$, the reduced subgraph of $Z$ is a union of supporting digraphs of minimal $(n-2)$-path in $\Supp(P)$ with the same starting point $S$.
\end{description}

We will prove the acyclic result of $Z$ via the third induction on the length of $Z$, that is $n-1$.

\textbf{Initial result (3).} Note that
\begin{itemize}
  \item For $n=2$, $Z=\{S\rightarrow E\}$. It is clearly acyclic.
  \item For $n=3$, $Z$ is a digraph of the form as follows.
  \begin{figure}[H]
	\centering
	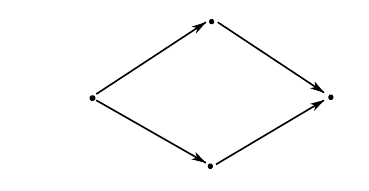
	%\caption[]{}
	%\label{Fig:Q1}
\end{figure}
  It is easy to check that $Z$ is acyclic in this case.
\end{itemize}

\textbf{Hypothesis (3).} We assume that $Z$ of the above form (satisfying (a) (b)) with $k<n-1$ is acyclic.

Now let us prove the case $k=n-1$. We repeat the above analysis and do the induction on the number of points in $E_2\cap V(Z)=d_E^{-1}(2)\cap V(Z)$.

\textbf{Initial result (2')} If $|E_2\cap V(Z)|=1$, the argument is the same as the proof of the initial result (2).

If $|E_2\cap V(Z)|=2$ and $E_2\cap V(Z)=\{\beta_1,\beta_2\}$. We embed $Z$ into a larger digraph $\widehat{Z}$ by adding new edges:
$$\gamma\rightarrow E,\quad \gamma\in d_E^{-1}(3)\cap V(Z).$$
Similarly, we have the MV pair $(Z_1,Z_2)$ for $\widehat{Z}$, where $Z_1$ and $Z_2$ are supporting digraphs of the minimal paths $P_{S,n-2,\beta_1}E$ and $P_{S,n-2,\beta_2}E$ respectively. By repeating using the Mayer-Vietoris sequence, and applying Hypothesis (3) to $Z_1\cap Z_2$, we obtain
$$\widetilde{H}_*(\widehat{Z})=0.$$

\begin{lemma}\label{ZhatZ} Let $Z$ and $\widehat{Z}$ be two digraphs as above corresponding to the case $E_2\cap V(Z)=\{\beta_1,\beta_2\}$. Then
$$\widetilde{H}_*(Z)=\widetilde{H}_*(\widehat{Z})=0.$$
\end{lemma}
\begin{proof} Let $u\in\Omega_p(Z)$, and $\partial u=0$. Since $\widetilde{H}_*(\widehat{Z})=0$, there exists $v\in\Omega_{p+1}(\widehat{Z})$, such that
$$\partial v=u.$$
Now we want to replace $v$ by $v'\in\Omega_{p+1}(Z)$, such that $\partial v'=u$.

By construction, $\Omega_*(\widehat{Z})$ only has several more paths ending with $E$ than $\Omega_*(Z)$. Thus, we only need to consider the path component ending with $E$ in $v$. Meanwhile, since $(Z_1,Z_2)$ is the Mayer-Vietoris pair for $\widehat{Z}$. Thus, we can represent $v\in\Omega_{p+1}(\widehat{Z})$ as
\begin{align*}
v=~&\left(\sum_{\gamma_1}p_{\gamma_1}+\sum_{\epsilon_1}p_{\epsilon_1}\right)\beta_1E+\sum_{\gamma_1}\sum_{\delta_1}p_{\delta_1}\gamma_1E \quad(:=v_1\in\Omega_{p+1}(Z_1))\\
 ~&+\left(\sum_{\gamma_2}p_{\gamma_2}+\sum_{\epsilon_2}p_{\epsilon_2}\right)\beta_2E+\sum_{\gamma_2}\sum_{\delta_2}p_{\delta_2}\gamma_2E \quad(:=v_2\in\Omega_{p+1}(Z_2)),
\end{align*}

where the summands are explained as follows:
\begin{itemize}
  \item $\left(\sum_{\gamma_1}p_{\gamma_1}+\sum_{\epsilon_1}p_{\epsilon_1}\right)\beta_1\in\Omega_p(Z)$ ends at $\beta_1$, $p_{\gamma_1}$ ends at $\gamma_1$, $p_{\epsilon_1}$ ends at $\epsilon_1$, where $\gamma_1\in d_E^{-1}(3)$, $\delta_1\in d_E^{-1}(4)$.
  \item $\left(\sum_{\gamma_2}p_{\gamma_2}+\sum_{\epsilon_2}p_{\epsilon_2}\right)\beta_2\in\Omega_p(Z)$ ends at $\beta_2$, $p_{\gamma_2}$ ends at $\gamma_2$, $p_{\epsilon_2}$ ends at $\epsilon_2$, where $\gamma_2\in d_E^{-1}(3)$, $\delta_2\in d_E^{-1}(4)$.
  \item $\sum_{\gamma_i}\sum_{\delta_i}p_{\delta_i}\gamma_i\in\Omega_{p}(Z)$, $\sum_{\delta_i}$ depends on $\gamma_i$, $i=1,2$.
  While $\sum_{\gamma_i}\sum_{\delta_i}p_{\delta_i}\gamma_i E\in\mathcal{A}_{p+1}(\widehat{Z})\setminus\mathcal{A}_{p+1}(Z)$.
\end{itemize}

Now let us analyze the condition $\partial v=u\in\Omega_p(Z)$ as follows.
\begin{align*}
\partial v=~&\left[\partial\left(\sum_{\gamma_1}p_{\gamma_1}+\sum_{\epsilon_1}p_{\epsilon_1}\right)\right]\beta_1E\pm \underline{\left(\sum_{\gamma_1}p_{\gamma_1}+\sum_{\epsilon_1}p_{\epsilon_1}\right)E}\pm \left(\sum_{\gamma_1}p_{\gamma_1}+\sum_{\epsilon_1}p_{\epsilon_1}\right)\beta_1\\
          ~&+\underline{\sum_{\gamma_1}\left(\partial\sum_{\delta_1}p_{\delta_1}\right)\gamma_1E}\pm \underline{\sum_{\gamma_1}\sum_{\delta_1}p_{\delta_1}E}\pm \sum_{\gamma_1}\sum_{\delta_1}p_{\delta_1}\gamma_1\\
          ~&+\left[\partial\left(\sum_{\gamma_2}p_{\gamma_2}+\sum_{\epsilon_2}p_{\epsilon_2}\right)\right]\beta_2E\pm \underline{\left(\sum_{\gamma_2}p_{\gamma_2}+\sum_{\epsilon_2}p_{\epsilon_2}\right)E}\pm \left(\sum_{\gamma_2}p_{\gamma_2}+\sum_{\epsilon_2}p_{\epsilon_2}\right)\beta_2\\
          ~&+\underline{\sum_{\gamma_2}\left(\partial\sum_{\delta_2}p_{\delta_2}\right)\gamma_2E}\pm \underline{\sum_{\gamma_2}\sum_{\delta_2}p_{\delta_2}E}\pm \sum_{\gamma_2}\sum_{\delta_2}p_{\delta_2}\gamma_2
\end{align*}

First, $v_1\in\Omega_{p+1}(Z_1)$ and $v_2\in\Omega_{p+1}(Z_2)$ imply that
\begin{equation}\label{Z1Z2}
\begin{aligned}
&\sum_{\epsilon_1}p_{\epsilon_1}E+\sum_{\gamma_1}\sum_{\delta_1}p_{\delta_1}E=0,\quad (\text{since they are non-allowed in $\partial v_1$});\\
&\sum_{\epsilon_2}p_{\epsilon_2}E+\sum_{\gamma_2}\sum_{\delta_2}p_{\delta_2}E=0,\quad (\text{since they are non-allowed in $\partial v_2$}).\\
\end{aligned}
\end{equation}

Second, $\partial v=u\in\Omega_{p}(Z)$ and the terms with underlines are non-allowed in $Z$. Combined with \eqref{Z1Z2}, we get
\begin{equation}\label{Z1Z2together}
\sum_{\gamma_1}p_{\gamma_1}+\sum_{\gamma_2}p_{\gamma_2}+\sum_{\gamma_1}\left(\partial\sum_{\delta_1}p_{\delta_1}\right)+\sum_{\gamma_2}\left(\partial\sum_{\delta_2}p_{\delta_2}\right)=0.
\end{equation}

We want to find $w\in\Omega_{p+2}(\widehat{Z})$ such that $v+\partial w\in\Omega_{p+1}(Z)$. To do this, let us study the above equations \eqref{Z1Z2}, \eqref{Z1Z2together} carefully.

Let us look at Equations $\eqref{Z1Z2}$. To cancel the terms $\sum_{\gamma_1}\sum_{\delta_1}p_{\delta_1}\gamma_1E$ and $\sum_{\gamma_2}\sum_{\delta_2}p_{\delta_2}\gamma_2E$ in $v$, we consider the following allowed path in $Z$:
$$w_1=\sum_{\gamma_1}\sum_{\delta_1}p_{\delta_1}\gamma_1\beta_1E+\sum_{\gamma_2}\sum_{\delta_2}p_{\delta_2}\gamma_2\beta_2E\in\mathcal{A}_{p+2}(Z).$$
Then we can compute $v+(-1)^{p}\partial w_1$ and obtain the corresponding non-allowed terms come from
\begin{equation}\label{nonallowed1}
\left[\sum_{\gamma_1}p_{\gamma_1}+(-1)^p\partial\bigg(\sum_{\gamma_1}\sum_{\delta_1}p_{\delta_1}\bigg)\gamma_1\right]\beta_1E+\left[\sum_{\gamma_2}p_{\gamma_2}+(-1)^p\partial\bigg(\sum_{\gamma_2}\sum_{\delta_2}p_{\delta_2}\bigg)\gamma_2\right]\beta_2E.
\end{equation}
To cancel the non-allowed path, let us read \eqref{Z1Z2together} by dividing $\gamma_1,\gamma_2\in d_E^{-1}(3)\cap V(Z)$ into three disjoint parts
\begin{itemize}
  \item $\gamma_1'$: $\gamma_1'\rightarrow\beta_1$, but $\gamma_1'\nrightarrow \beta_2$;
  \item $\gamma_2'$: $\gamma_2'\rightarrow\beta_2$, but $\gamma_2'\nrightarrow \beta_1$;
  \item $\gamma_{12}'$: $\gamma_{12}'\rightarrow\beta_1$, and $\gamma_{12}'\rightarrow\beta_2$.
\end{itemize}
Then $\eqref{Z1Z2together}$ implies, according to the endpoints,
\begin{align*}
&p_{\gamma_1'}+\partial\left(\sum_{\delta_1}p_{\delta_1}\right)\gamma_1'=0,\quad p_{\gamma_2'}+\partial\left(\sum_{\delta_2}p_{\delta_2}\right)\gamma_2'=0,\\
&p_{\gamma_{12}'}^{Z_1}+p_{\gamma_{12}'}^{Z_2}+\partial\left(\sum_{\delta_{12}}p_{\delta_{12}}\right)\gamma_{12}'=0.
\end{align*}
The last equation allows us to decompose the left hand side as follows:
\begin{align*}
  &p_{\gamma_{12}'}^{Z_1}+p_{\gamma_{12}'}^{Z_2}+\partial\left(\sum_{\delta_{12}}p_{\delta_{12}}\right)\gamma_{12}'\\
=~&\left[p_{\gamma_{12}'}^{Z_1}+\partial\left(\sum_{\delta_{12}}p_{\delta_{12}}'\right)\gamma_{12}'\right]+\left[p_{\gamma_{12}'}^{Z_2}+\partial\left(\sum_{\delta_{12}}\bigg(p_{\delta_{12}}-p_{\delta_{12}}'\bigg)\right)\gamma_{12}'\right]
\end{align*}
such that
$$p_{\gamma_{12}'}^{Z_1}+\partial\left(\sum_{\delta_{12}}p_{\delta_{12}}'\right)\gamma_{12}'=p_{\gamma_{12}'}^{Z_2}+\partial\left(\sum_{\delta_{12}}\bigg(p_{\delta_{12}}-p_{\delta_{12}}'\bigg)\right)\gamma_{12}'=0.$$
Then, we consider another allowed $(p+2)$-path in $Z$:
$$w_2=\sum_{\gamma_{12}'}\left(\sum_{\delta_{12}}(p_{\delta_{12}}-p_{\delta_{12}}')\right)\gamma_{12}'(\beta_1-\beta_2)E.$$
Then we can compute $\partial w_2$ and find that the corresponding non-allowed term cancels the one in \eqref{nonallowed1}. Thus we get
$$w_1-w_2\in\Omega_{p+1}(\widehat{Z})$$
and furthermore,
$$v':=v+(-1)^p\partial w_1-(-1)^p\partial w_2\in\Omega_{p+1}(Z).$$
Then we are done.
\end{proof}

Now let us return to the proof $\widetilde{H}_*(Z)=0$ for $|E_2\cap V(Z)|>2$.

\textbf{Hypothesis (2')} Assume when $|E_2\cap V(Z)|=k<m$, $\widetilde{H}_*(Z;\Z)=0$.

For $|E_2\cap V(Z)|=m$, assume that $E_2\cap V(Z)=\{\beta_1,\ldots,\beta_m\}$. Let
$$Z_1=\Supp(P_{S,n-2,\beta_1})E,\quad Z_2=\cup_{\beta_i,i=2,\ldots m}\Supp(P_{S,n-2,\beta_i})E.$$
Assume that
$$Z_1\cap Z_2\cap \{f_I^{-1}(n-3)\}_{e_I<P}=\{\gamma_1,\ldots,\gamma_l\},$$
then consider the new digraphs $\widehat{Z}$ by adding the edges $\gamma_i\rightarrow E$, $i=1,\ldots,l$.
Then $\widehat{Z_1}:=Z_1\bigcup\cup_{i=1}^l\{\gamma_i\rightarrow E\}$ and $\widehat{Z_2}:=Z_2\bigcup\cup_{i=1}^l\{\gamma_i\rightarrow E\}$ form a Mayer-Vietoris pair for $\widehat{Z}$, more explicitly
\begin{itemize}
  \item $\widehat{Z_1}\cup \widehat{Z_2}=\widehat{Z}$
  \item $\widehat{Z_1}\cap\widehat{Z_2}=\cup_j\Supp(P_{S,n-3,\gamma_j})E$, it is of the form in our second induction.
\end{itemize}

Also, we have the deformation retraction by mapping $E$ to $\beta_1$ and fixing other points. By Hypothesis (3), we known that
$$\widetilde{H}_*(\widehat{Z_1})=0.$$
Thus, again by Hypothesis (3) for $\widehat{Z_1}\cap\widehat{Z_2}$ and by Mayer-Vietoris exact sequence, we have
$$\widetilde{H}_*(\widehat{Z})\cong\widetilde{H}_*(\widehat{Z_2}).$$
Moreover, by Hypothesis (2') for $|E_2\cap V(\widehat{Z_2})|=m-1$, we have
$$\widetilde{H}_*(\widehat{Z})\cong\widetilde{H}_*(\widehat{Z_2})=0.$$

Let us repeat the argument for the case $|E_2\cap V(Z)|=2$, (note that the classification of $\gamma_i$ in the $|E_2|=2$ case still work for the general case by Theorem \ref{structurethm} (3)), thus we have
$$\widetilde{H}_*(Z)=\widetilde{H}_*(\widehat{Z})=0.$$
Then, we finish the proof of acyclic property for $Z$ of length $n-1$.

Thus, by \eqref{newsupp}, we obtain
$$\widetilde{H}_*(\widehat{\Supp(P)})=0,\quad \text{when } |E_1|=2.$$

Then, repeating the relative argument for the pair $(\widehat{\Supp(P)},\Supp(P))$ as Lemma \ref{ZhatZ}, we obtain,
$$\widetilde{H}_*(\Supp(P))=\widetilde{H}_*(\widehat{\Supp(P)})=0,\quad \text{when } |E_1|=2.$$

\textbf{Hypothesis (2).} Assume when $|E_1|=k<m$, $\widetilde{H}_*(\Supp(P);\Z)=0$. For $|E_1|=m$, assume that $E_1=\{\alpha_1,\ldots,\alpha_m\}$. To do the induction, we have
\begin{claim}\label{betterchoice} There exists $\alpha_i\in E_1$, such that if we add the edges
$$\gamma\rightarrow E, \quad\gamma\in d_E^{-1}(2)\cap V(P_{S,n-1,\alpha_i}),$$
then the two paths $P_1:=P_{S,n-1,\alpha_i}E$ and $P_2:=P-P_1$ become two minimal paths, with $|E_1\cap V(P_2)|=m-1$.
\end{claim}

\begin{proof}[Proof of Claim \ref{betterchoice}]
For any $\alpha_i$, if we add the edges from points in $d_{\alpha}^{-1}$ to $E$, we get the minimal path $P_1$. Meanwhile $P_2:=P-P_1$ becomes $\partial$-invariant, but may not be minimal. Then we can write $P_2$ as
$$P_2=P_2^1+P_2^2+\cdots+P_2^l, \text{ each $P_2^j$ is minimal}.$$
We keep $P_2^1$ minimal as $\tilde{P}_1$ and remove the redundant new edges such that $\tilde{P}_2=P-\tilde{P}_1$ is minimal. If $\tilde{P}_1$ satisfies our requirement, we are done. Otherwise, we continue the operation for $\tilde{P}_1$. After several times, we are done. An example of such an operation will be given in Appendix \ref{B}.
\end{proof}
It follows from this claim that we can repeat the above argument for the proof of $\tilde{H}_*(Z)=0$, then we are done.

\subsubsection{Proof idea via an example}

Let us use the idea in the proof to compute the path homology of the exotic cube. A more complicated example will be given in Appendix \ref{B}.
\begin{figure}[H]
	\centering
	\input{EX319.pdf_tex}
	%\caption[]{$\Supp(e_{0123})$}
	%\label{Fig:Q1}
\end{figure}

In this example, $E_1=\{5,6,7\}$, $|E_1|=3$. We add the directed edges $2\rightarrow 8$, and get the following larger digraph $\widehat{\Supp(P)}$
\begin{figure}[H]
	\centering
	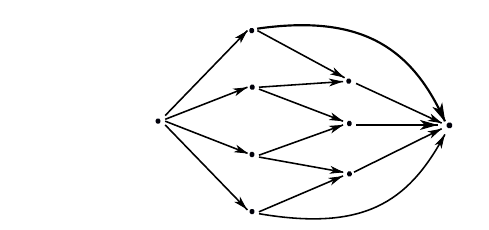
	%\caption[]{$\Supp(e_{0123})$}
	%\label{Fig:Q1}
\end{figure}

Then the digraphs $Y_1$ and $Y_2$ below form a Mayer-Vietoris pair for $\widehat{\Supp(P)}$.
\begin{figure}[H]
	\centering
	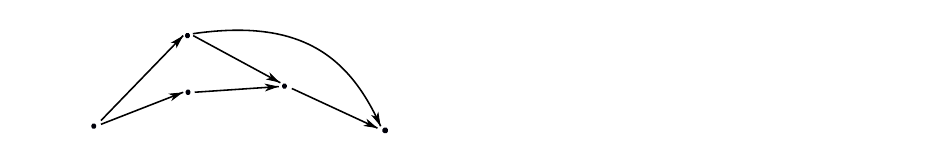
	%\caption[]{$\Supp(e_{0123})$}
	%\label{Fig:Q1}
\end{figure}
Here $Z=Y_1\cap Y_2=\{ 0\rightarrow 2\rightarrow 8\}$. Now we have,
\begin{description}
  \item[(a)] $Y_1$ can deform retract to a minimal $2$-path $e_{015}-e_{025}$. By Hypothesis (1), $\widetilde{H}(Y_1)=0$.
  \item[(b)] $Y_2$ is a minimal $3$-path with $\{6,7\}=E_1\cap V(Y_2)$, i.e. $|E_1\cap V(Y_2)|=2<3=|E_1|$. By Hypothesis (2), $\widetilde{H}(Y_2)=0$.
\end{description}

By Mayer-Vietoris exact sequence, we get
$$\widetilde{H}_*(\widehat{\Supp(P)})\cong\widetilde{H}_*(Z)=0.$$

It remains to prove that
$$\widetilde{H}_*(\Supp(P))\cong\widetilde{H}_*(\widehat{\Supp(P)}).$$

Note that
$$\Omega_2(\widehat{\Supp(P)})=\Omega_2(\Supp(P))\oplus\Span_{\Z}\{e_{018}-e_{028}, e_{258}\}.$$
If $u\in\Omega_1(\Supp(P))$, and $\partial u=0$. Then there exists $v\in\Omega_2(\widehat{\Supp(P)})$, such that $\partial v=u$. We write $v$ as
$$v=v_1+a(e_{018}-e_{028})+be_{258},\quad a,b\in\Z$$
where $v_1\in\Omega_2(\Supp(P))$.

We can compute
\begin{align*}
\partial v=~&\partial v_1+ a(e_{01}-e_{02}+e_{18}-e_{28})+b(e_{25}-e_{28}+e_{58})\\
          =~&\partial v_1+a(e_{01}-e_{02}+e_{18})+b(e_{25}+e_{58})-\underline{(a+b)e_{28}}.
\end{align*}
Then $\partial v=u\in\Omega_1(\Supp(P))$ implies that $a+b=0$.

Let $w=ae_{0158}-ae_{0258}\in\Omega_3(\widehat{\Supp(P)})$, we have
$$\partial w=a(e_{158}-e_{258}+e_{018}-e_{028}-e_{015}+e_{025}).$$
Furthermore, we define
$$v':=v+\partial w\in \mathcal{A}_{3}(\Supp(P)).$$
Since $\partial v'=\partial(v+\partial w)=\partial v=u\in\Omega_2(\Supp(P))$, we have $v'\in\Omega_3(\Supp(P))$. Thus $H_1(\Supp(P))=0$. Similarly, one check that
$H_2(\Supp(P))=H_3(\Supp(P))=0$.

\begin{remark} Note that one can also embed $\Supp(P)$ into two obvious acyclic digraphs.
\begin{enumerate}
  \item $\overline{\Supp(P)}$,  whose sets of vertices and edges are given by
  \begin{itemize}
  \item $V(\overline{\Supp(P)})=V(\Supp(P))$;
  \item $E(\overline{\Supp(P)})= \{i_a\rightarrow i_{b}~|\text{if }e_{i_0\cdots i_a\cdots i_b\cdots i_m}\in\mathcal{A}_m(\Supp(P))\}$.
\end{itemize}
  \item $\widetilde{\Supp(P)}$, by adding (if there are not such edges in $\Supp(P)$) all the edges
  $$v\rightarrow E,\quad \text{for any } v\in V(P)\setminus\{E\}.$$
\end{enumerate}
The two kinds of bigger digraphs of the exotic cube are given as follows.
\begin{figure}[H]
	\centering
	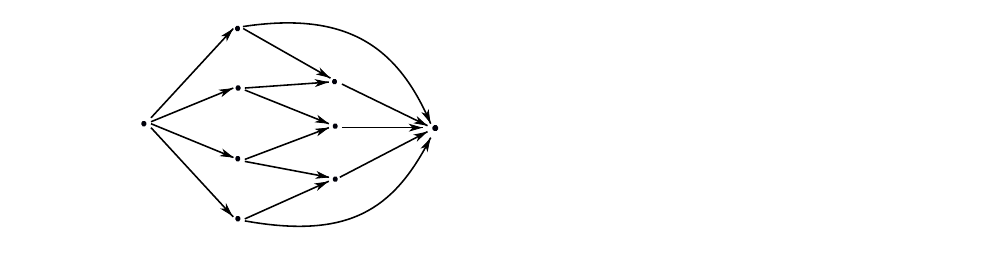
	%\caption[]{$\Supp(e_{0123})$}
	%\label{Fig:Q1}
\end{figure}
Both digraphs $\overline{\Supp(P)}$ and $\widetilde{\Supp(P)}$ can deform retract to its end vertex $E$. In particular, the first digraph is known as transitive closure \cite{BG}. Lin-Wang-Yau \cite{LWY} developed the Morse theory on the digraph and studied the homology relation between a digraph $G$ and its transitive closure $\bar{G}$ under some condition on discrete gradient vector field.

\end{remark}

\begin{remark} In fact, we prove a more general acyclic result of the digraphs of the form
$$G_1=\bigcup_i\Supp(P_{S,n-1,\alpha_i})E,\quad G_2=\bigcup_{j}\Supp(P_{S,n,E}^j),$$
where $P_{S,n-1,\alpha}$ is the minimal $(n-1)$-path and $P_{S,n,E}^j$'s are all minimal $n$-paths with the same starting vertex and ending vertex.

We obtain such a result by observing that the intersection digraph $Z$ and $\hat{Z}$ are of such two forms. But we can learn from our examples in this paper that we could only add one edge to split $\Supp(P)$ into two small supporting sub-digraphs, which means that $Z$ may have a simple look such as $\Supp(P_{S,n-1,\alpha})E$, then the proof of acyclic result will become much more simpler. To do that, we need to improve the structure theorem of minimal path $P$ and $\Supp(P)$. We leave such a question to the readers who are interested in it.
\end{remark}

\vskip 0.2cm

\section{Applications}\label{applicationsection}

In this section, we will briefly recall the path cohomology of a digraph \cite{GLMY2} and its cup product, here we switch to the field coefficient $\mathbb{K}$ ($\mathbb{K}=\mathbb{Q}$, $\mathbb{R}$ or $\mathbb{C}$). Next, we will use the method of acyclic models to construct the chain homotopy in the proof of the skew-symmetry of cup product.

\subsection{The cohomology of digraphs}\label{cohomology}
Let $V$ be a finite set. For any $n\geq0$ define a $n$-form on $V$ as any linear functional $\omega:\Lambda_n(V)\rightarrow\K$.
The linear space of all $n$-forms is denoted by $\Lambda^n(V)$. That is, $\Lambda^n$(V) is the dual linear space of $\Lambda_n(V)$.

For any elementary $n$-path $e_{i_0\cdots i_n}$, there is a dual elementary $n$-form $e^{i_0\cdots i_n}$ such that
$$(e^{i_0\cdots i_n},e_{j_0\cdots j_n})=\delta_{j_0\cdots j_n}^{i_0\cdots i_n}.$$
There is a linear operator $d:\Lambda^n(V)\rightarrow\Lambda^{n+1}(V)$, which is defined on the basis $\{e^{i_0\cdots i_n}\}$ by
$$de^{i_0\cdots i_n}=\sum_{k\in V}\sum_{p=0}^{n+1}(-1)^pe^{i_0\cdots i_{p-1}ki_p\cdots i_n}.$$

\begin{proposition}[\cite{GLMY2}] The pair $(\Lambda^*(V),d)$ is the dual complex of $(\Lambda_*(V),\partial)$. More explicitly,
\begin{enumerate}
  \item $d^2=0$;
  \item $(d\omega, u)=(\omega,\partial u)$, for any $\omega\in\Lambda^n(V), u\in\Lambda_{n+1}(V)$.
\end{enumerate}
\end{proposition}

\begin{remark}
Let $\mathcal{R}^n(V)$ be the subspace of $\Lambda^n(V)$ spanned by
$$\{e^{i_0\cdots i_n}:i_k\neq i_{k+1}, \text{ for all } k=0,1,\ldots,n-1\}.$$
Moreover, $(\mathcal{R}^*(V),d)$ is the dual complex of $(\mathcal{R}_*(V),\partial)$.
\end{remark}
To apply our result, we still restrict to the strongly regular condition, and set
$$e^{i_0\ldots i_p}=0, \quad\text{if } i_j=i_k \text{ for some }j\neq k.$$
The corresponding definitions and results hold in such a narrow case. We still denote by $\mathcal{R}^*(V)$ the corresponding subspace.

Let $G=(V, E)$ be a digraph. For any $n\geq 0$, consider the following subspaces of $\mathcal{R}^n(V)$:
\begin{itemize}
  \item $\mathcal{A}^n(G)=\Span\{e^{i_0\cdots i_n}\in\mathcal{R}^n(V): i_0\cdots i_p \text{ is allowed.}\}$
  \item $\mathcal{N}^n(G)=\Span\{e^{i_0\cdots i_n}\in\mathcal{R}^n(V): i_0\cdots i_p \text{ is not allowed.}\}$
\end{itemize}
such that
$$\mathcal{R}^n(V)=\mathcal{A}^n(G)\oplus\mathcal{N}^n(G).$$
Similarly, $\mathcal{A}^*(G)$ is not necessary preserved by the differential operator $d$. To obtain a cohomology theory of a digraph, set $\mathcal{J}^n(G)=\mathcal{N}^n(G)+d\mathcal{N}^{n-1}(G)\subset\mathcal{R}^n(V)$ and
$$\Omega^n(G)=\mathcal{R}^n(V)/\mathcal{J}^n(G)=\mathcal{A}^n(G)/(\mathcal{J}^n(G)\cap\mathcal{A}^n(G)).$$

\begin{proposition}[\cite{GLMY2}] The pair $(\Omega^*(G),d)$ forms a well defined quotient complex of $(\mathcal{R}^*(V),d)$, which is the dual complex of $(\Omega_*(G),\partial)$. The corresponding cohomology is denoted by $H^*(G)$, and it is the dual space of $H_*(G)$.
\end{proposition}

\subsection{Cup product}\label{cup}

The cup product $\cup:\Omega^p(G)\otimes\Omega^q(G)\rightarrow\Omega^{p+q}(G)$ is induced by the concatenation on the space $\Lambda^*(V)$:
$$e^{i_0\cdots i_p}e^{j_0\cdots j_q}=\delta_{i_pj_0}e^{i_0\cdots i_p j_1\cdots j_q}.$$

One can also formulate the cup product as follows.

%By definition, each $u\in \Omega_n(G)$ can be written as
For $u\in\Omega_n(G)$, we write it more explicitly as
%$$u=\sum_{I,|I|=n}c_{\mathbf{n}}e_{\mathbf{n}},\quad e_{\mathbf{n}}\in \mathcal{A}_n(G).$$
$$u=\sum_{I,|I|=n}c_{I}e_{I},\quad e_I=e_{i_0i_1\ldots i_n}\in \mathcal{A}_n(G).$$

And for each s-regular allowed $n$-path $e_{I}=e_{i_0i_1\ldots i_n}$, denote
$$e_I\big|_{0\ldots p}=e_{i_0\ldots i_p},\quad e_{I}\big|_{p\ldots n}=e_{i_p\ldots i_n}.$$

For $\alpha\in\Omega^p(G)$, $\beta\in\Omega^q(G)$ , we define $\alpha\cup\beta$ on $\Omega^{p+q}(G)$ as
\begin{equation}\label{cupprod}
(\alpha\cup\beta, u)=\sum_{I,|I|=p+q}c_{I}(\alpha, e_I\big|_{0\ldots p})(\beta,e_I\big|_{p\ldots p+q}),
\end{equation}
where $u=\sum_{I,|I|=p+q}c_{I}e_{I}\in\Omega_{p+q}(G)$.
\begin{proposition} (1) $\alpha\cup\beta$ in \eqref{cupprod} is well defined for $\alpha,\beta\in\Omega^*(G)$. We call it the cup product of $\alpha$ and $\beta$.

(2) The cup product can descend to $H^*(G)$. That is, for any $p,q\in\mathbb{N}$,
$$H^p(G)\otimes H^q(G)\stackrel{\cup}{\longrightarrow} H^{p+q}(G).$$

\end{proposition}

\begin{proof} (1) It suffices to check that for any $n^p\in \mathcal{N}^p(G)$, $n^{p-1}\in \mathcal{N}^{p-1}(G)$, the following identity holds
$$\sum_{I,|I|=p+q}c_{I}(n^p+dn^{p-1}, e_{I}\big|_{0\ldots p})(\beta,e_{I}\big|_{p\ldots p+q})=0.$$
Note that  $e_{I}\big|_{0\ldots p}$ is still allowed path in $G$, then
$$(n^p, e_{I}\big|_{0\ldots p})=0,\quad n^p\in \mathcal{N}^p(G).$$
then
\begin{align*}
  &\sum_{I,|I|=p+q}c_{I}(n^p+dn^{p-1}, e_{I}\big|_{0\ldots p})(\beta,e_{I}\big|_{p\ldots p+q})\\
=~&\sum_{I,|I|=p+q}c_{I}(n^{p-1}, \partial e_{I}\big|_{0\ldots p})(\beta,e_{I}\big|_{p\ldots p+q}).
\end{align*}

\begin{itemize}
  \item If $\delta_i e_{I}\big|_{0\ldots p}=(-1)^ie_{I}\big|_{0\ldots\hat{i}\ldots p}\in\mathcal{A}_{p-1}(G)$, the corresponding terms in the above identity vanish.
  \item If $\delta_i e_{I}\big|_{0\ldots p}=(-1)^ie_{I}\big|_{0\ldots\hat{i}\ldots p}$ is not allowed, it means
$\delta_ie_{I}$ is not allowed. Since
$$u=\sum_{I,|I|=p+q}c_{I}e_{I}\in\Omega_{p+q}(G),$$
it means that there exists another summand in $u$, and its corresponding term cancels the non-allowed path.  So the corresponding summation will be zero.
\end{itemize}

(2)It is obvious by definition.\qedhere
\end{proof}

\begin{remark} One can check that the cup product defined by \eqref{cupprod} coincides with concatenation on $\Omega^*(G)$. We leave it as an exercise for the readers.
\end{remark}

To study the skew-symmetry property of the cup product, we express the cup product in terms of the star product and diagonal map. Now we recall the necessary definitions and basic properties, and refer to \cite{GLMY2} for more details.

\begin{definition}[cross product and star product, \cite{Grigoryan}] Let $X$ and $Y$ be two digraphs.
\begin{enumerate}
  \item For any s-regular elementary $p$-path $e_x$ in $X$ and $q$-path $e_y$ in $Y$, the cross product $e_x\times e_y$ is defined as a $(p+q)$-path in $X\boxdot Y$ by
$$e_x\times e_y=\sum_{z\in\Sigma_{x,y}}(-1)^{L(z)}e_z,$$
where $\Sigma_{x,y}$  is the set of all stair-like paths $z$ on $X\boxdot Y$ whose projections on $X$ and $Y$ are respectively $x$ and $y$, $L(z)$ is the number of cells in $\mathbb{N}_+^2$ below the staircase $S(z)$. See the following figure\footnote{We refer to the figure in \cite{GMY}.} for an explanation.
\begin{figure}[H]
  \centering
  % Requires \usepackage{graphicx}
  \includegraphics[width=14cm]{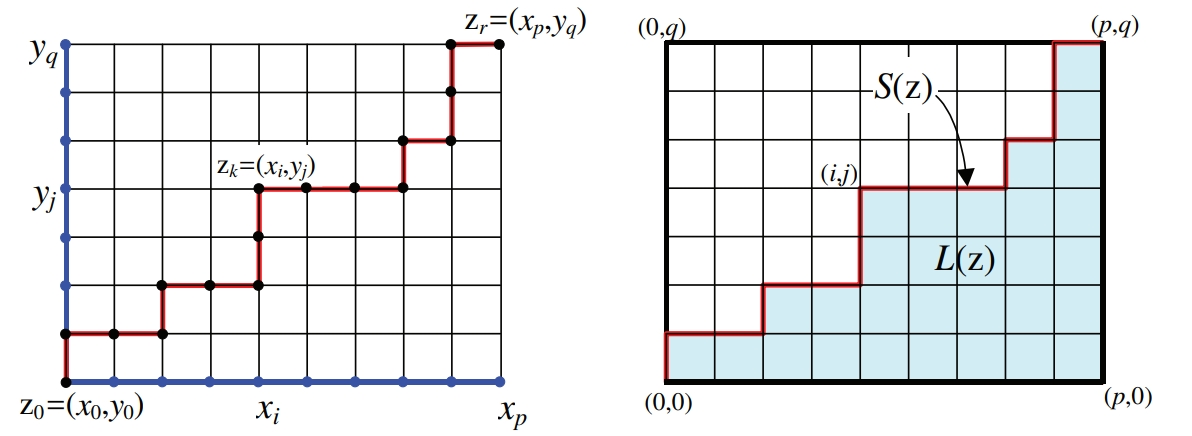}\\
  %\caption{}\label{}
\end{figure}

And it extends by linearity to all $u\in\mathcal{R}_p(X)$ and $v\in\mathcal{R}_q(Y)$ so that $u\times v\in\mathcal{R}_{p+q}(X\boxdot Y)$.
  \item For any s-regular elementary $p$-form $e^{i_0i_1\ldots i_p}$ on $X$ and $q$-form $e^{j_0j_1\ldots j_q}$ on $Y$, the star product $e^{i_0i_1\ldots i_p}\star e^{j_0j_1\ldots j_q}$ is defined as a $(p+q)$-form on $Z=X\boxdot Y$ by
$$e^{i_0i_1\ldots i_p}\star e^{j_0j_1\ldots j_q}=e^{(i_0j_0)(i_1j_0)\ldots(i_pj_0)(i_pj_1)\ldots(i_pj_q)}.$$
And it extends by linearity to all $\alpha\in\mathcal{R}^p(X)$ and $\beta\in\mathcal{R}^q(Y)$ so that $\alpha\star\beta\in\mathcal{R}^{p+q}(X\boxdot Y)$.
\end{enumerate}
\end{definition}

%We will not use the explicit expression of the cross product, but the following simple examples.
\begin{example}\label{crossexample} Let $a,b\in V(X)$, $1,2\in V(Y)$, we have
\begin{enumerate}
  \item $e_a\times e_{12}=e_{(a1)(a2)}$, $e_{ab}\times e_1=e_{(a1)(b1)}$;
  \item $e_{ab}\times e_{12}=e_{(a1)(b1)(b2)}-e_{(a1)(a2)(b2)}$.
\end{enumerate}
\end{example}

We will not use the explicit expression of the cross product, but the following important properties.

\begin{lemma}[\cite{Grigoryan}]\label{crossstarprop} (1) Let $\alpha\in\mathcal{R}^p(X)$, $\beta\in\mathcal{R}^q(Y)$, $u\in\mathcal{R}_p(X)$ and $v\in\mathcal{R}_q(Y)$, then
$$(\alpha\star\beta,u\times v)=(\alpha,u)(\beta,v).$$

(2) If $a\in \mathcal{A}_p(X)$, $b\in\mathcal{A}_q(Y)$, then $a\times b\in\mathcal{A}_{p+q}(X\boxdot Y)$. Moreover, if $u\in\Omega_p(X)$, and $v\in\Omega_q(Y)$, then $u\times v\in\Omega_{p+q}(X\boxdot Y)$. In particular,
$$\partial(u\times v)=(\partial u)\times v+(-1)^pu\times(\partial v).$$

(3) The star product is well defined for all $\alpha\in\Omega^p(X)$ and $\beta\in\Omega^q(Y)$, and $\alpha\star\beta\in\Omega^{p+q}(X\boxdot Y)$. In particular,
$$d(\alpha\star\beta)=(d\alpha)\star\beta+(-1)^p\alpha\star(d\beta).$$
\end{lemma}

Now for any $p,q\geq0$, any $\alpha\in\Omega^p(G)$, $\beta\in\Omega^q(G)$, $u=\sum_{I,|I|=p+q}c_Ie_I\in\Omega_{p+q}(G)$, we have
\begin{equation}\label{cupstar}
\begin{aligned}
(\alpha\cup\beta, u)=&\sum_{I,|I|=p+q}c_{I}(\alpha, e_{I}\big|_{0\ldots p})(\beta,e_{I}\big|_{p\ldots p+q})\\
                    =&\sum_{I,|I|=p+q}c_{I}(\alpha\star\beta, e_{I}\big|_{0\ldots p}\times e_{I}\big|_{p\ldots p+q})\\
                    =&~(\alpha\star\beta, \sum_{I,|I|=p+q}c_{I}e_{I}\big|_{0\ldots p}\times e_{I}\big|_{p\ldots p+q}).
\end{aligned}
\end{equation}

\begin{lem/def}\label{diag} Let $G=(V,E)$ be a digraph, we define the following two linear operators.
\begin{itemize}
  \item $\widetilde{\Delta}_{\sharp}: \mathcal{R}_n(G)\rightarrow\mathcal{R}_n(G\boxdot G)$: for s-regular elementary $n$-path $e_I$, set
$$\widetilde{\Delta}_{\sharp}e_I=\sum_{i=0}^{|I|}e_I\big|_{0\ldots i}\times e_I\big|_{i\ldots n}.$$
  \item the transposition operator $t_{\sharp}$:
  $$t_{\sharp}(e_{i_0i_1\ldots i_p}\times e_{j_0j_1\ldots j_q})=(-1)^{pq}e_{j_0j_1\ldots j_q}\times e_{i_0i_1\ldots i_p}.$$
\end{itemize}
\begin{enumerate}
  \item If $u=\sum_{I,|I|=n}c_{I}e_{I}\in\Omega_n(G)$, then
\begin{align*}
&\widetilde{\Delta}_{\sharp}(u)=\sum_{i=0}^{n}\sum_{I,|I|=n}c_{I}e_{I}\big|_{0\ldots i}\times e_{I}\big|_{i\ldots n}\in\Omega_n(G\boxdot G).\\
&t_{\sharp}\circ\widetilde{\Delta}_{\sharp}(u)=\sum_{i=0}^{n}(-1)^{i(n-i)}\sum_{I,|I|=n}c_{I}e_{I}\big|_{i\ldots n}\times e_I\big|_{0\ldots i}\in\Omega_n(G\boxdot G).
\end{align*}
  \item Both $\widetilde{\Delta}_{\sharp}$ and $t_{\sharp}\circ\widetilde{\Delta}_{\sharp}$ commute with the boundary operator $\partial$.
\end{enumerate}
We call $\widetilde{\Delta}_{\sharp}$ and $t_{\sharp}\circ\widetilde{\Delta}_{\sharp}$ the diagonal approximations.
\end{lem/def}

\begin{proof} (1) By Lemma \ref{crossstarprop} (2), we have
\begin{align*}
\partial\big(\widetilde{\Delta}_{\sharp}(u)\big)=~&\sum_{i=0}^{n}\left(\sum_{I,|I|=n}c_{I}\partial e_{I}\big|_{0\ldots i}\right)\times e_{I}\big|_{i\ldots n}\\
                                &+\sum_{i=0}^n(-1)^i\sum_{I,|I|=n}c_Ie_{I}\big|_{0\ldots i}\times \left(\partial e_{I}\big|_{i\ldots n}\right)
\end{align*}
Since $u=\sum_Ic_Ie_I\in\Omega_n(G)$, by Lemma \ref{deltaiP}, we have
$$\delta_ju=\sum_Ic_I\delta_je_I\in\mathcal{A}_{n-1}(G),\quad j=0,1,\ldots n.$$
It implies that, with a shift of index,
\begin{align*}
&\sum_Ic_I\delta_j\left(e_I\big|_{0,\ldots,i}\right)\in\mathcal{A}_{i-1}(G),\quad j=0,1,\ldots,i;\\
&\sum_Ic_I\delta_j\left(e_I\big|_{i,\ldots,n}\right)\in\mathcal{A}_{n-i-1}(G),\quad j=0,\ldots,n-i.
\end{align*}
Thus, we have $\partial\big(\widetilde{\Delta}_{\sharp}(u)\big)\in\mathcal{A}_{n-1}(G\boxdot G)$, i.e., $\widetilde{\Delta}_{\sharp}(u)\in\Omega_n(G\boxdot G)$. The similar argument implies that $t_{\sharp}\circ\widetilde{\Delta}_{\sharp}(u)\in\Omega_n(G\boxdot G)$.

(2) By definition, one can check directly, for $u\in\Omega_n(G)$,
$$\widetilde{\Delta}_{\sharp}(\partial u)=\partial \big(\widetilde{\Delta}_{\sharp}(u)\big),\quad t_{\sharp}\circ\widetilde{\Delta}_{\sharp}(\partial u)=\partial \big(t_{\sharp}\circ\widetilde{\Delta}_{\sharp}(u)\big).$$
\end{proof}

In particular, we learn from the proof that for each fixed $i$,
$$\sum_{I,|I|=n}c_{I}e_{I}\big|_{0\ldots i}\times e_{I}\big|_{i\ldots n}\in\Omega_n(G\boxdot G).$$

Dually we have two linear operators $\widetilde{\Delta}^{\sharp}:\Omega^n(G\boxdot G)\rightarrow \Omega^n(G)$, and $t^{\sharp}$. By \eqref{cupstar}, for $\alpha\in\Omega^p(G)$, $\beta\in\Omega^q(G)$,
$$\alpha\cup\beta=\widetilde{\Delta}^{\sharp}(\alpha\star\beta),\quad \beta\star\alpha=(-1)^{pq}t^{\sharp}(\alpha\star\beta).$$
which induce, for $\varphi\in H^p(G)$, $\psi\in H^q(G)$,
$$\varphi\cup\psi=\widetilde{\Delta}^*(\varphi\star\psi),\quad (-1)^{pq}\psi\cup\varphi=t^*\widetilde{\Delta}^*(\varphi\star\psi).$$

\begin{theorem}\label{skewsymm} The two chain maps $\widetilde{\Delta}_{\sharp}$ and $t_{\sharp}\circ\widetilde{\Delta}_{\sharp}$ are chain homotopic. Furthermore, for $\varphi\in H^p(G)$, $\psi\in H^q(G)$, we have
$$\varphi\cup\psi=(-1)^{pq}\psi\cup\varphi.$$
\end{theorem}

\begin{proof}
We will construct the chain homotopy $F: \Omega_*(G)\rightarrow \Omega_{*+1}(G\boxdot G)$ one by one. For simplicity, assume that $G$ is connected.

First, for the case $k=0$, $\Omega_0(G)=\Span\{e_i\}_{i\in V(G)}$, we have
$$t_{\sharp}\circ\widetilde{\Delta}_{\sharp}(e_i)-\widetilde{\Delta}_{\sharp}(e_i)=e_i\times e_i-e_i\times e_i=0.$$
Thus, $F(e_i)$ is a cycle in $\Omega_1(G\boxdot G)$. For simplicity, we fix $e_{ab}\in \Omega_1(G)$, and set
$$F(e_i)=e_a\times e_{ab}-e_b\times e_{ab}+e_{ab}\times e_b-e_{ab}\times e_a\quad \text{for any }e_i\in\Omega_0(G).$$
Clearly, $\partial F(e_i)=0$ and $t_{\sharp}\circ\widetilde{\Delta}_{\sharp}(e_i)-\widetilde{\Delta}_{\sharp}(e_i)=\partial F(e_i)+F(\partial e_i)$ holds naturally.

By our choice of $F$ on $\Omega_0(G)$, furthermore, we can obtain
\begin{align*}
&t_{\sharp}\circ\widetilde{\Delta}_{\sharp}(e_{ij})-\widetilde{\Delta}_{\sharp}(e_{ij})-F(\partial e_{ij})\\
=&~e_{ij}\times e_i+e_j\times e_{ij}-e_i\times e_{ij}-e_{ij}\times e_j\\\
=&~\partial (e_{ij}\times e_{ij}).
\end{align*}
So we can define $F(e_{ij})=e_{ij}\times e_{ij}\in\Omega_2(G\boxdot G).$

Now to construct $F:\Omega_n(G)\rightarrow\Omega_{n+1}(G\boxdot G)$, first, let us choose an integral basis of $\Omega_n(G)$ consisting of minimal $n$-paths $\{u\}$, then we only need to define $F(u)$ for this basis $\{u\}$. Note that
$$t_{\sharp}\circ\widetilde{\Delta}_{\sharp}(u)-\widetilde{\Delta}_{\sharp}(u)-F(\partial u)$$
is a cycle in $\Omega_n\bigg(\Supp(u)\boxdot\Supp(u)\bigg)\subset\Omega_n(G\boxdot G)$: by Lemma \ref{diag},
\begin{align*}
&\partial \left(t_{\sharp}\circ\widetilde{\Delta}_{\sharp}(u)\right)-\partial\widetilde{\Delta}_{\sharp}(u)-\partial F(\partial u)\\
=~&t_{\sharp}\circ\widetilde{\Delta}_{\sharp}(\partial u)-\widetilde{\Delta}_{\sharp}(\partial u)-\left(-F(\partial^2 u)+t_{\sharp}\circ\widetilde{\Delta}_{\sharp}(\partial u)-\widetilde{\Delta}_{\sharp}(\partial u)\right)=0.
\end{align*}

Since $H_{i>0}\bigg(\Supp(u)\bigg)=0$, then by K\"{u}nneth formula\footnote{One can also use the acyclic result to prove the K\"{u}nneth formula, see the last subsection.} in \cite{GMY},
$$H_{i>0}\bigg(\Supp(u)\boxdot\Supp(u)\bigg)=0.$$
Then there exists $v\in \Omega_{n+1}\bigg(\Supp(u)\boxdot\Supp(u)\bigg)\subset\Omega_{n+1}(G\boxdot G)$ such that
$$\partial v=t_{\sharp}\circ\widetilde{\Delta}_{\sharp}(u)-\widetilde{\Delta}_{\sharp}(u)-F(\partial u),$$
then we can define
$$F(u)=v.$$
By $\mathbb{K}$-linear extension, for any $p_n\in\Omega_n(G)$,
$$t_{\sharp}\circ\widetilde{\Delta}_{\sharp}(p_n)-\widetilde{\Delta}_{\sharp}(p_n)-F(\partial p_n)=\partial F(p_n).$$
After descending to homology, we have
$$t_*\circ\widetilde{\Delta}_{*}(q_n)=\widetilde{\Delta}_{*}(q_n)\quad \text{for } q_n\in H_n(G).$$
Furthermore, for $\varphi\in H^p(G),\psi\in H^q(G)$, with $p+q=n$ we have
\begin{align*}
(\psi\cup\varphi, q_n)&=(\widetilde{\Delta}^*(\psi\star\varphi),q_n)=(\psi\star\varphi, \widetilde{\Delta}_*(q_n))\\
                      &=((-1)^{pq}t^*(\varphi\star\psi),\widetilde{\Delta}_*(q_n))\\
                      &=(-1)^{pq}(\varphi\star\psi,t_*\widetilde{\Delta}_*(q_n))\\
                      &=(-1)^{pq}(\varphi\star\psi,\widetilde{\Delta}_*(q_n))\\
                      &=(-1)^{pq}(\widetilde{\Delta}^*(\varphi\star\psi),q_n)=(-1)^{pq}(\varphi\cup\psi,q_n).\qedhere
\end{align*}

\end{proof}

\subsection{Other applications}\label{other}

The acyclic model of a homology theory has a wide application, such as K\"{u}nneth formula, the universal coefficient theorem and so on. In \cite{GMY}, they show the K\"{u}nneth formula by building the isomorphism between
$$\Omega_*(X\boxdot Y)\quad \text{and}\quad \Omega_*(X)\otimes\Omega_*(Y).$$
We can also consider the strong product $X\boxtimes Y$, which is a digraph defined as follows:
\begin{itemize}
  \item $V(X\boxtimes Y)=V(X)\times V(Y)$;
  \item the arrows are defined as follows:
  \begin{itemize}
    \item horizontal edge: $(x,y)\rightarrow (x',y)$, where $x\rightarrow x'$,
    \item vertical edge: $(x,y)\rightarrow (x,y')$, where $y\rightarrow y'$,
    \item diagonal edge: $(x,y)\rightarrow (x',y')$, where $x\rightarrow x'$ and $y\rightarrow y'$.
  \end{itemize}
\end{itemize}
In this case, the chains $\Omega_*(X\boxtimes Y)$ and $\Omega_*(X)\otimes\Omega_*(Y)$ are not isomorphic, but there exist a chain map and its homotopic inverse, where the existence of the chain homotopy can be shown by our acyclic result.

\vskip 0.2cm

\begin{appendices}

\section{An example: $H_*^{\sr}(G)\ncong H_*^{\rr}(G)$}\label{A}

In this section, we will give an example to see that $H_*^{\sr}(G)\ncong H_*^{\rr}(G)$, where we re-denote
\begin{itemize}
  \item $H_*^{\sr}(G)$: the path homologies under the s-regular condition;
  \item $H_*^{\rr}(G)$: the path homologies under the regular condition in Remark \ref{regular} defined by \cite{GLMY3}.
\end{itemize}

As we write in Remark \ref{regular}, for any simple finite digraph $G$, $H_*^{\sr}(G)$ is bounded above. We will give an example $G$, satisfying
$$\dim H_*^{\rr}(G)=\infty.$$
Thus the two kinds of path homologies are not always the same.

This example\footnote{The authors thank A. Grigoryan and I. Sergei for sharing such an interesting example.} is due to Gabor Lippner and Paul Horn, 2012.
\begin{figure}[H]
	\centering
	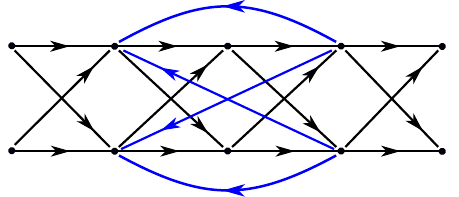
	%\caption[]{$\Supp(e_{0123})$}
	%\label{Fig:Q1}
\end{figure}
One can understand this example as the generalized join\footnote{That is, for each vertex $i$ in $H$, we associate to a digraph $G_i$, then $X_G=(V(X_G),E(X_G))$ is defined to be
\begin{itemize}
  \item $V(X_H)=\cup_{i\in V(H)}V(G_i)$;
  \item $E(X_H)$ consists of two kinds of edges:
  \begin{itemize}
    \item the edges in $G_i$,
    \item the edges among $\{G_i\}$: for any $a\in V(G_i)$, $b\in V(G_j)$, if $i\rightarrow j$ in $H$, then $a\rightarrow b$.
  \end{itemize}
\end{itemize}}
(which is also known as graph extension) $G:=X_H$ of the digraphs $G_i=\{i,-i\}_{i=1,2,\ldots,5}$ along the following digraph $H$.
\begin{figure}[H]
	\centering
	%% Creator: Inkscape 1.0.1 (3bc2e813f5, 2020-09-07), www.inkscape.org
%% PDF/EPS/PS + LaTeX output extension by Johan Engelen, 2010
%% Accompanies image file 'append2.pdf' (pdf, eps, ps)
%%
%% To include the image in your LaTeX document, write
%%   \input{<filename>.pdf_tex}
%%  instead of
%%   \includegraphics{<filename>.pdf}
%% To scale the image, write
%%   \def\svgwidth{<desired width>}
%%   \input{<filename>.pdf_tex}
%%  instead of
%%   \includegraphics[width=<desired width>]{<filename>.pdf}
%%
%% Images with a different path to the parent latex file can
%% be accessed with the `import' package (which may need to be
%% installed) using
%%   \usepackage{import}
%% in the preamble, and then including the image with
%%   \import{<path to file>}{<filename>.pdf_tex}
%% Alternatively, one can specify
%%   \graphicspath{{<path to file>/}}
%% 
%% For more information, please see info/svg-inkscape on CTAN:
%%   http://tug.ctan.org/tex-archive/info/svg-inkscape
%%
\begingroup%
  \makeatletter%
  \providecommand\color[2][]{%
    \errmessage{(Inkscape) Color is used for the text in Inkscape, but the package 'color.sty' is not loaded}%
    \renewcommand\color[2][]{}%
  }%
  \providecommand\transparent[1]{%
    \errmessage{(Inkscape) Transparency is used (non-zero) for the text in Inkscape, but the package 'transparent.sty' is not loaded}%
    \renewcommand\transparent[1]{}%
  }%
  \providecommand\rotatebox[2]{#2}%
  \newcommand*\fsize{\dimexpr\f@size pt\relax}%
  \newcommand*\lineheight[1]{\fontsize{\fsize}{#1\fsize}\selectfont}%
  \ifx\svgwidth\undefined%
    \setlength{\unitlength}{245.38610503bp}%
    \ifx\svgscale\undefined%
      \relax%
    \else%
      \setlength{\unitlength}{\unitlength * \real{\svgscale}}%
    \fi%
  \else%
    \setlength{\unitlength}{\svgwidth}%
  \fi%
  \global\let\svgwidth\undefined%
  \global\let\svgscale\undefined%
  \makeatother%
  \begin{picture}(1,0.10355591)%
    \lineheight{1}%
    \setlength\tabcolsep{0pt}%
    \put(0,0){\includegraphics[width=\unitlength,page=1]{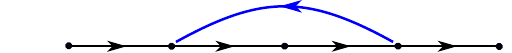}}%
    \put(0.10909691,0.04326378){\makebox(0,0)[lt]{\lineheight{1.25}\smash{\begin{tabular}[t]{l}$1$\end{tabular}}}}%
    \put(0.3059301,0.04134048){\makebox(0,0)[lt]{\lineheight{1.25}\smash{\begin{tabular}[t]{l}$2$\end{tabular}}}}%
    \put(0.52989797,0.0402189){\makebox(0,0)[lt]{\lineheight{1.25}\smash{\begin{tabular}[t]{l}$3$\end{tabular}}}}%
    \put(0.75695554,0.04240877){\makebox(0,0)[lt]{\lineheight{1.25}\smash{\begin{tabular}[t]{l}$4$\end{tabular}}}}%
    \put(0.94833692,0.04224259){\makebox(0,0)[lt]{\lineheight{1.25}\smash{\begin{tabular}[t]{l}$5$\end{tabular}}}}%
    \put(-0.00185056,0.00706349){\makebox(0,0)[lt]{\lineheight{1.25}\smash{\begin{tabular}[t]{l}$H=$\end{tabular}}}}%
  \end{picture}%
\endgroup%

	%\caption[]{$\Supp(e_{0123})$}
	%\label{Fig:Q1}
\end{figure}

Applying the computation method in \cite{GTY}, one can see that $\dim H_*^{\rr}(X_H)=\infty$. In particular, for any $N\in\mathbb{N}$,
$$u=(e_1-e_{-1})\bullet_j[(e_2-e_{-2})\bullet_j(e_3-e_{-3})\bullet_j(e_4-e_{-4})]^{\bullet_j N}\bullet_j(e_5-e_{-5})$$
is a non-trivial $(3N+1)$-cycle, where $\bullet_j$ is the induced operation from the bilinear operation
\begin{align*}
\bullet_j : &\Lambda_p(V_1)\times \Lambda_q(V_2)\rightarrow \Lambda_{p+q+1}(V_1\cup V_2)\\
&e_{i_0i_1\ldots i_p}\bullet_j e_{j_0j_1\ldots j_q}=e_{i_0\ldots i_pj_0\ldots j_q}.
\end{align*}
We refer to \cite{GMY} for more details.

\section{A minimal $4$-path and some remarks}\label{B}

In this section, we will explain our proof of acyclic model theorem \ref{acyclicresult} again via a minimal $4$-path.

Let us consider the following path
\begin{align*}
P=~&e_{S159E}-e_{S169E}+e_{S269E}\\
   &+e_{S16(10)E}-e_{S26(10)E}+e_{S27(10)E}-e_{S37(10)E}\\
   &-e_{S27(11)E}+e_{S37(11)E}-e_{S38(11)E}+e_{S48(11)E}.
\end{align*}
It could be a minimal $4$-path with the following supporting digraph
\begin{figure}[H]
	\centering
	\input{length4.pdf_tex}
	%\caption[]{$\Supp(e_{0123})$}
	%\label{Fig:Q1}
\end{figure}
Now we add the new directed edges $6\rightarrow E$ which make the paths
$$P_1=e_{S159E}-e_{S169E}+e_{S269E}$$
and
$$P_2=e_{S16(10)E}-e_{S26(10)E}+e_{S27(10)E}-e_{S37(10)E}-e_{S27(11)E}+e_{S37(11)E}-e_{S38(11)E}+e_{S48(11)E}$$
$\partial$-invariant in the new larger digraph $\widehat{\Supp(P)}$.
\begin{figure}[H]
	\centering
	\input{length4new.pdf_tex}
	%\caption[]{$\Supp(e_{0123})$}
	%\label{Fig:Q1}
\end{figure}

In particular, both $P_1$ and $P_2$ are minimal as Claim \ref{betterchoice} says. Denote $X_i=\Supp(P_i)$, $i=1,2$. Then $(X_1, X_2)$ form a Mayer-Vietoris pair for $\widehat{\Supp(P)}$.
\begin{figure}[H]
	\centering
	\input{length4step1.pdf_tex}
	%\caption[]{$\Supp(e_{0123})$}
	%\label{Fig:Q1}
\end{figure}

Let us redo the computation similar to the exotic cube. First, by applying the Mayer-Vietoris exact sequence for the pair $(X_1,X_2)$ and hypothesis, we have
$$\widetilde{H}_*(\widehat{\Supp(P)})\cong\widetilde{H}_*(Z)\cong 0.$$
Second, one can show that the embedding
$$(\Omega_*(\Supp(P)),\partial) \hookrightarrow(\Omega_*(\widehat{\Supp(P)}),\partial)$$
induces the isomorphism on path homologies. Then we are done.

\begin{remark}
Actually, one can further add new edge $7\rightarrow E$ to embed $X_2$ into a larger graph $\widehat{X_2}$ such that the paths
\begin{align*}
&P_{21}=e_{S16(10)E}-e_{S26(10)E}+e_{S27(10)E}-e_{S37(10)E},\\
&P_{22}=-e_{S27(11)E}+e_{S37(11)E}-e_{S28(11)E}+e_{S48(11)E}
\end{align*}
are $\partial$-invariant in $\widehat{X_2}$. Similarly, it is obvious that $(\Supp(P_{21}),\Supp(P_{22}))$ form a Mayer-Vietoris pair for $\widehat{X_2}$. By the same reason, we have
$$\widetilde{H}_*(\widehat{X_2})=0.$$
Then one can show that in this case the two embedding maps
$$(\Omega_*(X_2),\partial) \hookrightarrow(\Omega_*(\widehat{X_2}),\partial)$$ induce the isomorphisms on homologies.
\end{remark}

\begin{remark}[Related to Claim \ref{betterchoice}]
Alternatively, one can first add two edges $6\rightarrow E$ and $7\rightarrow E$ to get $\widehat{\widehat{\Supp(P)}}$, which makes the paths
$$P_1'=e_{S16(10)E}-e_{S26(10)E}+e_{S27(10)E}-e_{S37(10)E},\quad\text{and}\quad P_2'=P-P_1$$
$\partial$-invariant in $\widehat{\widehat{\Supp(P)}}$. One can also consider the other pair $(X_1',X_2')$ related to $P_1'$, $P_2'$ as follows:
\begin{figure}[H]
	\centering
	\input{length4step2.pdf_tex}
	%\caption[]{$\Supp(e_{0123})$}
	%\label{Fig:Q1}
\end{figure}
However, in this case, $P_2'$ is not a minimal path but the summand of the two minimal paths, then we can not directly apply our induction assumption to say that $\widetilde{H}_*(X_2')=0$. And in this case, $Z$ also becomes more complicated than the previous case. Thus, we will remove the redundant new edges $6\rightarrow E$ or $7\rightarrow E$ as we prove in Claim \ref{betterchoice}.
\end{remark}

\end{appendices}

\end{document}